\documentclass[11pt,a4paper]{article}
\usepackage[T1]{fontenc}
\usepackage[margin=1in]{geometry}
\usepackage{microtype}
\usepackage{needspace}
\usepackage[numbers]{natbib}
\usepackage{hyperref}

\usepackage{amsmath,amssymb,amsthm,mathtools,booktabs,enumitem}
\usepackage{xcolor}
\usepackage{placeins}
\usepackage{graphicx}
\usepackage{upquote}
\usepackage{aliascnt}
\usepackage{longtable}
\usepackage{array}
\usepackage{tikz}
\hypersetup{hidelinks,pdftitle={Euler discretizations, invariant faces and brackets for lower exponents of positive switched systems},pdfauthor={Francesco Paolo Maiale}}
\usepackage[capitalise,noabbrev]{cleveref}
\usepackage{bookmark}
\bookmarksetup{depth=2}
\makeatletter\g@addto@macro\UrlBreaks{\do\_}\makeatother

\newcommand{\sharedthm}[3]{\newaliascnt{#1}{theorem}\newtheorem{#1}[#1]{#2}\aliascntresetthe{#1}\crefname{#1}{#2}{#3}}
\theoremstyle{plain}
\newtheorem{theorem}{Theorem}[section]
\sharedthm{proposition}{Proposition}{Propositions}
\sharedthm{lemma}{Lemma}{Lemmas}
\sharedthm{corollary}{Corollary}{Corollaries}
\sharedthm{conjecture}{Conjecture}{Conjectures}
\theoremstyle{definition}
\sharedthm{definition}{Definition}{Definitions}
\sharedthm{example}{Example}{Examples}
\sharedthm{observation}{Numerical observation}{Numerical observations}
\sharedthm{computation}{Computation}{Computations}
\theoremstyle{remark}
\sharedthm{remark}{Remark}{Remarks}

\newcommand{\cA}{\mathcal A}
\newcommand{\cB}{\mathcal B}
\newcommand{\cR}{\mathcal R}
\newcommand{\R}{\mathbb R}
\newcommand{\one}{\mathbf 1}
\newcommand{\lsr}{\check\rho}
\newcommand{\lexp}{\check\lambda}

\AddToHook{env/proposition/before}{\Needspace{6\baselineskip}}

\begin{document}

\title{Euler discretizations, invariant faces and brackets for lower exponents of positive switched systems}
\author{Francesco Paolo Maiale\\[3pt]
\small Gran Sasso Science Institute\\
\small Viale F. Crispi 7, 67100 L'Aquila, Italy\\
\small\href{mailto:francescopaolo.maiale@gssi.it}{\texttt{francescopaolo.maiale@gssi.it}}}
\date{}
\maketitle
\begin{abstract}
In this paper, we prove that explicit and implicit Euler discretizations yield monotone brackets for the lower Lyapunov exponent of positive switched systems. We establish convergence of the implicit bound for every Metzler family. For the explicit bound, we derive convergence criteria in reducible systems from a face formula for the uniform lower exponent on an invariant cone. The formula takes the minimum of the rates on faces that meet the cone nontrivially and are preserved by at least one mode, using only their preserving modes. For individually irreducible modes, we prove $O(h)$ discretization errors, where $h$ is the step size, and show that this order is sharp in general. At interior two-mode chattering optima with an irreducible optimal mixture, we obtain one-sided first-order bounds. We also construct a three-mode example in which the explicit bound converges while finite polytope lower bounds remain capped strictly below the true exponent. We prove convergence of fixed-antinorm Gripenberg brackets under individual irreducibility as the step and tolerance decrease and the search depth increases appropriately. Finally, we obtain independently verified continuous-time brackets by exact rational checks, with lower bounds that remain valid at smaller steps.
\end{abstract}
\noindent\textbf{Keywords:} Positive switched systems, Lower Lyapunov exponent, Lower spectral radius, Euler discretization, Invariant faces, Polytope antinorms.
\raggedbottom

\section{Introduction}

Let $\cA=\{A_1,\dots,A_m\}$ be a family of Metzler matrices.
Then the lower exponent $\lexp(\cA)$ is the best exponential growth rate achievable by switching in $\dot x=A_{\sigma(t)}x$ from interior initial states (\cref{prop:basic}). The computation of rigorous bounds requires estimating discretization error and verifying the discrete bounds. For reducible systems, boundary dynamics can limit the rates established by uniform lower-growth functions even when the discretization converges to $\lexp$.

The discrete lower spectral radius admits polytope bounds~\cite{GLP} and Gripenberg-type bounds~\cite{GM}. For continuous-time systems, Guglielmi, Laglia and Protasov~\cite[Sec.~4]{GLP} obtain a $C\tau+\varepsilon$ gap using $\mathrm e^{\tau A}$ and suitable polytope antinorms for individually $K$-irreducible generators, where $\varepsilon$ is the polytope approximation error. 

Earlier methods use irreducibility or secondary-eigenvector conditions~\cite{GLP,GM}; boundary eigenvectors can obstruct polytope methods~\cite[Sec.~6.4]{GP13}.

\subsection{Main results and conventions}

We first summarize the main results obtained in this paper, first theoretical and then numerical ones.

\begin{enumerate}[label=(R\arabic*),leftmargin=*]
\item \emph{Face formula} (\cref{thm:faceformula}). The uniform lower exponent on an invariant cone equals the minimum of the exponents of the faces it touches that are preserved by at least one mode, each restricted to its preserving modes. This identifies the boundary rates above. Occupation measures~\cite{Vin93,Gai04,GR06} make the face decomposition explicit within the ergodic variational framework~\cite{CK96,CK00,Jen19}. Reduced graph classes give the one-matrix analogue~\cite[Thms.~5.10 and~6.8]{FS80}.

\item \emph{Monotone Euler brackets} (\cref{thm:explicit,thm:implicit}). Explicit and implicit Euler steps give $g(h)\le\lexp\le G(h)$ that improve monotonically as $h\downarrow0$. The implicit bound always converges. The explicit limit lies between the boundary exponent and $\lexp$; cone criteria based on the face formula give convergence in reducible systems, verified on five test pairs (\cref{cor:eulerface,obs:cone}). Compared with the exponential discretization of~\cite[Sec.~4]{GLP}, Euler steps add monotonicity in $h$, convergence of the implicit bound without irreducibility, and polytope lower bounds that remain valid at smaller steps. The explicit bound parallels classical stability under smaller Euler steps~\cite{Opo,MoP86,PJ}.

\item \emph{Rates} (\cref{thm:eulerrate,thm:rateconst}). For individually irreducible modes both errors are $O(h)$, with this order optimal in general. At interior two-mode chattering optima with an irreducible optimal mixture we obtain one-sided first-order coefficients. Related work treats polyhedral Lyapunov functions shared by linear systems and their Pad\'e discretizations~\cite{RCS}, and second-order bounds for the upper Lyapunov exponent under positive dwell-time constraints~\cite{MeP}.

\item \emph{Limits of polytope bounds} (\cref{lem:bounded,ex:bcv}). Lower-growth functions comparable to $\one^\top x$ cannot establish rates above the boundary exponent. In a three-mode example, we show that the explicit Euler bound converges while lower bounds obtained from finite polytopes remain capped strictly below $\lexp$.

\item \emph{Rigorous brackets} (\cref{prop:polytope,prop:upper,cor:gripenberg}). Building on polytope antinorm methods~\cite{GP13,GLP,GZ20,GM}, rational polytopes, vector chains and exponential enclosures yield continuous-time brackets verified independently, also for reducible families; exact verification also appears for the joint spectral radius~\cite{Mej22}. Fixed-antinorm Gripenberg brackets converge for individually irreducible modes as $h\downarrow0$, with tolerance $\delta=o(h)$ and sufficient search depth.
\end{enumerate}

The convergence results are summarized below; $\lexp_{\min}\le\lexp$ denotes the boundary exponent (\cref{subsec:coneexponents}).

\begin{center}\small
\begin{tabular}{@{}>{\raggedright\arraybackslash}p{0.13\textwidth}>{\raggedright\arraybackslash}p{0.44\textwidth}>{\raggedright\arraybackslash}p{0.35\textwidth}@{}}
\toprule
& general Metzler family & individually irreducible modes\\
\midrule
expl. $g(h)$ & $g(h)\le\lexp$, monotone; $\lexp_{\min}\le\lim g\le\lexp$; $\lim g=\lexp$ under the cone criteria; open in general (\cref{thm:explicit,cor:eulerface}) & $\lexp-Ch\le g(h)$ (\cref{thm:eulerrate}); order $h$ sharp (\cref{ex:sharp})\\
\addlinespace
impl. $G(h)$ & $\lexp\le G(h)$, monotone, $G(h)\to\lexp$; no rate proved (\cref{thm:implicit}) & $G(h)\le\lexp+Ch$ (\cref{thm:eulerrate}); order $h$ sharp\\
\bottomrule
\end{tabular}
\end{center}

Experiments up to dimension~$20$ distinguish discretization and search effects and compare with Algorithm~(L) of~\cite{GLP}. Numerical schedule rates can fall below numerically minimized constant-mixture rates. Note, however, that these comparisons do not prove strict gaps from all constant mixtures.

\paragraph{Organization of the paper}

\Cref{sec:prelim} develops the face formula, \cref{sec:euler} proves the Euler bounds, \cref{sec:gripenberg} treats rigorous bounds and numerics, and \cref{sec:open} discusses some open problems. The appendices give the test matrices (\cref{app:data}), Gripenberg analysis and detailed runs (\cref{app:gm}), first-order diagnostics (\cref{app:tables}), and the two technical proofs (\cref{app:rateproof,app:faceproof}). Experimental protocols are collected in \cref{app:protocols}, and the larger numerical samples in \cref{app:samples}.

\section{Lower exponents and invariant faces}\label{sec:prelim}

In this section, we define the lower exponent, construct the lower-growth functions used in \cref{sec:euler}, and prove our main result (the face formula).

\subsection{Switched systems and lower exponents}\label{subsec:lowerexponents}

Let $\cA=\{A_1,\dots,A_m\}\subset\R^{d\times d}$ be a family of Metzler matrices, i.e. with nonnegative off-diagonal entries, and let $K=\R^d_+$ be the non-negative cone. We denote the simplex associated to $K$ by $\Delta=\{x\in K \ :\ \one^\top x=1\}$.

\begin{definition}
A \emph{switching signal} $\sigma$ is any measurable function $[0,\infty)\to\{1,\dots,m\}$. 
\end{definition}

For $t\ge s\ge0$, let $\Phi_\sigma(t,s)$ be the transition matrix of $\dot x=A_{\sigma(t)}x$ from time $s$ to time $t$, and write $\Phi_\sigma(t)=\Phi_\sigma(t,0)$. Notice that every $\Phi_\sigma(t,s)$ is nonnegative. Set
\[
 \bar\mu(t)=\inf_\sigma\|\Phi_\sigma(t)\|_1,\qquad \lexp(\cA)=\lim_{t\to\infty}\frac{\log \bar\mu(t)}{t}.
\]
For a Metzler matrix $S$, let $\eta(S)$ denote its spectral abscissa. Let $S_\alpha=\sum_i\alpha_iA_i$ for $\alpha$ in the control simplex $\Delta_m=\{u\in\R^m_+ \ :\ \sum_i u_i=1\}$. Throughout, $m$ is the number of modes, and
\[
 K_A=\max_{1\le i\le m}\|A_i\|_1,\qquad c_A=\max_{\substack{1\le i\le m\\1\le j\le d}}|(A_i)_{jj}| .
\]

\begin{proposition}\label{prop:basic} The following properties hold: \mbox{}
\begin{enumerate}[label=\textup{(\roman*)}]
\item The limit $\lexp$ exists, equals $\inf_{t>0}t^{-1}\log \bar\mu(t)$ and does not depend on the norm. In addition, $|\lexp|\le K_A$.
\item For every $x_0>0$, $\lexp=\lim_t t^{-1}\log\inf_\sigma\|\Phi_\sigma(t)x_0\|_1$. Hence $\lexp$ is the best exponential rate achievable from any interior initial state, by open-loop or by state-feedback switching.
\item For every $\alpha\in\Delta_m$, $\lexp\le\eta(S_\alpha)$. For every rational $x>0$,
\[
 \eta(S_\alpha)\le\max_j\frac{(S_\alpha x)_j}{x_j} \qquad \textup{(Collatz--Wielandt inequality)}.
\]
This upper bound is rational if the matrices $A_i$ and the chosen weights $\alpha_i$ are rational.
\end{enumerate}
\end{proposition}

\begin{proof}
Concatenating switching signals shows $\bar\mu(t+s)\le \bar\mu(t)\bar\mu(s)$. Since $\|\Phi_\sigma(t)^{\pm1}\|_1\le e^{K_At}$, we have $e^{-K_At}\le \bar\mu(t)\le e^{K_At}$. Fekete's lemma gives the limit and the infimum formula, while equivalence of norms gives norm independence, proving (i).

For (ii), notice that given $\Phi\ge0$ and $x_0>0$, it holds
\[
 \min_j(x_0)_j\,\|\Phi\|_1\le\one^\top\Phi x_0\le d\max_j(x_0)_j\,\|\Phi\|_1 .
\]
A feedback law produces, for a fixed initial state, an open-loop signal. For (iii), the Lie--Trotter formula (see~\cite{Trotter}) gives
\[
\left(\prod_i e^{(t/n)\alpha_iA_i}\right)^n\xrightarrow{n\to\infty} e^{tS_\alpha},
\]
and each left-hand side is a transition matrix of a switching signal. Hence $\bar\mu(t)\le\|e^{tS_\alpha}\|_1$, and $t^{-1}\log\|e^{tS_\alpha}\|_1\to\eta(S_\alpha)$. The Collatz--Wielandt bound follows from the Perron--Frobenius theorem applied to $S_\alpha+c_AI\ge0$ in the weighted norm $\max_j|y_j|/x_j$.
\end{proof}

\subsection{Uniform exponents on cones}\label{subsec:coneexponents}

The exponent $\lexp$ describes optimal growth from each fixed interior state. To describe growth uniformly over normalized states in a cone, let $C\subseteq K$ be a closed convex cone, $C\ne\{0\}$, that every $e^{tA_i}$ maps into itself, and define
\[
 \mu_C(t)=\inf_\sigma\min_{x\in C\cap\Delta}\one^\top\Phi_\sigma(t)x,\qquad \lexp_C=\lim_{t\to\infty}\frac{\log\mu_C(t)}t ,
\]
where the infima range over both signals and initial states in $C\cap\Delta$. Every $\Phi_\sigma(t)$ maps $C$ into $C$: for piecewise-constant $\sigma$ directly, and for measurable $\sigma$ by approximation, since $C$ is closed and $\Phi_\sigma(t)$ depends continuously on $\sigma$ in $L^1$. So, for $x\in C$, we have
\[
\one^\top\Phi_\sigma(t+s)x\ge\mu_C(t)\,\one^\top\Phi_\sigma(s)x.
\]
Hence $\mu_C$ is supermultiplicative, the limit exists and equals $\sup_tt^{-1}\log\mu_C(t)$. We thus write $\lexp_{\min}=\lexp_K$ and call it the boundary exponent since it controls the decay of every initial state, including those on the boundary of $K$.

\begin{lemma}\label{lem:bounded} The following properties hold:
\begin{enumerate}[label=\textup{(\roman*)}]
\item $\lexp_{\min}\le\lexp_C\le\lexp$. Moreover $\lexp_{\min}\le\eta(A_i[J,J])$ whenever $J\ne\emptyset$ and $A_i$ leaves the face $K\cap\operatorname{span}\{e_j:j\in J\}$ invariant.
\item Let $a:C\to[0,\infty)$ satisfy $c_1\one^\top x\le a(x)\le c_2\one^\top x$ on $C$ for constants $0<c_1\le c_2$, and $a(\Phi_\sigma(t)x)\ge e^{\nu t}a(x)$ for all $\sigma$, $t\ge0$ and $x\in C$. Then $\nu\le\lexp_C$.
\end{enumerate}
\end{lemma}

\begin{proof}
The first inequality holds because $C\cap\Delta\subseteq\Delta$. For $x\in C\cap\Delta$, $\one^\top\Phi_\sigma(t)x\le\|\Phi_\sigma(t)\|_1$, so $\mu_C(t)\le\bar\mu(t)$ and $\lexp_C\le\lexp$. If the face is invariant, a nonnegative Perron vector $x_J$ of $A_i[J,J]$, extended by zero and normalized, satisfies $A_ix=\eta(A_i[J,J])\,x$, so 
\[\mu_K(t)\le\one^\top e^{tA_i}x=e^{\eta(A_i[J,J])t}.
\]
For (ii), any $x\in C\cap\Delta$ satisfies $\one^\top\Phi_\sigma(t)x\ge a(\Phi_\sigma(t)x)/c_2\ge e^{\nu t}a(x)/c_2\ge(c_1/c_2)\,e^{\nu t}$.
\end{proof}

\subsection{Positivity and lower-growth functions}\label{subsec:lowergrowth}

Individual irreducibility removes the boundary obstruction by making transition matrices uniformly positive over each fixed positive time interval.

\begin{lemma}\label{lem:positive}
If every $A_i$ is irreducible, then for every $\tau>0$ there is $\delta>0$ with $\Phi_\sigma(\tau)\ge\delta\,\one\one^\top$ for every switching signal $\sigma$. Consequently, 
\[
\lexp_{\min}=\lexp.
\]
\end{lemma}

\begin{proof}
We may assume $d\ge2$. Write $A_i=N_i-c_AI$, so that $N_i\ge0$, and $\Phi_\sigma(t)=e^{-c_At}\Psi_\sigma(t)$ where $\dot\Psi=N_\sigma\Psi$. The Dyson series has nonnegative terms, so $\Psi_\sigma(b,a)\ge I+\int_a^bN_\sigma$. The proof now follows from this sequence of steps:
\begin{itemize}[leftmargin=*]
\item Split $[0,\tau]$ into $d-1$ intervals of length $\tau/(d-1)$. On each, some mode $i_k$ is active on a set of measure at least a fraction $1/m$ of the interval, so $\Psi_\sigma(\tau)\ge\prod_k(I+\kappa N_{i_k})$ with $\kappa=\tau/((d-1)m)$.
\item Each $N_{i_k}$ is irreducible, so at every factor the support of each column of the product gains at least one index. After $d-1$ factors the product is positive.
\item Every positive entry is a sum of products of at most $d-1$ entries of the $\kappa N_i$ and of $I$. So it is at least $(\min\{1,\kappa\}\nu_0)^{d-1}$, where $\nu_0=\min\{1,\text{smallest positive entry of the }N_i\}$.
\end{itemize}
So the lemma holds with $\delta:=e^{-c_A\tau}(\min\{1,\kappa\}\nu_0)^{d-1}$. Conversely, if some $A_i$ is reducible, then $e^{\tau A_i}$ has a zero entry, and no such $\delta$ exists.
To conclude, for $x\in\Delta$ and $t\ge\tau$, 
\[\one^\top\Phi_\sigma(t)x\ge\delta\,\one^\top\Phi_\sigma(t,\tau)\one\ge\delta\,\bar\mu(t-\tau),
\]
so $\lexp_{\min}\ge\lexp$, and \cref{lem:bounded}(i) gives the reverse inequality.
\end{proof}

The next lemma supplies the lower-growth functions used in the Euler estimates. Part~(i) uses a finite time horizon for any $\nu<\lexp_C$. Under irreducibility, part~(ii) attains the exact rate $\lexp$ and provides the Lipschitz control needed for the $O(h)$ error bound.

\begin{lemma} \leavevmode \label{lem:lower}
\begin{enumerate}[label=(\roman*)]
\item For every $\nu<\lexp_C$ there is a concave, positively homogeneous function $a_T$ on $\R^d$, Lipschitz with some constant $L$ in the $\ell^1$ norm, such that $0\le a_T(x)\le C'\one^\top x$ on $K$, $a_T\ge c_0>0$ on $C\cap\Delta$, and $a_T(e^{sA_i}y)\ge e^{\nu s}a_T(y)$ for $y\in C\cap\Delta$ and $s\in[0,1]$.
\item If every $A_i$ is irreducible, then for every $\tau>0$ the function $a(x)=\inf\{e^{-\lexp t}\one^\top\Phi_\sigma(t)x:\ t\ge\tau,\ \sigma\}$ satisfies $a(e^{sA_i}x)\ge e^{\lexp s}a(x)$ for all $x\in K$ and $s\ge0$, and $c_a\one^\top x\le a(x)\le C_a\one^\top x$ on $K$ for some constants $0<c_a\le C_a$. Moreover $a$ is finite and Lipschitz on an $\ell^1$-neighbourhood $N'$ of $\Delta$, while outside $K\cup N'$ it may take the value $-\infty$.
\end{enumerate}
\end{lemma}

\begin{proof}
Set $\xi=(\lexp_C-\nu)/2$. Since $\Phi_\sigma(t)\ge e^{-c_At}I$, we have $\mu_C(t)\ge e^{-c_At}$. By the convergence $t^{-1}\log\mu_C(t)\to\lexp_C$, there is $c_0>0$ such that 
\[
e^{-\nu t}\one^\top\Phi_\sigma(t)y\ge c_0e^{\xi t} \qquad \text{for all $y\in C\cap\Delta$, all $\sigma$ and all $t\ge0$}.
\]
Fix $\tau>0$, let $C_0=\max_{y\in\Delta}e^{-\nu\tau}\one^\top e^{\tau A_1}y$, and choose $T\ge\tau$ with $c_0e^{\xi T}\ge C_0+1$. Define
\[
 a_T(x)=\inf\bigl\{e^{-\nu t}\one^\top\Phi_\sigma(t)x \ : \ \tau\le t\le T,\ \sigma\bigr\},\qquad x\in\R^d .
\]
It is an infimum of linear functionals whose coefficient vectors $e^{-\nu t}\Phi_\sigma(t)^\top\one$ are nonnegative and bounded by some $C'$. So $a_T$ is finite, concave, positively homogeneous and Lipschitz on $\R^d$, $0\le a_T(x)\le C'\one^\top x$ on $K$, and $c_0\le a_T\le C_0$ on $C\cap\Delta$. Concatenating signals gives
\[
 a_T(e^{sA_i}x)\ge e^{\nu s}\min\{a_T(x),\,b_s(x)\},\quad b_s(x)=\inf\{e^{-\nu t}\one^\top\Phi_\sigma(t)x \ :\ T<t\le T+s,\ \sigma\}
\]
for $s\in(0,1]$ and every $x$. On $C\cap\Delta$ we have $b_s\ge c_0e^{\xi T}\ge a_T+1$, so $a_T(e^{sA_i}y)\ge e^{\nu s}a_T(y)$ there, proving (i).

\smallskip
By \cref{lem:positive} there is $\delta>0$ with $\Phi_\sigma(\tau)\ge\delta\,\one\one^\top$ for every $\sigma$. Each coefficient vector $w=e^{-\lexp t}\Phi_\sigma(t)^\top\one$, $t\ge\tau$, has the form $w=\Phi_\sigma(\tau)^\top u$ with $u\ge0$. So $\delta\|u\|_1\le w_j\le M\|u\|_1$, where $M$ bounds $\Phi_\sigma(\tau)$ entrywise, and $\max_jw_j\le\kappa_0\min_jw_j$ with $\kappa_0=M/\delta$. Then:
\begin{itemize}[leftmargin=*]
\item \emph{Lower bound.} By \cref{prop:basic}(i), $\bar\mu(t)\ge e^{\lexp t}$ so for $y\in\Delta$ we have
\[
w\cdot y\ge\min_jw_j\ge(d\kappa_0)^{-1}\one^\top w=(d\kappa_0)^{-1}e^{-\lexp t}\one^\top\Phi_\sigma(t)\one\ge(d\kappa_0)^{-1}e^{-\lexp t}\bar\mu(t)\ge(d\kappa_0)^{-1}.
\]
So $c_a=(d\kappa_0)^{-1}$ works, while for $C_a$ one can choose $t=\tau$ and $\sigma\equiv1$.
\item \emph{Lipschitz continuity.} Let $C_0'$ bound the functional with $t=\tau$ and $\sigma\equiv1$ on
\[
N' = \bigl\{ x \ : \ \| x - y\|_1 \le 1/(2 \kappa_0) \text{ for some } y \in \Delta \bigr\}.
\]
For each $x\in N'$, choose $y\in\Delta$ with $\|x-y\|_1\le1/(2\kappa_0)$. Then 
\[
w\cdot x\ge\min_jw_j-\|w\|_\infty\|x-y\|_1\ge\tfrac12\min_jw_j.
\]
So the near-minimizers, with $w\cdot x\le a(x)+1$, satisfy $\max_jw_j\le2\kappa_0(C_0'+1)=:B$. Hence on $N'$, $a$ coincides with the infimum over coefficient vectors with $\max_jw_j\le B$, which is Lipschitz with constant $B$.
\item \emph{Invariance.} Concatenation gives $a(e^{sA_i}x)\ge e^{\lexp s}a(x)$ for all $s\ge0$ and $x\in K$.\qedhere
\end{itemize}
\end{proof}

\subsection{Invariant faces and the face formula}\label{subsec:faceformula}

We now identify the boundary rates through the invariant faces and their preserving modes.
For nonempty $J\subseteq\{1,\dots,d\}$ let $F_J=\{x\in K \ :\ x_k=0\text{ for }k\notin J\}$, and let
\[
 S(J)=\{i \ :\ e^{tA_i}F_J\subseteq F_J\text{ for all }t\ge0\} =\{i \ :\ (A_i)_{kj}=0\text{ for all }k\notin J,\ j\in J\}
\]
be the set of modes that leave $F_J$ invariant. If $S(J)\neq\emptyset$, let $\lexp_J$ be the lower exponent of the restricted family $\{A_i[J,J] \ :\ i\in S(J)\}$, on $\R^J$. Thus $\lexp_{\{1,\dots,d\}}=\lexp$ and
\begin{align*}
    S(J)=\{i\} & \implies \lexp_J=\eta(A_i[J,J])
    \\ J=\{j\} & \implies \lexp_J=\min_{i\in S(J)}(A_i)_{jj}.
\end{align*}
For a cone $C$ let $\Sigma(C)=\{\operatorname{supp}x \ :\ x\in C\setminus\{0\}\}$ be the set of supports of its points. 

\begin{definition}
We call a face $F_J$ \emph{fast} if $S(J)\ne\emptyset$ and $\lexp_J<\lexp$, and \emph{slow} otherwise.
\end{definition}

\begin{theorem}[Face formula]\label{thm:faceformula}
Let $C\subseteq K$ be a closed convex cone, $C\ne\{0\}$, invariant under every matrix exponential $e^{tA_i}$. Then
\[
 \lexp_C=\min\bigl\{\lexp_J \ :\ J\in\Sigma(C),\ S(J)\ne\emptyset\bigr\},
\]
and the set on the right is nonempty. In particular, for $C=K$,
\[
 \lexp_{\min}=\min\bigl\{\lexp_J \ :\ J\ne\emptyset,\ S(J)\ne\emptyset\bigr\},
\]
which corresponds to the minimum of the lower exponents of the faces of $K$ that are invariant under at least one mode, $K$ itself included.
\end{theorem}

This formula involves the invariant faces that $C$ touches. Indeed, the minimum is unchanged if $\Sigma(C)$ is replaced by the set of all $J$ with $C\cap F_J\ne\{0\}$. If $y\in C\cap F_J\cap\Delta$ and $x\in F_J$ has support $J$, then $x\ge cy$ for some $c>0$, so the exponent $\lexp_J$ of $x$ under signals in $S(J)$ is at least $\lexp_C$. For a single mode, this reduces to the classical description of the growth of $e^{tA}x$, $x\ge0$, through the classes of the reduced graph~\cite{FS80}.

As a consequence, by \cref{lem:bounded}(ii), lower-growth functions comparable to $\one^\top x$ on $K$, including piecewise-linear ones, cannot establish rates above $\lexp_{\min}$.

Equivalently, $\lexp_C$ is the minimal exponent of $\dot x=(\sum_i u_iA_i)x$, $u(t)\in\Delta_m$, over initial states in $C$, attained by an ergodic control-flow measure~\cite{CK96,CK00,Jen19}. An associated extremal trajectory has nondecreasing support, eventually some $J$. Then, its subsequent control uses only $S(J)$, giving $\lexp_C\ge\lexp_J$.

\begin{figure}[!htbp]
\centering
\begin{tikzpicture}[scale=0.85,
  fastface/.style={circle,fill=red!75!black,inner sep=2.2pt},
  slowface/.style={line width=2.2pt,blue!70!black}]
\begin{scope}
\coordinate (a1) at (0,0); \coordinate (a2) at (4,0); \coordinate (a3) at (2,3.464);
\draw[gray] (a1)--(a2)--(a3)--cycle;
\draw[slowface] (a1)--(a2);
\node[fastface] at (a2) {};
\node[below left] at (a1) {$e_1$}; \node[below right] at (a2) {$e_2$}; \node[above] at (a3) {$e_3$};
\node[below] at (2,-0.08) {\small $F_{\{1,2\}}$ (slow)};
\node[right] at (3.9,0.45) {\small $F_{\{2\}}$ (fast)};
\node at (2,-1.05) {(a) as given};
\end{scope}
\begin{scope}[xshift=7cm]
\coordinate (b1) at (0,0); \coordinate (b2) at (4,0); \coordinate (b3) at (2,3.464);
\fill[gray!25] (b1)--(b2)--(2.4,2.771)--(1.6,2.771)--cycle;
\draw[gray] (b1)--(b2)--(b3)--cycle;
\draw[dashed] (1.6,2.771)--(2.4,2.771);
\draw[slowface] (b1)--(b3);
\node[fastface] at (b3) {};
\node at (2,1.05) {$C\cap\Delta$};
\node[below left] at (b1) {$e_1$}; \node[below right] at (b2) {$e_2$}; \node[above left] at (b3) {$e_3$};
\node[left] at (0.95,1.75) {\small $F_{\{1,3\}}$ (slow)};
\node[right] at (2.2,3.464) {\small $F_{\{3\}}$ (fast)};
\node at (2,-1.05) {(b) transposed};
\end{scope}
\end{tikzpicture}
\caption{Invariant faces of case~148 (\cref{ex:case148}) on the simplex $\Delta$: slow faces as thick edges, fast faces as filled vertices. The shaded set $C\cap\Delta$ for the cone of \cref{obs:cone} touches the slow edge but not the fast vertex; it is drawn with $\beta=1/5$ instead of $1/20$ for visibility.}
\label{fig:faces}
\end{figure}
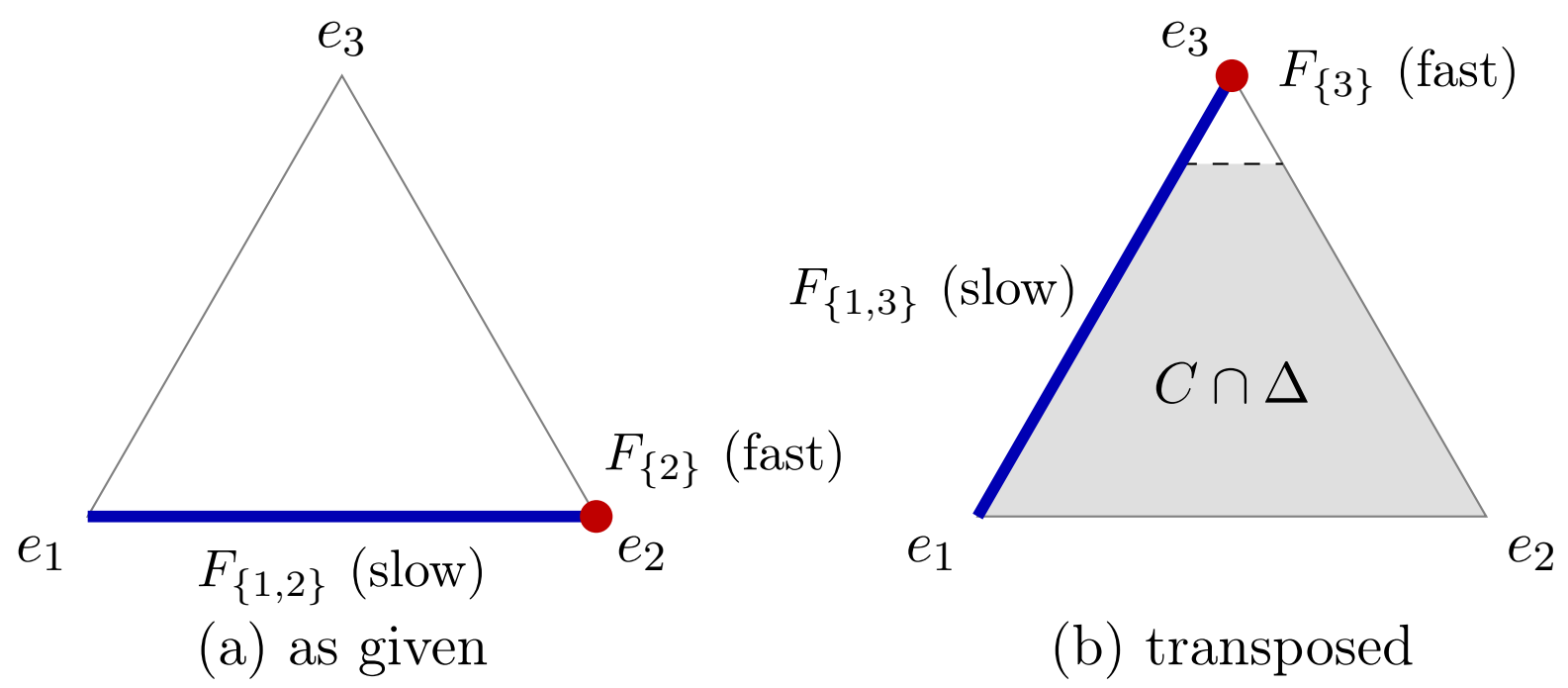

\begin{example}\label{ex:case148}
Case~148 of \cref{app:data} has two reducible modes, \cref{fig:faces} shows its invariant faces. As given, $A_1$ leaves $F_{\{1,2\}}$ invariant, with exponent $\eta(A_1[\{1,2\},\{1,2\}])\approx-1.220815$, and $A_2$ leaves $F_{\{2\}}$ invariant, with exponent $(A_2)_{22}\approx-1.368010$; no proper face is invariant under both modes. By \cref{thm:faceformula}, $\lexp_{\min}=(A_2)_{22}$. Since $-1.247834\le\lexp\le-1.241336$ (\cref{rem:certified}), $F_{\{2\}}$ is fast and $F_{\{1,2\}}$ is slow. After transposition, which changes neither $\lexp$ nor the Euler bounds (\cref{prop:cone}(iv)), $A_1^\top$ leaves $F_{\{3\}}$ invariant, with exponent $(A_1)_{33}\approx-3.669288$ (fast), and $A_2^\top$ leaves $F_{\{1,3\}}$ invariant, with exponent $\approx0.147543$ (slow). The case recurs three times: a cone touching only slow faces proves explicit Euler convergence (\cref{obs:cone}); a polytope avoiding the fast face certifies a lower bound above $\lexp_{\min}$ (\cref{rem:certified}); and the fast face caps every positive antinorm used by the Gripenberg algorithm (\cref{app:gm}).
\end{example}

The proof makes this decomposition explicit through projective occupation measures, as in long-run average control~\cite{Vin93,Gai04,GR06}. The key step, \cref{lem:restriction}, restricts such measures to the points of support $J$, where only modes in $S(J)$ occur and mean growth is at least $\lexp_J$.

\subsubsection{Projective dynamics and occupation measures} 

Throughout the rest of this subsection, fix $C$ as in \cref{thm:faceformula} and $X=C\cap\Delta$. For $x\in K\setminus\{0\}$ let $p=x/\one^\top x\in\Delta$. Under mode $i$, $p$ satisfies
\[
 \dot p=F_i(p):=A_ip-(\one^\top A_ip)\,p,\qquad \tfrac{\mathrm d}{\mathrm d t}\log\one^\top x=\gamma_i(p):=\one^\top A_ip,
\]
so that, for $x\in\Delta$ and the projective solution $p(s)$, it holds
\begin{equation}\label{eq:growth}
 \one^\top\Phi_\sigma(t)x=\exp\int_0^t\gamma_{\sigma(s)}(p(s))\, \mathrm ds , \qquad  |\gamma_i|\le K_A \ \text{ on } \Delta.
\end{equation}
The set $X$ is compact and forward invariant under the projective dynamics. For $x\in X$, a switching signal $\sigma$ and $T>0$, the measure $\frac1T\int_0^T\delta_{(p(s),\sigma(s))}\,\mathrm ds$ is a probability measure on $X\times\{1,\dots,m\}$. Let $\mathcal M$ be the set of weak$^*$ limits of such measures along sequences $(x_n,\sigma_n,T_n)$ with $T_n\to\infty$. Since the probability measures on the compact set $X\times\{1,\dots,m\}$ are weak$^*$ compact, every such sequence has a convergent subsequence. Since $|f(p(T))-f(p(0))|/T\le2\|f\|_{\infty,X}/T\to0$ and $(p,i)\mapsto\nabla f(p)\cdot F_i(p)$ is continuous, every $\mu\in\mathcal M$ satisfies
\begin{equation}\label{eq:invariance}
 \int\nabla f(p)\cdot F_i(p)\, \mathrm d\mu(p,i)=0\qquad\text{for every }f\in C^1(\R^d).
\end{equation}
For $J\ne\emptyset$ let $X_J :=\{p\in\Delta \ :\ \operatorname{supp}p=J\}$, and write $\mu(X_J)$ for $\mu(X_J\times\{1,\dots,m\})$. The sets $X_J$ partition $\Delta$, and $X\subseteq\bigcup_{J\in\Sigma(C)}X_J$.

\begin{lemma}\label{lem:restriction}
Let $\mu\in\mathcal M$ and let $J\subseteq\{1,\dots,d\}$ be nonempty. Then:
\begin{enumerate}[label=\textup{(\alph*)}]
\item for $\mu$-a.e. $(p,i)$ with $p\in X_J$ we have $i\in S(J)$; in particular $\mu(X_J)=0$ if $S(J)=\emptyset$;
\item $\int_{X_J}\nabla f\cdot F\,\mathrm d\mu=0$ for every $f\in C^1(\R^d)$;
\item if $S(J)\ne\emptyset$, then $\int_{X_J}\gamma\,\mathrm d\mu\ge\lexp_J\,\mu(X_J)$.
\end{enumerate}
\end{lemma}

Here $\int_{X_J}\nabla f\cdot F\,\mathrm d\mu$ denotes the integral $\int_{X_J\times\{1,\dots,m\}}\nabla f(p)\cdot F_i(p)\,\mathrm d\mu(p,i)$, and similarly for $\gamma$. The proof is technical and thus given in \cref{app:faceproof}.

\Needspace{6\baselineskip}
\subsubsection{Proof of the face formula}

\begin{proof}[Proof of \cref{thm:faceformula}]
For simplicity, we divide the proof into upper and lower bound respectively. 

\smallskip
\textbf{Upper bound.} Let $J\in\Sigma(C)$ with $S(J)\ne\emptyset$, and choose $x\in X$ with support $J$. Signals with values in $S(J)$ keep $\Phi_\sigma(t)x$ in $F_J$, so $\one^\top\Phi_\sigma(t)x\le\|\Phi_\sigma(t)|_{F_J}\|_1$. Taking their infimum gives $\mu_C(t)\le\bar\mu_J(t)$, where $\bar\mu_J(t)=\inf_{\sigma:\,\sigma(s)\in S(J)}\|\Phi_\sigma(t)|_{F_J}\|_1$, and hence $\lexp_C\le\lexp_J$.

We claim that the set in the minimum is nonempty. Fix a mode $i$. Its projective flow $P_i^t(v)=e^{tA_i}v/\one^\top e^{tA_i}v$ maps the compact convex set $X$ into itself. By Brouwer's theorem, the sets $E_n=\{v\in X:\ P_i^{2^{-n}}(v)=v\}$, $n\ge0$, are nonempty and compact. The semigroup property gives $E_{n+1}\subseteq E_n$, so their intersection contains a point $v$. It is fixed at every nonnegative dyadic time, and hence at every $t\ge0$ by continuity. Differentiation at $t=0$ gives $A_iv=\beta v$ for $\beta=\one^\top A_iv$. If $J=\operatorname{supp}v$ and $k\notin J$, then $0=(A_iv)_k=\sum_{j\in J}(A_i)_{kj}v_j$. Finally, nonnegativity and $v_j>0$ imply $i\in S(J)$.

\textbf{Lower bound.} Let $\lambda_*$ be the minimum in the statement. For $\mu\in\mathcal M$, \cref{lem:restriction}(a),(c) and $X\subseteq\bigcup_{J\in\Sigma(C)}X_J$ give
\[
 \int\gamma\, \mathrm d\mu=\sum_{J\in\Sigma(C),\ S(J)\ne\emptyset}\int_{X_J}\gamma\, \mathrm d\mu\ge\sum_J\lexp_J\,\mu(X_J)\ge\lambda_*.
\]
Now let $\nu<\lambda_*$. If for every $n$ there were $x_n\in X$ and $\sigma_n$ with $\int_0^n\gamma_{\sigma_n}(p_n(s))\, \mathrm ds<\nu n$, a weak$^*$ limit point $\mu\in\mathcal M$ of their occupation measures on $[0,n]$ would satisfy $\int\gamma\, \mathrm d\mu\le\nu<\lambda_*$, a contradiction. So there is $T_0$ with $\int_0^{T_0}\gamma_\sigma(p(s))\, \mathrm ds\ge\nu T_0$ for all $x\in X$ and all $\sigma$. Since $X$ is forward invariant, \eqref{eq:growth} applied to consecutive blocks of length $T_0$ gives 
\[
\one^\top\Phi_\sigma(t)x\ge e^{\nu t-(|\nu|+K_A)T_0} \quad \text{for } x\in X \implies \lexp_C\ge\nu.
\]
Letting $\nu\uparrow\lambda_*$ gives $\lexp_C\ge\lambda_*$.
\end{proof}

\section{Discretization bounds}\label{sec:euler}

Explicit and implicit Euler steps give lower and upper bounds for $\lexp$. We first prove general bounds and implicit convergence, then an $O(h)$ bracket for irreducible modes and cone criteria for reducible families. Sharpness examples and one-sided first-order coefficients complete the section.

\subsection{General Euler bounds}\label{subsec:eulerbounds}

For a finite family $\mathcal F=\{M_1,\ldots,M_r\}$ of square matrices, its lower spectral radius is
\[
 \lsr(\mathcal F)=\lim_{n\to\infty}\min_{|w|=n}\|M_w\|^{1/n}.
\]
For a word $w=w_1\cdots w_k$ we write $M_w=M_{w_k}\cdots M_{w_1}$, so that the first letter acts first. The set of products of a given length, and hence $\lsr$, does not depend on this convention; for instance, in~\cite{GM} we multiply in the other order.

For the Metzler family $\cA$ of \cref{sec:prelim}, set $\cB_h=\{I+hA_i\}$ and $\cR_h=\{(I-hA_i)^{-1}\}$ on the admissible step ranges specified below. Their associated rates are 
\[
g(h)=(\lsr(\cB_h)-1)/h, \qquad G(h)=(1-\lsr(\cR_h)^{-1})/h.
\]
Under the hypotheses of \cref{cor:gripenberg}, the fixed-antinorm Gripenberg algorithm~\cite{GM} approximates these radii to prescribed accuracy.

The explicit bound follows by expressing each flow as a nonnegative sum of Euler iterates and transferring their antinorm growth inequality to the flow.

\Needspace{6\baselineskip}
\begin{theorem}[Explicit Euler lower bound]\label{thm:explicit}
Let $0<h<h_0:=1/c_A$, with $h_0=+\infty$ when $c_A=0$, and define $g(h):=\bigl(\lsr(\cB_h)-1\bigr)/h$. Then, the following properties hold:
\begin{enumerate}[label=\textup{(\roman*)}]
\item $g(h)\le\lexp$ and $g$ is nonincreasing in $h$ (i.e., smaller steps give better bounds);
\item $\lexp_{\min}\le\lim_{h\downarrow0}g(h)\le\lexp$, with $\lexp_{\min}$ as in \cref{sec:prelim}. In particular $g(h)\to\lexp$ when $\lexp_{\min}=\lexp$, for instance when every $A_i$ is irreducible (\cref{lem:positive});
\item if $m=1$, then $g(h)=\eta(A_1)$ for every such $h$.
\end{enumerate}
\end{theorem}

\begin{proof}
Write $B_i=I+hA_i$. These matrices are nonnegative with positive diagonal, so the lower spectral radius is positive:
\[
\rho^\ast:=\lsr(\cB_h)>0.
\]
Define
\[
 a(x)=\inf_w\frac{\one^\top B_wx}{(\rho^\ast)^{|w|}},\qquad x\in K.
\]
As the infimum of linear functionals, $a$ is concave, positively homogeneous, superadditive and monotone on $K$. In addition,
\[
a(x)\le\one^\top x, \quad a(B_ix)\ge\rho^\ast a(x), \quad a(\one)\ge1,
\]
where the last property follows from $\one^\top B_w\one\ge\|B_w\|_1\ge(\rho^\ast)^{|w|}$ by Fekete's lemma.
By uniformization, $e^{tA_i}=e^{-t/h}\sum_k\frac{(t/h)^k}{k!}B_i^k$, a series of nonnegative matrices. Superadditivity and monotonicity therefore give
\[
 a(e^{tA_i}x)\ge e^{-t/h}\sum_k\frac{(t\rho^\ast/h)^k}{k!}a(x)=e^{t(\rho^\ast-1)/h}a(x).
\]
Composing along any switching signal yields 
\[
d\|\Phi_\sigma(t)\|_1\ge\one^\top\Phi_\sigma(t)\one\ge a(\Phi_\sigma(t)\one)\ge e^{tg(h)}.
\]
Now let $h'=\theta h$ with $\theta\in(0,1)$. Then $I+h'A_i=(1-\theta)I+\theta B_i$, and the same $a$ satisfies $a((I+h'A_i)x)\ge(1-\theta+\theta\rho^\ast)a(x)$. Now use the following fact: if $a(B'_ix)\ge\gamma a(x)$ for all $i$ and $x$ and some $\gamma>0$, then $\lsr(\{B'_i\})\ge\gamma$. Indeed, $d\|B'_w\|_1\ge a(B'_w\one)\ge\gamma^{|w|}$. Hence $\lsr(\cB_{h'})-1\ge\theta(\rho^\ast-1)$, that is, $g(h')\ge g(h)$, proving (i).

\smallskip
By (i), $g(h)$ increases as $h\downarrow0$ and is bounded by $\lexp$, so the limit exists and is at most $\lexp$. Fix $\nu<\lexp_{\min}$. \Cref{lem:lower}(i), with $C=K$, gives a concave, positively homogeneous function $a_T$ on $\R^d$, Lipschitz with some constant $L$ in the $\ell^1$ norm, such that:
\begin{itemize}
\item $a_T\ge c_0>0$ on $\Delta$, and $0\le a_T(x)\le C'\one^\top x$ on $K$, because its coefficient vectors are nonnegative and bounded by some $C'$;
\item $a_T(e^{sA_i}y)\ge e^{\nu s}a_T(y)$ for $y\in\Delta$ and $s\in[0,1]$.
\end{itemize}
Let $h\le1$. For $y\in\Delta$, $\|(I+hA_i)y-e^{hA_i}y\|_1\le C_1h^2$ with $C_1=\tfrac12K_A^2e^{K_A}$. Hence
\[
 a_T(B_iy)\ge e^{\nu h}a_T(y)-LC_1h^2\ge\gamma\,a_T(y),\qquad\gamma=e^{\nu h}-LC_1h^2/c_0 ,
\]
and by homogeneity the same holds for every $y\in K$. Take $h$ sufficiently small that $\gamma>0$, which is possible since $\gamma\to1$ as $h\downarrow0$. Then $C'd\,\|B_w\|_1\ge a_T(B_w\one)\ge\gamma^{|w|}a_T(\one)$ for every word, with $a_T(\one)\ge dc_0>0$. Therefore 
\[
\lsr(\cB_h)\ge\gamma \stackrel{\raisebox{3pt}{$\scriptstyle e^x-1\ge x$}}{\implies} g(h)\ge\nu-LC_1h/c_0,
\]
Letting $h\downarrow0$ and then $\nu\uparrow\lexp_{\min}$ gives $\lim g\ge\lexp_{\min}$. When every mode is irreducible, $\lexp_{\min}=\lexp$ by \cref{lem:positive}, concluding the proof of (ii).

Finally, the last property (iii) is an immediate consequence of the fact that the Perron root of $I+hA$ is $1+h\eta(A)$.
\end{proof}

\begin{remark}
By \cref{thm:faceformula}, $\lexp_{\min}$ is the smallest exponent of an invariant face, so $g(h)\to\lexp$ whenever every invariant face is slow. \Cref{prop:cone,cor:eulerface} give sufficient conditions for $g(h)\to\lexp$ that also cover families with $\lexp_{\min}<\lexp$; however, in general the question is open. Refer to \cref{sec:open} for a more thorough discussion.
\end{remark}

\begin{lemma}\label{lem:shadow}
Let $T_0>0$. There is $C$, depending only on $m$, $K_A$ and $T_0$, with the following property. Let $0<h\le1$ with $hK_A\le\frac12$, let $q\ge1$ be an integer with $qh\le T_0$, and let $\sigma$ be a switching signal on $[0,qh]$. Then there is a word $w$ of length $q$ with 
\[
\|R_w-\Phi_\sigma(qh)\|_1\le Ch, \qquad \text{where } R_i=(I-hA_i)^{-1}.
\]
\end{lemma}

\begin{proof}
Let $U_i(t)$ be the time that $\sigma$ spends in mode $i$ up to $t$ and choose $w=w_1\cdots w_q$ so that $w_k$ is a mode maximizing 
\[
U_i(kh)-W_i((k-1)h), \qquad \text{where } \ W_i(t):=|\{s\le t \ : \ w(s)=i\}|.
\]
Here $w(s)=w_k$ for $s\in[(k-1)h,kh)$ is the piecewise-constant signal of the word. At grid times, $W_i((k-1)h)$ is $h$ times the number of letters $i$ among $w_1,\dots,w_{k-1}$. For simplicity's sake, we split the proof into several small steps:

\begin{itemize}[leftmargin=*]
\item \textbf{Discrepancy.} The deficits $U_i-W_i$ at grid times sum to $0$. A deficit decreases only when its mode is chosen, and the chosen mode has deficit at least $h/m$ before the step. So every deficit at grid times lies in $(-h,(m-1)h]$. Between grid times $U_i$ and $W_i$ are nondecreasing and grow by at most $h$, so $D_i=W_i-U_i$ satisfies $|D_i(t)|\le mh$ for all $t$ (and $D_i\equiv0$ if $m=1$). Both $W_i$ and $U_i$ are Lipschitz, and $\mathrm d D_i=(\mathbf 1[w(s)=i]-\mathbf 1[\sigma(s)=i])\, \mathrm ds$.

\item \textbf{Exponential words.} Let $\Phi^w$ be the transition matrix of the signal $w$. For $x(t)=\Phi_\sigma(t)x_0$, variation of constants gives 
\[
\Phi^w(t)x_0-x(t)=\sum_i\int_0^t\Phi^w(t,s)A_ix(s)\, \mathrm d D_i(s).
\]
Integrate by parts, using $D_i(0)=0$, $|D_i|\le mh$, $\partial_s\Phi^w(t,s)=-\Phi^w(t,s)A_{w(s)}$ and $\dot x=A_\sigma x$. This gives $\|\Phi^w(qh)-\Phi_\sigma(qh)\|_1\le C_4h$, with $C_4 := C_4(m, K_A,T_0)$.

\item \textbf{Resolvent words.} For $hK_A\le 1/2$ the Neumann series gives $\|R_i\|_1\le e^{2hK_A}$ and 
\[
\|R_i-e^{hA_i}\|_1\le\sum_{k\ge2}(hK_A)^k\le2K_A^2h^2.
\]
A telescoping sum over the $q\le T_0/h$ factors gives $\|R_w-\Phi^w(qh)\|_1\le2K_A^2T_0e^{2K_AT_0}h$.
\end{itemize}
So the lemma holds with $C:=C_4+2K_A^2T_0e^{2K_AT_0}$.
\end{proof}

For implicit Euler, instead, the resolvent averages the flow over an exponential waiting time. This representation gives the upper bound. \cref{lem:shadow} supplies the approximation needed for convergence.

\begin{theorem}[Implicit Euler upper bound]\label{thm:implicit}
Let $h>0$ satisfy $h\,\eta(A_i)<1$ for all $i$, and define the rate $G(h)=\bigl(1-\lsr(\cR_h)^{-1}\bigr)/h$. Then:
\begin{enumerate}[label=\textup{(\roman*)}]
\item $\lexp\le G(h)$. In particular, every word $w$ of length $k\ge1$ gives the upper bound
\[
 \lexp\le\frac{1-\rho(R_w)^{-1/k}}{h};
\]
\item $G$ is nondecreasing in $h$;
\item $G(h)\to\lexp$ as $h\downarrow0$;
\item if $m=1$, then $G(h)=\eta(A_1)$.
\end{enumerate}
\end{theorem}

\begin{proof}
(i) Since $h\eta(A)<1$, the resolvent can be written in integral form as
\[
(I-hA)^{-1}=\int_0^\infty e^{-s}e^{shA}\, \mathrm ds.
\]
For independent $S_1,\dots,S_k\sim\mathrm{Exp}(1)$, it follows that $R_w=\mathbb E\,\Phi_\sigma(T)$, where $\sigma$ uses the letters of $w$ for the durations $hS_j$ and $T=h\sum_jS_j\sim h\,\Gamma(k,1)$. By \cref{prop:basic}(i), $\|\Phi_\sigma(T)\|_1\ge e^{\lexp T}$. Moreover $h\lexp\le h\min_i\eta(A_i)<1$. Therefore
\[
 d\|R_w\|_1\ge\mathbb E\,\one^\top\Phi_\sigma(T)\one\ge\mathbb E\,e^{\lexp T}=(1-h\lexp)^{-k}.
\]
Hence $\lsr(\cR_h)\ge(1-h\lexp)^{-1}$. The word bound follows from $\lsr\le\rho(R_w)^{1/k}$.

\smallskip
(ii) Let $0<h'<h$ and $\theta=h'/h\in(0,1)$. Write $R_i=(I-hA_i)^{-1}$, $R'_i=(I-h'A_i)^{-1}$ and $\gamma=\lsr(\cR_{h'})$. Then, since $I-hA_i=\theta^{-1}\bigl((I-h'A_i)-(1-\theta)I\bigr)$, we have
\[
 R_i=\theta R'_i\bigl(I-(1-\theta)R'_i\bigr)^{-1}=\theta\sum_{k\ge0}(1-\theta)^k(R'_i)^{k+1}.
\]
The series converges. By the proof of (iv), $\rho(R'_i)=(1-h'\eta(A_i))^{-1}$, and $(1-\theta)<1-h'\eta(A_i)$ because $h\eta(A_i)<1$. For the same reason $(1-\theta)\gamma<1$, since $\gamma\le\min_i\rho(R'_i)$.

Set $a(x):=\inf_w\one^\top R'_wx/\gamma^{|w|}$ for $x\in K$. As in the proof of \cref{thm:explicit}(i), $a$ is concave, positively homogeneous, superadditive and monotone on $K$. It satisfies $a(x)\le\one^\top x$, $a(\one)\ge1$ and $a(R'_ix)\ge\gamma a(x)$. Superadditivity extends to convergent series $\sum_ky_k$ with $y_k\in K$:
\[
a\left(\sum_ky_k\right)\ge a\left(\sum_{k\le n}y_k\right)\ge\sum_{k\le n}a(y_k) \qquad \text{for every } n \in \mathbb N.
\]
Hence
\[
 \begin{aligned}
 a(R_ix)&\ge\theta\sum_{k\ge0}(1-\theta)^k a\bigl((R'_i)^{k+1}x\bigr)\\
 &\ge\theta\sum_{k\ge0}(1-\theta)^k\gamma^{k+1}a(x)=\beta\,a(x),\qquad
 \beta:=\frac{\theta\gamma}{1-(1-\theta)\gamma}.
 \end{aligned}
\]
By the fact used in the proof of \cref{thm:explicit}(i), $\lsr(\cR_h)\ge\beta$. Then 
\[
\lsr(\cR_h)^{-1}\le\beta^{-1}=(\theta\gamma)^{-1}-(1-\theta)/\theta \implies 1-\lsr(\cR_h)^{-1}\ge\theta^{-1}(1-\gamma^{-1}).
\]
Dividing by $h$ gives $G(h)\ge(1-\gamma^{-1})/h'=G(h')$.

\smallskip
(iii) By (i) and (ii), $\lim_{h\downarrow0}G(h)$ exists and is at least $\lexp$. Fix $\varepsilon>0$. By \cref{prop:basic}(i) there are $T>0$ and a switching signal $\sigma$ with $\|\Phi_\sigma(T)\|_1\le e^{(\lexp+\varepsilon)T}$. Take $h=T/N$ with $N$ large. Once $h\le1$ and $hK_A\le 1/2$, \cref{lem:shadow} with $T_0=T$ and $q=N$ gives a word $w$ of length $N$ with $\|R_w-\Phi_\sigma(T)\|_1\le C_Th$, where $C_T$ depends only on $m$, $K_A$ and $T$. By the word bound in (i) and the inequality $1-e^{-y}\le y$,
\[
 G(T/N)\le\frac{1-\rho(R_w)^{-1/N}}{h}\le\frac{\log\rho(R_w)}{Nh}\le\frac1T\log\bigl(\|\Phi_\sigma(T)\|_1+C_Th\bigr).
\]
Letting $N\to\infty$ gives $\lim_{h\downarrow0}G(h)\le\lexp+\varepsilon$. Since $\varepsilon$ is arbitrary, the limit is $\lexp$.

(iv) The eigenvalues of $(I-hA)^{-1}$ are $(1-h\mu)^{-1}$ for $\mu$ in the spectrum of $A$. Write $A=N-c_AI$ with $N\ge0$. Then $|\mu+c_A|\le\rho(N)=\eta(A)+c_A$, so 
\[
|1-h\mu|\ge1+hc_A-h|\mu+c_A|\ge1-h\eta(A)>0,
\]
with equality at $\mu=\eta(A)$. Hence $\rho\bigl((I-hA)^{-1}\bigr)=(1-h\eta(A))^{-1}$.
\end{proof}

\subsection{First-order convergence}\label{subsec:eulerrates}

For irreducible modes, local $O(h^2)$ errors measured by a Lipschitz antinorm give the explicit estimate. Uniform positivity converts the resolvent approximation on fixed time blocks into multiplicative bounds, yielding the implicit estimate.

\begin{theorem}\label{thm:eulerrate}
Suppose every $A_i$ is irreducible. There are $h_1>0$ and $C>0$ such that, for $0<h\le h_1$,
\[
 \lexp-Ch \le g(h) \le \lexp \le G(h) \le \lexp+Ch .
\]
\end{theorem}

\begin{proof}
Fix $\tau>0$. We take $h_1$ so small that every $h\in(0,h_1]$ satisfies:
\begin{itemize}
    \item $h<h_0$ and $h\,\eta(A_i)<1$ for all $i$, so that $g(h)$ and $G(h)$ are defined;
    \item $h\le\min(1,\tau)$; $hK_A\le\tfrac12$; $2C_1h^2<\rho'$;
    \item $C_2h^2<e^{\lexp h}$, with $C_1$, $C_2$ and $\rho'$ as below.
\end{itemize}
The two inner inequalities are \cref{thm:explicit}(i) and \cref{thm:implicit}(i). Now, let $a$ be the function provided by \cref{lem:lower}(ii). We recall that it satisfies:
\begin{itemize}
\item $a(e^{sA_i}x)\ge e^{\lexp s}a(x)$ for all $x\in K$ and $s\ge0$;
\item $c_a\one^\top x\le a(x)\le C_a\one^\top x$ on $K$, for some constants $0<c_a\le C_a$;
\item $a$ is Lipschitz, with some constant $L$, on an $\ell^1$-neighbourhood $N'$ of $\Delta$ of radius $\rho'>0$.
\end{itemize}
Let $y\in\Delta$, $w=e^{hA_i}y$, $z=(I+hA_i)y$ and $s=\one^\top w\in[e^{-K_Ah},e^{K_Ah}]$. Then, for $h \le 1$,
\[
\|z-w\|_1\le C_1h^2, \qquad C_1:=\tfrac12K_A^2e^{K_A}.
\]
So $\|z/s-w/s\|_1\le2C_1h^2<\rho'$, and $z/s\in N'$. By homogeneity, $|a(z)-a(w)|=s\,|a(z/s)-a(w/s)|\le LC_1h^2$. Hence
\[
 a\bigl((I+hA_i)y\bigr)\ge e^{\lexp h}a(y)-LC_1h^2\ge\bigl(e^{\lexp h}-C_2h^2\bigr)a(y),\qquad C_2:=LC_1/c_a ,
\]
and by homogeneity the same holds for every $y\in K$. Set $\gamma=e^{\lexp h}-C_2h^2>0$. As in the proof of \cref{thm:explicit}(i), $C_ad\,\|B_w\|_1\ge a(B_w\one)\ge\gamma^{|w|}a(\one)$ for every word, so 
\[
\lsr(\cB_h)\ge\gamma \implies  g(h)\ge\frac{e^{\lexp h}-1}h-C_2h\ge\lexp-C_2h.
\]
By \cref{lem:positive} and $\Phi_\sigma(t,\tau)\ge e^{-c_A(t-\tau)}I$, there is $\delta_*>0$ with $\Phi_\sigma(t)\ge\delta_*\one\one^\top$ for every switching signal and every $t\in[\tau,2\tau]$. Let $q=\lceil\tau/h\rceil$, so that $qh\in[\tau,2\tau)$. By \cref{lem:shadow} with $T_0=2\tau$, every switching signal $\sigma$ on $[0,qh]$ has a word $w$ of length $q$ with $\|R_w-\Phi_\sigma(qh)\|_1\le C'h$. So every entry of $R_w-\Phi_\sigma(qh)$ is at most $C'h$, and since $\Phi_\sigma(qh)\ge\delta_*\one\one^\top$, we get $0\le R_w\le(1+C_3h)\,\Phi_\sigma(qh)$ entrywise, with $C_3:=C'/\delta_*$.

For $N\ge1$ and $T=Nqh$, choose $\sigma$ with $\|\Phi_\sigma(T)\|_1\le2\bar\mu(T)$ and apply the block construction to each block. Products of nonnegative matrices preserve entrywise inequalities, so the resulting word $w$ of length $Nq$ satisfies $R_w\le(1+C_3h)^N\Phi_\sigma(T)$. By Fekete's lemma $\lsr(\cR_h)\le\|R_w\|_1^{1/(Nq)}$, and $Nq=T/h$ with $T\ge N\tau$, so
\[
 \log\lsr(\cR_h)\le\frac hT\Bigl(N\log(1+C_3h)+\log\bigl(2\bar\mu(T)\bigr)\Bigr)\le\frac{C_3}{\tau}h^2+\frac hT\log\bigl(2\bar\mu(T)\bigr).
\]
Letting $N\to\infty$ gives $\log\lsr(\cR_h)\le h(\lexp+C_3h/\tau)$. Since $1-e^{-y}\le y$, we finally conclude that 
\[
G(h)\le\lexp+(C_3/\tau)h.
\]
\end{proof}

\begin{remark}\label{rem:rates}
The order is \textbf{sharp}. In \cref{ex:sharp}, both errors equal, including for irreducible modes,
\[
(1-\sqrt{1-h^2})/h=h/2+O(h^3),
\]
and under \cref{thm:rateconst} the explicit error, when chattering is optimal, is at least 
\[
Qh/2-O(h^2).
\]
The constant $C_2=LC_1/c_a$, with $c_a=(d\kappa_0)^{-1}$, $L=2\kappa_0(C_0'+1)$ and $\kappa_0=M/\delta$ from \cref{lem:lower}(ii), can exceed the observed widths $0.02h$--$1.5h$ (\cref{tab:gripenberg}); \cref{thm:rateconst} refines it. Without irreducibility, implicit convergence remains general, without a rate. Norms absorb the additive errors in its proof, whereas the lower argument needs a flow-invariant function bounded above and below on $\Delta$, capped at $\lexp_{\min}$ (\cref{lem:bounded}(ii)).
\end{remark}

\subsection{Invariant cones and reducible families}\label{subsec:eulercones}

Without irreducibility, the general estimate only places the explicit limit between $\lexp_{\min}$ and $\lexp$. The next criteria establish convergence when the relevant invariant faces are slow, using cones that the Euler matrices preserve.

\begin{proposition}[Invariant cones]\label{prop:cone}
Let $C\subseteq K$ be a closed convex cone, $C\ne\{0\}$, and assume that $(I+h_CA_i)C\subseteq C$ for all $i$ and some $h_C\in(0,h_0)$. Then $C$ is invariant under every flow $e^{tA_i}$, so $\lexp_C$ is defined (\cref{sec:prelim}). In addition:
\begin{enumerate}[label=\textup{(\roman*)}]
\item $\lim_{h\downarrow0}g(h)\ge\lexp_C$.
\item If $C\setminus\{0\}\subset\operatorname{int}K$, then $\lexp_C\ge\lexp$, so $g(h)\to\lexp$.
\item Suppose that, for some $x_0>0$, the vectors $B_wx_0/\one^\top B_wx_0$, with $w$ ranging over all words of $\cB_{h_C}$, stay in a compact subset of $\operatorname{int}\Delta$. Then the closed convex cone they generate satisfies the hypotheses of \textup{(ii)}.
\item Replacing every $A_i$ by $A_i^\top$ changes neither $g$ nor $\lexp$, so \textup{(i)--(iii)} may be applied to the transposed family.
\end{enumerate}
A cone as in \textup{(ii)} exists only if every $A_i$ has a positive eigenvector. Analogously, for the transposed family, only if every $A_i$ has a positive left eigenvector.
\end{proposition}

\begin{proof}
For $h=\theta h_C$ with $\theta\in(0,1]$, $I+hA_i=(1-\theta)I+\theta(I+h_CA_i)$ maps $C$ into $C$. Since $e^{tA_i}=\lim_n(I+(t/n)A_i)^n$ and $C$ is closed, every $e^{tA_i}$ maps $C$ into $C$, and hence so does every $\Phi_\sigma(t)$ (\cref{sec:prelim}). We now divide the proof into the four statements:

\smallskip
(i) Repeat the proof of \cref{thm:explicit}(ii) with $\Delta$ replaced by $C\cap\Delta$. \Cref{lem:lower}(i) gives, for every $\nu<\lexp_C$, a function $a_T$ with the properties listed there. For sufficiently small $0<h\le h_C$ we have $\gamma>0$. Since $B_iy\in C$, the comparison with $e^{hA_i}y$ gives $a_T(B_iy)\ge\gamma a_T(y)$ on $C$; iteration from any $x_0\in C\setminus\{0\}$ gives $\lsr(\cB_h)\ge\gamma$.

(ii) By compactness, $\varepsilon:=\min\{x_j \ : \ x\in C\cap\Delta,\ j\le d\}$ is strictly positive. For $x\in C\cap\Delta$, $\one^\top\Phi_\sigma(t)x\ge\varepsilon\,\one^\top\Phi_\sigma(t)\one\ge\varepsilon\,\bar\mu(t)\ge\varepsilon e^{\lexp t}$ by \cref{prop:basic}(i). So $\lexp_C\ge\lexp$, and (i) together with \cref{thm:explicit}(i) gives $g(h)\to\lexp$.

(iii) Each $B_i$ maps the ray of each of these vectors onto the ray of another one, so it maps the closed convex cone they generate into itself. If the vectors lie in a compact set $X_0\subset\operatorname{int}\Delta$, the base of that cone, the closed convex hull of the vectors, lies in $\operatorname{conv}X_0$, which is compact and contained in $\operatorname{int}\Delta$.

(iv) If $\tilde\sigma(s)=\sigma(t-s)$, the transition matrix of the transposed family is $\Phi_{\tilde\sigma}(t)^\top$, and products of the matrices $I+hA_i^\top$ are transposes of products of the $B_i$ in reverse order. Since $\|M^\top\|_1\le d\,\|M\|_1$, the exponential rates are unchanged.

The fixed-point argument in the proof of \cref{thm:faceformula} gives a rest point of the projective flow of each $A_i$ in $C\cap\Delta$, that is,
\[
A_iv=\mu v \quad \text{with } v\in C\setminus\{0\}.
\]
Under (ii), $v>0$. If $u\ge0$, $u\neq0$, is a left eigenvector of $A_i$ for $\eta(A_i)$, then $\eta(A_i)u^\top v=u^\top A_iv=\mu u^\top v$ with $u^\top v>0$, so $\mu=\eta(A_i)$.
\end{proof}

\begin{corollary}\label{cor:eulerface}
Let $C$ be as in \cref{prop:cone}. Then
\[
 \lim_{h\downarrow0}g(h)\ \ge\ \min\bigl\{\lexp_J \ :\ J\in\Sigma(C),\ S(J)\ne\emptyset\bigr\}.
\]
In particular, if every face $F_J$ with $J\in\Sigma(C)$ is slow, then $g(h)\to\lexp$.
\end{corollary}

\begin{proof}
Combine \cref{prop:cone}(i) and \cref{thm:faceformula}. If every face in the (nonempty) minimum is slow, the minimum is at least $\lexp$, so $\lexp\ge\lim g(h)\ge\lexp$ by \cref{thm:explicit}(i).
\end{proof}

\Cref{prop:cone}(ii) corresponds to $\Sigma(C)=\{\{1,\dots,d\}\}$; the corollary also permits contact with slow boundary faces. For $C=K$, it recovers \cref{thm:explicit}(ii) through \cref{thm:faceformula}, giving convergence when all invariant faces are slow.

An Euler-invariant cone containing $x>0$ contains the limits of $B_i^nx/(\one^\top B_i^nx)$, including the dominant eigenvector when unique. If that vector lies on a fast face, the criterion cannot apply and convergence remains undecided. A cone inside $\partial K$ lies in a common invariant face and helps only when that face has exponent $\lexp$.

\begin{computation}\label{obs:cone}
For cases 77, 148, 677, 1195 and 1440 (\cref{app:data}), only the last has two irreducible modes. In both orientations, every proper invariant face has $S(J)=\{i\}$. Its exponent $\eta(A_i[J,J])$ is bounded below exactly by Collatz--Wielandt with a nonnegative eigenvector. We prove slowness by comparison with a certified upper bound for $c=\min_\alpha\eta(S_\alpha)\ge\lexp$, using \cref{prop:basic}(iii) at rational $\alpha$ near the minimizer.
\begin{itemize}
\item \emph{All invariant faces slow.} In case~1195 as given, the only invariant proper face is $F_{\{2,3\}}$, invariant under $A_1$, with exponent at least $-0.289318$, while $c\le-0.532919$. 

In case~677 as given it is $F_{\{2\}}$, invariant under $A_2$, with exponent $-0.237097$, while $c\le-0.394492$. 

In case~77 after transposition it is $F_{\{2\}}$, invariant under $A_2^\top$, with exponent $0.024892$, while $c\le-0.148186$. Case~1440 has none. By \cref{thm:faceformula}, $\lexp_{\min}=\lexp$ for these families, and \cref{thm:explicit}(ii) gives $g(h)\to\lexp$.

The orientation matters, although $g$ and $\lexp$ do not depend on it (\cref{prop:cone}(iv)): case~77 as given and cases~677 and~1195 transposed each have an invariant face with exponent below $c$, which this comparison cannot show to be slow. For cases~77 and~677 in these orientations, the invariant polygons of \cref{prop:cone}(ii) (with $h_C=0.2$, of 21 and 53 vertices, coordinates at least $0.0796$ and $0.1137$) give the same conclusion.

\item In case~148 both orientations have a fast face (\cref{ex:case148}), so we use \cref{cor:eulerface} for the transposed family, with $h_C=1/5<h_0$ (\cref{fig:faces}(b)). Its other invariant face, $F_{\{1,3\}}$ under $A_2^\top$, has exponent at least $0.147542$, while $c\le-1.236979$, so it is slow. The cone
\[
C=\{x\ge0 \ :\ (1-\beta)(x_1+x_2)\ge\beta x_3\}, \qquad \beta=1/20,
\]
has extreme rays $e_1$, $e_2$, $(\beta,0,1-\beta)$ and $(0,\beta,1-\beta)$, and cuts off only $e_3$. The supports of its points are all nonempty subsets of $\{1,2,3\}$ except $\{3\}$, and the only proper invariant face among them is the slow face $F_{\{1,3\}}$. Both Euler matrices map the four extreme rays into $C$. This is checked in exact rational arithmetic, with binary64 data, and replayed by an independent method (barycentric coordinates in a fan triangulation). By \cref{cor:eulerface}, $g(h)\to\lexp$. So $g(h)\to\lexp$ for all five test pairs.

\item \emph{Case 148 without a cone.} In the given orientation $\lexp_{\min}=(A_2)_{22}$ (\cref{ex:case148}). The polytopes of \cref{rem:certified}, which have no vertex on the axis of $e_2$, certify $\lim_{h\downarrow0}g(h)\ge-1.247834>\lexp_{\min}$ directly.
\end{itemize}
\end{computation}

\subsection{Sharpness and refined coefficients}\label{subsec:eulercoefficients}

We first give an exact example showing that the linear rate is sharp, then derive one-sided coefficients for two-mode families.

\begin{example}\label{ex:sharp}
Let $\varepsilon\ge0$ and $A_{1,2}=\begin{psmallmatrix}\mp1&\varepsilon\\\varepsilon&\pm1\end{psmallmatrix}$. For $\varepsilon=0$ the modes are diagonal; for $\varepsilon>0$ they are irreducible. For $x>0$, $a(x)=\sqrt{x_1x_2}$ satisfies, along both modes,
\[
\frac{\mathrm d}{\mathrm d t}\log a(x(t))=\frac\varepsilon2\Bigl(\frac{x_2}{x_1}+\frac{x_1}{x_2}\Bigr)\ge\varepsilon,
\]
so $\lexp\ge\varepsilon$. Since $\eta(S_{1/2})=\varepsilon$, \cref{prop:basic}(iii) gives $\lexp=\varepsilon$. Since $\eta(A_i)=\sqrt{1+\varepsilon^2}$, let $0<h<(1+\varepsilon^2)^{-1/2}$, so that $h<h_0=1$ and $h\,\eta(A_i)<1$, and put $r=\sqrt{1-h^2}+\varepsilon h$.
\begin{itemize}
\item \emph{Explicit family.} Let $a(x)=\sqrt{x_1x_2}$ and $t=x_2/x_1$. By AM--GM,
\[
 \frac{a(B_1x)^2}{a(x)^2}=\frac{(1-h)h\varepsilon}{t}+1-h^2+h^2\varepsilon^2+h\varepsilon(1+h)\,t\ \ge\ r^2 ,
\]
and the same holds for $B_2$. Since $a(x)\le\tfrac12\one^\top x$ and $a(\one)=1$, the fact used in the proof of \cref{thm:explicit}(i) gives $\lsr(\cB_h)\ge r$. Conversely, $\lsr(\cB_h)\le\rho(B_1B_2)^{1/2}$ and 
\[
\rho(B_1B_2)=1-h^2(1-\varepsilon^2)+2h\varepsilon\sqrt{1-h^2}=r^2 \implies \lsr(\cB_h)=r.
\]
\item \emph{Implicit family.} For a $2\times2$ matrix $A$ with $\operatorname{tr}A=0$, $\operatorname{adj}(I-hA)=I+hA$. So $(I-hA_i)^{-1}=B_i/D$ with $D=\det(I-hA_i)=1-h^2(1+\varepsilon^2)=(\sqrt{1-h^2}-\varepsilon h)\,r$, the same for both modes. Hence $\lsr(\cR_h)=r/D$.
\end{itemize}
Therefore
\[
 g(h)=\lexp-\frac{1-\sqrt{1-h^2}}{h},\quad G(h)=\lexp+\frac{1-\sqrt{1-h^2}}{h},\quad\frac{1-\sqrt{1-h^2}}{h}=\frac h2+O(h^3).
\]
For every $\varepsilon\ge0$, and in particular for irreducible modes, both one-sided errors are of exact order $h$, and the midpoint $(g+G)/2$ equals $\lexp$.
\end{example}

\begin{remark}\label{rem:comm}
Let $m=2$, and suppose $\eta=\eta(S_\alpha)$ is a simple eigenvalue with left and right eigenvectors $u,v$. Writing $F=A_1-A_2$, one has $[A_1,A_2]=[F,S_\alpha]$ and hence 
\[
u^\top[A_1,A_2]v=\eta\,u^\top Fv-\eta\,u^\top Fv=0.
\]
At $\alpha=1/2$, under the same simplicity assumption, Baker--Campbell--Hausdorff gives 
\[
\log\rho(e^{hA_2}e^{hA_1})/(2h)=\eta(S_{1/2})+O(h^2).
\]
The closest statement is~\cite[Prop.~10]{CGG14}: periodic perturbations of the best constant control change the Floquet eigenvalue only at second order. For general noncommuting matrices, splitting formulas with positive durations have order at most two~\cite{Sheng,Suzuki}. Note that special families, such as commuting ones, can be exact.
\end{remark}

We next derive explicit one-sided first-order bounds for the Euler errors at a two-mode chattering optimum.

\begin{theorem}\label{thm:rateconst}
Let $m=2$, $F=A_1-A_2$ and $S_a=aA_1+(1-a)A_2$ for $a\in[0,1]$. Suppose that $a^\ast\in(0,1)$ minimizes $\eta(S_a)$ over $[0,1]$ and that $S=S_{a^\ast}$ is irreducible. Let $c=\eta(S)$, $u^\top S=cu^\top$, $Sv=cv$, $u^\top v=1$, and
\[
 Q=a^\ast(1-a^\ast)\,u^\top F^2v=u^\top\bigl(a^\ast A_1^2+(1-a^\ast)A_2^2\bigr)v-c^2 .
\]
There are $h_2>0$ and $C>0$ such that, for $0<h\le h_2$,
\[
 g(h)\le c-\tfrac Q2h+Ch^2,\qquad G(h)\le c+\tfrac Q2h+Ch^2 .
\]
If moreover chattering is optimal, $\lexp=c$, as for every planar pair~\cite[Thm.~1]{BCV11}, then
\[
 \lexp-g(h)\ge\tfrac Q2h-Ch^2,\qquad 0\le G(h)-\lexp\le\tfrac Q2h+Ch^2,
\]
and $Q\ge0$. For $d=2$, $Q=\tfrac12a^\ast(1-a^\ast)(c-\lambda_2)\,\eta''(a^\ast)$, where $\lambda_2$ is the other eigenvalue of $S$ and $\eta(a)=\eta(S_a)$. Hence, $Q>0$ when the minimum is nondegenerate.
\end{theorem}

\begin{proof}
Every $S_a$, $0<a<1$, has the off-diagonal pattern of $S$, so $\eta(a)$ is a simple eigenvalue of $S_a$, analytic in $a$, with Perron vectors $u_a,v_a$ ($u_a^\top v_a=1$) and $\eta'(a^\ast)=u^\top Fv=0$. Set 
\[
\kappa(a)=u_a^\top(aA_1^2+(1-a)A_2^2)v_a, \qquad Q(a)=\kappa(a)-\eta(a)^2,
\]
so that $Q=Q(a^\ast)$. For $1\le p<q$ let $w(p,q)=w_1\cdots w_q$ be such that
\begin{equation} \label{eq:chr_word}
w_j = \begin{cases}1 & \text{if } \lfloor jp/q\rfloor>\lfloor(j-1)p/q\rfloor,
\\ 2 & \text{otherwise}.\end{cases}
\end{equation}
Fix a compact subinterval $\mathcal I\subset(0,1)$ with $a^\ast\in\operatorname{int}\mathcal I$. \Cref{lem:wordexp} in \cref{app:rateproof} shows that, uniformly for $b\in \mathcal I$ and $T=qh\le T_1$, with $T_1$ fixed,
\[
 \frac{\log\rho(B_w)}T=\eta(b)-\frac h2\kappa(b)+O(h^2),\qquad\frac{\log\rho(R_w)}T=\eta(b)+\frac h2\kappa(b)+O(h^2).
\]
Take $q=\lfloor T_1/h\rfloor$ and $p$ nearest to $qa^\ast$. Reduce $h_2$ so that $h\le T_1/2$, $1\le p<q$ and $b\in \mathcal I$; then $|b-a^\ast|\le h/T_1$. With $\zeta=\log\rho(B_w)/T$ we have $\lsr(\cB_h)\le\rho(B_w)^{1/q}=e^{h\zeta}$, so
\[
 g(h)\le\frac{e^{h\zeta}-1}h\le\zeta+\frac h2\zeta^2+O(h^2)=\eta(b)-\frac h2Q(b)+O(h^2).
\]
Since $\eta'(a^\ast)=0$, $\eta(b)=c+O(h^2)$, and $Q(b)=Q+O(h)$. The implicit side is the same. Indeed, by applying \cref{thm:implicit}(i) we get
\[
G(h)\le(1-e^{-h\zeta'})/h\le\zeta'-\frac h2\zeta'^2+O(h^2), \qquad \zeta'=\log\rho(R_w)/T.
\]
When $\lexp=c$, the second pair of inequalities in the statement follows from these bounds and $g(h)\le\lexp\le G(h)$. Dividing the upper bound for $G(h)-\lexp$ by $h$ and letting $h\downarrow0$ gives $Q\ge0$. For $d=2$, write $I=vu^\top+v_2u_2^\top$ with the eigenvectors of $\lambda_2$. Here
\[
u^\top Fv=0, \qquad u^\top F^2v=(u^\top Fv_2)(u_2^\top Fv),
\]
and second-order perturbation gives $\eta''(a^\ast)=2(u^\top Fv_2)(u_2^\top Fv)/(c-\lambda_2)$.
\end{proof}

When $\lexp=c$ and $Q>0$, the explicit error is at least $Qh/2-O(h^2)$; for irreducible modes, \cref{thm:eulerrate} makes it of order $h$ (as in \cref{ex:sharp}, with $Q=1$). The implicit conclusion is only $\limsup_{h\downarrow0}(G(h)-\lexp)/h\le Q/2$; matching lower estimates remain open (\cref{sec:open} and \cref{app:tables}).

\begin{observation}\label{obs:sandwich}
We test \cref{ex:sharp} with $\varepsilon=0.2$ and the asymmetric planar pair given by
\[
A_1=\begin{psmallmatrix}-2&0.3\\0.5&1\end{psmallmatrix}, \qquad A_2=\begin{psmallmatrix}0.5&0.4\\0.1&-1.5\end{psmallmatrix}.
\]
For the latter,~\cite[Thm.~1]{BCV11} gives $\lexp=\min_\alpha\eta(S_\alpha)\approx-0.1949763083$. In \cref{tab:sandwich}, the finite-word gap satisfies $G_{\rm up}-g_{\rm est}\approx hQ$, with 
\[
Q=u^\top\left(\sum_i\alpha_i^\ast(A_i-S_{\alpha^\ast})^2\right)v
\]
as in \cref{thm:rateconst} (predicted slope $1.478$, observed $1.470$). The computed midpoint $(g_{\rm est}+G_{\rm up})/2$ suggests second-order cancellation, while the exact midpoint $(g+G)/2$ equals $\lexp$ for the symmetric pair. These finite-search data do not establish matching coefficients for $g,G$. Exponential grid schedules converge quadratically in the symmetric case; for the asymmetric pair, the word-length cap limits approximation of $\alpha^\ast\approx0.3905$.
\end{observation}

\begin{table}[!htbp]
\centering\small
\caption{Euler word estimates and exponential schedules, words of length at most~16. The implicit-word value $G_{\rm up}$ and the schedule rate bound $\lexp$ from above up to rounding; $g_{\rm est}$ estimates $g(h)$. For the symmetric pair, the values match \cref{ex:sharp}.}
\label{tab:sandwich}
\begin{tabular}{@{}r|ccc|ccc@{}}
\toprule
& \multicolumn{3}{c|}{symmetric ($\lexp=0.2$)} & \multicolumn{3}{c}{asymmetric ($\lexp\approx-0.194976$)}\\
$h$ & $g_{\rm est}$ & $G_{\rm up}$ & \shortstack{exp. schedule} & $g_{\rm est}$ & $G_{\rm up}$ & \shortstack{exp. schedule}\\
\midrule
0.2    & 0.098979 & 0.301021 & 0.201335 & $-0.358479$ & $-0.056019$ & $-0.188068$\\
0.1    & 0.149874 & 0.250126 & 0.200333 & $-0.271832$ & $-0.123611$ & $-0.192983$\\
0.05   & 0.174984 & 0.225016 & 0.200083 & $-0.232157$ & $-0.158550$ & $-0.194214$\\
0.025  & 0.187498 & 0.212502 & 0.200021 & $-0.213186$ & $-0.176426$ & $-0.194521$\\
0.0125 & 0.193750 & 0.206250 & 0.200005 & $-0.203857$ & $-0.185482$ & $-0.194598$\\
\bottomrule
\end{tabular}
\end{table}

\begin{remark}\label{rem:freqavg}
Let $\psi(\omega)=\operatorname{Re}u^\top F(c+\imath \omega-S)^{-1}Fv$ for $\omega>0$. Then (\cref{app:rateproof})
\[
 Q=\frac{2a^\ast(1-a^\ast)}\pi\int_0^\infty\psi(\omega)\, \mathrm d\omega.
\]
To derive the associated Floquet formula, let $f$ be real, bounded, measurable and $P$-periodic with zero mean, where $P=2\pi/\omega_0$ with $\omega_0>0$, and set 
\[
\hat f_k=P^{-1}\int_0^P f(t)e^{-\imath k\omega_0t}\, \mathrm dt.
\]
Write $\langle\cdot\rangle$ for the period average. For fixed $f,P$ and small $\varepsilon$, the system $\dot x=(S+\varepsilon f(t)F)x$ has a Floquet solution $x(t)=e^{\lambda_\varepsilon t}p_\varepsilon(t)$, with $p_\varepsilon$ periodic and $\langle u^\top p_\varepsilon\rangle=1$. The monodromy is analytic in $\varepsilon$ and its Perron eigenvalue at $\varepsilon=0$ is simple, so write
\[
 p_\varepsilon=v+\varepsilon z+\varepsilon^2r+o(\varepsilon^2),\qquad
 \lambda_\varepsilon=c+\varepsilon\lambda_1+\varepsilon^2\lambda_2+o(\varepsilon^2).
\]
Substitution into $\dot p_\varepsilon=(S+\varepsilon fF-\lambda_\varepsilon I)p_\varepsilon$ and averaging against $u^\top$ give $\lambda_1=\langle f\rangle u^\top Fv=0$ and
\[
 \dot z=(S-cI)z+fFv,\qquad \lambda_2=\langle f\,u^\top Fz\rangle.
\]
The normalization and zero mean give $\hat z_0=0$; for $k\ne0$,
$\hat z_k=\hat f_k(c+\imath k\omega_0-S)^{-1}Fv$. These inverses exist because all other eigenvalues of the irreducible Metzler matrix $S$ have real part below $c$. Parseval's identity and conjugate pairing therefore yield
\[
 \lambda_\varepsilon-c
 =2\varepsilon^2\sum_{k\ge1}|\hat f_k|^2\psi(k\omega_0)+o(\varepsilon^2).
\]
The series is absolutely convergent. Indeed, $f\in L^2(0,P)$ and the resolvents are bounded over the nonzero integer frequencies. Thus $Q$ is a frequency average of this quadratic coefficient. If chattering is optimal, periodic convexified controls cannot improve on $c$ (by approximation with switching). Taking $f(t)=\cos(\omega_0t)$, with $a^\ast+\varepsilon f(t)\in[0,1]$, gives $\psi(\omega_0)\ge0$ for every $\omega_0>0$, hence $Q\ge0$.
\end{remark}

\section{Rigorous brackets}\label{sec:gripenberg}

In this section, we discuss lower bounds obtained from polytopes and upper bounds obtained from implicit words, mixtures or schedules. We distinguish algorithmic guarantees, exact verification and numerical performance.

\subsection{Algorithmic guarantees}\label{subsec:algorithmic}

The Euler bounds turn discrete lower-spectral-radius algorithms into continuous-time brackets. For individually irreducible modes, the fixed-antinorm Gripenberg algorithm of~\cite{GM} yields a convergent scheme.

\begin{corollary}[Continuous-time Gripenberg bounds]\label{cor:gripenberg}
Suppose every $A_i$ is irreducible, and let $0<h<h_0$ with $h\,\eta(A_i)<1$ for all $i$.
\begin{enumerate}[label=(\roman*)]
\item The families $\cB_h$ and $\cR_h$ are primitive. In other words, every product of $d-1$ matrices of $\cB_h$, and every matrix of $\cR_h$, is positive.

\item Run the fixed-antinorm Gripenberg algorithm~\cite[Sec.~3.1]{GM}, with a fixed positive antinorm and accuracy $\delta>0$, on $\cB_h$ and on $\cR_h$. Let $t_k$ be the lower bound of the first run and $s_k$ the upper bound of the second, after $k$ levels. Then, for every $k$,
\[
 \frac{t_k-1}h\ \le\ \lexp\ \le\ \frac{1-s_k^{-1}}h .
\]
Write $W_k(h,\delta)$ for the width of this bracket. If $h\le h_1$, with $h_1$ and $C$ as in \cref{thm:eulerrate}, then, for each fixed $h$ and $\delta>0$,
\[
 \limsup_{k\to\infty}W_k(h,\delta)\le2Ch+5\delta h^{-1}.
\]
For $\delta=\delta(h)=o(h)$, choose $k(h)$ so that $W_{k(h)}(h,\delta(h))\le2Ch+5\delta(h)/h+h$. The brackets then converge to $\lexp$ as $h\downarrow0$; with $\delta=h^2$ their width is at most $(2C+6)h$.
\end{enumerate}
\end{corollary}

\begin{proof}
For $h<h_0$, the matrix $I+hA_i$ is nonnegative, has positive diagonal and has the irreducible pattern of $A_i$. The argument in the proof of \cref{lem:positive} shows that every product of $d-1$ such matrices is positive. For the resolvents, write $I-hA_i=(1+hc_A)I-hN_i$ with $N_i=A_i+c_AI\ge0$ irreducible. Since $\rho(hN_i)=h(\eta(A_i)+c_A)<1+hc_A$, we have 
\[
(I-hA_i)^{-1}=(1+hc_A)^{-1}\sum_{k\ge0}\bigl(hN_i/(1+hc_A)\bigr)^k,
\]
which is positive and proves (i). For (ii), Guglielmi and Maiale~\cite[Secs.~2.3 and~3.1]{GM} state that primitive families satisfy the hypotheses of their convergence theorem; for completeness we verify the cone condition. Every product $P$ is primitive, hence asymptotically rank-one, and its secondary spectral subspace is the kernel of its positive left Perron vector $u_P$. For a product of length at least $d-1$, write $P=P'M$ with $M$ one of the finitely many positive products of length $d-1$. Then $u_P^\top$ is proportional to $(u_P^\top P')M$, so the normalized vectors $u_P$ have entries bounded below uniformly. Hence the closure of the union of the secondary subspaces meets $K$ only at $0$, which is the cone condition. So the main theorem applies to both families. The first run gives $t_k\le\lsr(\cB_h)\le s'_k$ and the second $t'_k\le\lsr(\cR_h)\le s_k$, with $\limsup_{k\to\infty}(s'_k-t_k)\le\delta$ and $\limsup_{k\to\infty}(s_k-t'_k)\le\delta$, with bounds held constant once the prescribed accuracy is attained. The bracket follows from \cref{thm:explicit}(i), \cref{thm:implicit}(i), and the monotonicity of $x\mapsto(x-1)/h$ and $x\mapsto(1-x^{-1})/h$. For the width, we have
\[
t_k\ge\lsr(\cB_h)-\delta-o(1) \implies(t_k-1)/h\ge g(h)-\delta/h-o(1)\ge\lexp-Ch-\delta/h-o(1).
\]
For the upper end, $s_k\le\lsr(\cR_h)+\delta+o(1)$, and $\lsr(\cR_h)\ge(1-h\lexp)^{-1}\ge\tfrac12$ by the proof of \cref{thm:implicit}(i), since $h|\lexp|\le hK_A\le\frac12$ by \cref{prop:basic}(i) and the choice of $h_1$. Finally 
\[
\rho^{-1}-(\rho+\delta)^{-1}\le \frac{\delta}{\rho^2} \implies \frac{1-s_k^{-1}}{h} \le G(h)+4 \frac{\delta}{h} + o(1) \le \lexp+Ch+4 \frac{\delta}{h}+o(1).
\]
\end{proof}

\begin{remark}\label{rem:adaptive}
\Cref{cor:gripenberg} covers a fixed antinorm. The adaptive variant~\cite[Sec.~3.2]{GM} enlarges the unit polytope; its spectral-radius upper bounds remain valid in exact arithmetic. Its lower-bound argument assumes that enlargement cannot decrease the operator antinorm $a(P)=\min_{v\in V}a(Pv)$. This monotonicity fails in general.
\end{remark}

\subsection{Construction and limitations of the bounds}\label{subsec:certificates}

The next two propositions also apply to reducible families. Their bounds are valid whenever the stated inequalities hold, independently of how the candidate vertices or words were found. We use the word \emph{certificate} for finite data whose defining inequalities can be checked exactly to prove a bound.

\begin{proposition}[Polytope certificates]\label{prop:polytope}
Let $0<h<h_0$ and $B_i=I+hA_i$, let $V\subset K\setminus\{0\}$ be finite and nonempty, and let $t>0$. Suppose that for every $v\in V$ and every $i$ there are weights $\lambda_u\ge0$, $u\in V$, with $\sum_u\lambda_u=1$ and
\begin{equation}\label{eq:polycert}
 B_iv\ \ge\ t\sum_{u\in V}\lambda_u\,u .
\end{equation}
\begin{enumerate}[label=\textup{(\roman*)}]
\item Then $\lsr(\cB_h)\ge t$, and hence $\lexp\ge(t-1)/h$.
\item Let $0<h'\le h$ and $\theta=h'/h$. The same set $V$ satisfies \eqref{eq:polycert} for the matrices $I+h'A_i$, with $t'=1-\theta+\theta t$ in place of $t$. So it yields the same bound $(t'-1)/h'=(t-1)/h$ at every smaller step.
\item Replacing every $A_i$ by $A_i^\top$ does not change $\lsr(\cB_h)$, so the same inequalities may be checked for the transposed family instead.
\end{enumerate}
\end{proposition}

\begin{proof}
(i) Let $\mathcal P=\operatorname{conv}V+K$. It is convex with $\mathcal P+K=\mathcal P$. By \eqref{eq:polycert} and $B_i\ge0$, $B_i\mathcal P\subseteq t\mathcal P$, so $B_wv\in t^{|w|}\mathcal P$ for every word $w$ and every $v\in V$. Since $\one^\top y\ge c_V:=\min_{u\in V}\one^\top u>0$ on $\mathcal P$,
\[
 \|B_w\|_1\ge\frac{\one^\top B_wv}{\one^\top v}\ge\frac{c_V}{\one^\top v}\,t^{|w|}.
\]
Hence $\lsr(\cB_h)\ge t$, and $\lexp\ge(t-1)/h$ by \cref{thm:explicit}(i). This is the lower bound of~\cite[Prop.~6]{GP13} (on antinorms see also~\cite{GZ20}), or the argument of \cref{thm:explicit}(i), for the antinorm with unit ball $\mathcal P$. We include the proof because this antinorm need not be positive on $K$; it vanishes on the faces of $K$ that $\mathcal P$ does not meet, as in case~148 (\cref{rem:certified}).

(ii) We have $0<h'<h_0$ and
\[
 (I+h'A_i)v=(1-\theta)v+\theta B_iv\ \ge\ (1-\theta)v+\theta t\sum_u\lambda_uu .
\]
The right-hand side is $t'$ times a convex combination of points of $V$, with weight $(1-\theta)/t'$ on $v$ and $\theta t\lambda_u/t'$ on each $u$.

(iii) Products of the matrices $I+hA_i^\top$ are transposes of products of the $B_i$ in reverse order, and $\|M^\top\|_1\le d\,\|M\|_1$.
\end{proof}

\begin{proposition}[Certified upper bounds]\label{prop:upper}
Let $h>0$, let $w=w_1\cdots w_k$ be a word of length $k\ge1$, and assume that the following holds:
\begin{enumerate}[label=\textup{(\alph*)}]
\item for every $i$ there is $z_i>0$ with $(I-hA_i)z_i>0$;
\item there are vectors $x_0>0$ and $x_1,\dots,x_k$ with $(I-hA_{w_j})x_j\ge x_{j-1}$ for $j=1,\dots,k$.
\end{enumerate}
Let $q=\max_\ell(x_k)_\ell/(x_0)_\ell$, and let $r>0$ satisfy $r^k\ge q$. Then $h\,\eta(A_i)<1$ for every $i$, and 
\[
\lexp\le \frac1h (1-r^{-1}).
\]
\end{proposition}

\begin{proof}
By (a) and the Collatz--Wielandt bound of \cref{prop:basic}(iii), applied to $hA_i$, we have $h\,\eta(A_i)\le\max_\ell(hA_iz_i)_\ell/(z_i)_\ell<1$. So, as in the proof of \cref{thm:implicit}(i),
\[
R_i=(I-hA_i)^{-1}=\int_0^\infty e^{-s}e^{shA_i}\, \mathrm ds\ge0.
\]
By (b), $x_j=R_{w_j}(I-hA_{w_j})x_j\ge R_{w_j}x_{j-1}$. Since every $R_i\ge0$, induction gives $x_k\ge R_wx_0\ge0$. Collatz--Wielandt for the nonnegative matrix $R_w$ and $x_0>0$ gives $\rho(R_w)\le\max_\ell(R_wx_0)_\ell/(x_0)_\ell\le q\le r^k$. Finally, \cref{thm:implicit}(i) gives 
\[
\lexp\le(1-\rho(R_w)^{-1/k})/h\le(1-r^{-1})/h.
\]
\end{proof}

\subsubsection{Numerical search and exact verification} 

Floating-point searches propose candidates, which are then checked in exact rational arithmetic and verified independently. \Cref{app:protocols} gives the polytope, implicit-word, mixture and schedule procedures. Its primal--dual precision concerns the chosen polytope, not $\lexp$; unfinished searches are inconclusive.

\begin{remark}\label{rem:certified}
The bounds in \cref{tab:gripenberg} hold for exact-decimal and binary64 data. Case~1440 at $h=0.025$ reuses its $h=0.05$ polytope.

For case~148, we close $e_1,e_3$. Since $B_2e_2=(1+h(A_2)_{22})e_2$ with $(A_2)_{22}\approx-1.3680097$, a vertex on the axis of $e_2$ would cap the level at $t\le1+h(A_2)_{22}$ (\cref{app:gm}). The $1836$-vertex polytope at $h=0.1$ avoids that axis and certifies $g(0.1)\ge-1.247834$. Transfer by \cref{prop:polytope}(ii), also checked directly at $h=0.05,0.025$, gives this bound for every $h\le0.1$. Yet $\one^\top e^{tA_2}e_2=e^{(A_2)_{22}t}$ implies $\lexp_{\min}\le(A_2)_{22}<-1.368009$. Thus $\lim g>\lexp_{\min}$ and both $\lim g$ and $\lexp$ lie in $[-1.247834,-1.241336]$ (\cref{tab:glp}); the adaptive lower-estimate stall does not reflect failure of Euler convergence.
\end{remark}

In contrast with case~148, every finite polytope in the next example is capped by the minimum face rate, although the Euler bound converges.

\begin{example}\label{ex:bcv}
Blanchini, Colaneri and Valcher~\cite{BCV13} consider the three modes
\[
 A_1=\begin{psmallmatrix}-1&0&0\\0&0&0\\0&1&\varepsilon\end{psmallmatrix},\qquad
 A_2=\begin{psmallmatrix}\varepsilon&0&1\\0&-1&0\\0&0&0\end{psmallmatrix},\qquad
 A_3=\begin{psmallmatrix}0&0&0\\1&\varepsilon&0\\0&0&-1\end{psmallmatrix},\qquad\varepsilon=\tfrac1{10}.
\]
No convex combination is Hurwitz ($\min_\alpha\eta(S_\alpha)=1/120$), yet a cyclic schedule stabilizes the system~\cite{BCV13}. All three modes are reducible.
\begin{itemize}[leftmargin=*]
\item For $x>0$, $a(x)=(x_1x_2x_3)^{1/3}$ satisfies $\frac{\mathrm d}{\mathrm dt}\log a(x(t))=\frac13(-1+\varepsilon+x_2/x_3)$ along $A_1$, and the cyclic analogues along $A_2$ and $A_3$. Hence $a(x(t))\ge e^{(\varepsilon-1)t/3}a(x(0))$ under every switching signal. Since $3a(x)\le\one^\top x$, \cref{prop:basic}(ii) gives 
\[
\lexp\ge(\varepsilon-1)/3.
\]
For the reverse inequality, use the cyclic schedule $A_3$, $A_2$, $A_1$, with equal dwell time $t>0$. Put $a_t=e^{-t}$, $b_t=e^{\varepsilon t}$ and $f_t=(b_t-1)/\varepsilon$. Direct multiplication gives
\[
 P_t=e^{tA_1}e^{tA_2}e^{tA_3}
 =\begin{pmatrix}
 a_tb_t&0&f_ta_t^2\\
 f_ta_t&a_tb_t&0\\
 f_t^2a_t&f_ta_tb_t&a_tb_t
 \end{pmatrix},\qquad
 \det(zI-P_t)=(z-a_tb_t)^3-f_t^3a_t^3z.
\]
Since $P_t\ge a_tb_tI$, its Perron root is $\rho(P_t)=a_tb_t(1+k_t)$ with $k_t\ge0$. The characteristic equation becomes
\[
 \frac{k_t^3}{1+k_t}=\frac{f_t^3a_t}{b_t^2}
 \le\frac{a_tb_t}{\varepsilon^3}\longrightarrow0.
\]
The function $k\mapsto k^3/(1+k)$ is strictly increasing on $[0,\infty)$ and vanishes only at zero, so $k_t\to0$. Repeating this schedule and using Gelfand's formula yields
\[
 \lexp\le\frac{\log\rho(P_t)}{3t}
 =\frac{\varepsilon-1}{3}+\frac{\log(1+k_t)}{3t}
 \longrightarrow\frac{\varepsilon-1}{3}.
\]
Thus $\lexp=-0.3$, approached by these schedules as their dwell times tend to infinity. Note that at dwell time $30$, the exponential enclosure and independent replay of \cref{sec:gripenberg} give the finite-schedule bound $\lexp\le-0.299986$.

\item Let $0<h<1$, $B_i=I+hA_i$ and $r_h=(1-h)(1+\varepsilon h)\in(0,1)$. The geometric mean satisfies $a(B_ix)^3\ge r_h a(x)^3$, and iteration from $x=\one$ gives $\lsr(\cB_h)\ge r_h^{1/3}$. For the reverse inequality, the block product $P_n=B_1^nB_2^nB_3^n$ has the form of $P_t$, with
\[
 a_n=(1-h)^n,\qquad b_n=(1+\varepsilon h)^n,\qquad f_n=(b_n-1)/\varepsilon.
\]
Its Perron root is $a_nb_n(1+k_n)$, where $k_n\ge0$ and $k_n^3/(1+k_n)\le r_h^n/\varepsilon^3\to0$. Consequently
\[
 \lsr(\cB_h)\le\rho(P_n)^{1/(3n)}
 =r_h^{1/3}(1+k_n)^{1/(3n)}\longrightarrow r_h^{1/3}.
\]
We have therefore proved the exact identities
\[
 \lsr(\cB_h)=r_h^{1/3},\qquad
 g(h)=\frac{r_h^{1/3}-1}{h}\longrightarrow\frac{\varepsilon-1}{3}=\lexp.
\]
The values of $g$ round to $-0.313030$ at $h=0.1$ and $-0.306335$ at $h=0.05$. The axes are invariant faces with exponent $-1$, each under two modes, and the two-dimensional faces have exponent $1/120$; hence $\lexp_{\min}=-1$ by \cref{thm:faceformula}.

\item Every finite $V$ that satisfies \eqref{eq:polycert} has $t\le1-h$, so \cref{prop:polytope} certifies at most $\lexp\ge-1$, the exponent of the axes, at every step. Indeed, suppose $t>1-h$. Then:
 \begin{itemize}
 \item Let $m_3=\min_{v\in V}v_3$. Since $(B_3x)_3=(1-h)x_3$, applying \eqref{eq:polycert} at a minimizer gives $(1-h)m_3\ge tm_3$, so $m_3=0$ and $V_0=\{v\in V:v_3=0\}$ is nonempty.
 \item Let $m_2=\min_{u\in V_0}u_2$. Since $(B_2x)_3=x_3$, the weights for a minimizer $u$ are supported on $V_0$. Since $(B_2x)_2=(1-h)x_2$, we get $(1-h)m_2\ge tm_2$, so $m_2=0$. As all vertices are nonzero, $V$ has a point on the positive axis of $e_1$.
 \item For the smallest of them, $me_1$, \eqref{eq:polycert} for $(1,me_1)$ reads $(1-h)me_1\ge t\sum_w\lambda_ww$. So the weights live on that axis, and $(1-h)m\ge tm$, a contradiction.
 \end{itemize}
 Applied to the families $\{A_i-\alpha I\}$ as in \cref{obs:glp}, the same argument caps the tangent bound of Algorithm~(L) at $-1$; the positive antinorms used by the Gripenberg algorithm~\cite{GM} are capped at $-1$ by the bound $\lsr_{\min}$ in \cref{app:gm}. In our runs at $h=0.1$, closures seeded by $\one$, by points of the two-dimensional faces, or by the coordinate vectors finish at the level $g=-1.01$, and none finishes at $-0.99$ within $120$\,s.
\end{itemize}
So finite polytope bounds do not approach the true exponent here, although the Euler bound converges. The optimal trajectories dwell near the axes, which are fast faces; the extremal antinorm, the geometric mean, vanishes on all of $\partial K$ and is not polyhedral.
\end{example}

\subsection{Numerical results}\label{subsec:numerical}

The experiments are designed to answer three questions: whether certified brackets behave as the rate theorems predict (\cref{obs:gripenberg}), how the polytope search scales with the dimension (\cref{obs:scaling}), and how the Euler approach compares with Algorithm~(L) of~\cite{GLP} (\cref{obs:glp}). 

\Cref{app:protocols,app:samples} record generators, budgets, hardware and search protocols. Adaptive Gripenberg lower estimates guide searches but require a posteriori verification: every test-pair estimate is verified, whereas seven of $54$ random-pair estimates remain unverified within the allocated budget (see \cref{app:gm}).

\begin{observation}\label{obs:gripenberg}
On test pairs, the bracket widths in \cref{tab:gripenberg} are approximately $1.48h$ for the planar pair (whose known exponent is enclosed) and $0.55h$ for case~1, matching $Qh$ with the coefficients $Q=1.478$ and $0.548$ of \cref{thm:rateconst}; case~1440's width halves from $h=0.1$ to $0.05$. Subsequent improvements for cases~1440 and~148 reuse polytopes from larger steps. Their certified upper bounds beat numerical mixture minima by up to $2.2\times10^{-3}$ and $3.1\times10^{-3}$, respectively; proving $\lexp<c$ still requires a global lower bound for $c$. For the fourteen pairs of \cref{app:tables}, finite-word slopes also approach the conjectured coefficients.

The rates cannot, however, be used directly to sharpen these bounds. For dimension-$5$ pair~4 of \cref{tab:switching}, the certified endpoints $L_h$ at target $10^{-5}$, rounded outward, satisfy $2L_{0.05}-L_{0.1}=-0.426998>-0.430082$, the certified schedule upper bound, so the extrapolated lower endpoint is not valid.
\end{observation}

\begin{table}[!htbp]
\centering\small
\setlength{\tabcolsep}{12pt}
\caption{Certified Euler brackets (\cref{rem:certified}); $c$ is the numerical mixture minimum, rounded up. $^\dagger$Lower bound obtained by reusing the polytope from $h=0.1$ (148) or $0.05$ (1440). Upper endpoints use implicit words; width is upper minus lower. The planar pair is that of \cref{obs:sandwich}.}
\label{tab:gripenberg}
\begin{tabular}{@{}lcccc@{}}
\toprule
pair ($c$) & $h$ & lower & upper & width\\
\midrule
planar & 0.1 & $-0.272230$ & $-0.123619$ & $1.5\times10^{-1}$\\
($-0.194976$) & 0.05 & $-0.232706$ & $-0.158706$ & $7.4\times10^{-2}$\\
 & 0.025 & $-0.213641$ & $-0.176676$ & $3.7\times10^{-2}$\\
 & 0.0125 & $-0.204261$ & $-0.185783$ & $1.8\times10^{-2}$\\
\midrule
case 1 & 0.1 & $-1.044234$ & $-0.988047$ & $5.6\times10^{-2}$\\
($-1.012395$) & 0.05 & $-1.027040$ & $-0.999505$ & $2.8\times10^{-2}$\\
 & 0.025 & $-1.019518$ & $-1.005757$ & $1.4\times10^{-2}$\\
\midrule
1440 & 0.1 & $\phantom{-}0.344598$ & $\phantom{-}0.346673$ & $2.1\times10^{-3}$\\
($0.348177$) & 0.05 & $\phantom{-}0.345146$ & $\phantom{-}0.346182$ & $1.0\times10^{-3}$\\
 & 0.025 & $\phantom{-}0.345146^\dagger$ & $\phantom{-}0.345931$ & $7.9\times10^{-4}$\\
\midrule
148 & 0.1 & $-1.247834$ & $-1.236958$ & $1.1\times10^{-2}$\\
($-1.236979$) & 0.05 & $-1.247834^\dagger$ & $-1.238961$ & $8.9\times10^{-3}$\\
 & 0.025 & $-1.247834^\dagger$ & $-1.240094$ & $7.7\times10^{-3}$\\
\bottomrule
\end{tabular}
\end{table}

\begin{observation}\label{obs:scaling}
The performance degrades quickly with the dimension. Indeed, at $h=0.1$, certified lower bounds are within $4.4\times10^{-3}$ of the word estimates for $18/20$ dimension-$5$ pairs of the scaling sample, but for only $1/20$ dimension-$10$ pairs, and dimension-$20$ widths are between $0.33$ and $1.6$.

The scaling sample has ten dense and ten sparse pairs in each of dimensions $5,10$, plus five dimension-$20$ pairs, with irreducible modes and interior numerical mixture minima. Chattering remains unproved; dimension-$10$ sparse~8 admits an improving schedule. Results are in \cref{tab:scaling,tab:scalingnew} and the repository records. Transfer improves every direct $h=0.05$ bound for the original six dimension-$5$ pairs, and mixture bounds improve on implicit Euler for the original twelve pairs (\cref{tab:glp}).

The switching sample (\cref{tab:switching} and \cref{tab:switchingnew}) contains $30$ pairs in each of dimensions $5,10$ and five in dimension~$20$, selected by schedules improving the numerical mixture minimum by $10^{-3}$; this does not prove $\lexp<c$. At $h=0.1$, a target $10^{-3}$ below the word estimate is reached for $30/30$ dimension-$5$ pairs and $14/30$ dimension-$10$ pairs, while dimension~$20$ closes only at target $0.1$. Target tolerance is measured below a word estimate, not relative to $\lexp$; outward rounding adds at most $10^{-6}$ per endpoint, and the bounds hold for both data readings. Finer steps increase search cost, and bracket widths also include schedule suboptimality. \Cref{tab:switching} and \cref{tab:switchingnew} report separate searches at each step; transfer preserves the best lower bound as $h$ decreases (\cref{prop:polytope}(ii)). Detailed success counts and widths are in \cref{app:samples}.
\end{observation}

\begin{table}[!htbp]
\centering\footnotesize
\setlength{\tabcolsep}{3pt}
\caption{Pairs selected by improving schedules (\cref{obs:scaling}). Estimates: $c$ (mixture) and $g(h)\le$ (explicit word), rounded up. Bounds: the schedule bound ``upper'' and $g(h)\ge$, the better lower bound from both orientations, rounded down; both hold for both data readings. $^\top$Transpose; $^+$extended search. Width is upper minus the larger lower bound, rounded up.}
\label{tab:switching}
\begin{tabular}{@{}rrccccccc@{}}
\toprule
& & & & \multicolumn{2}{c}{$h=0.1$} & \multicolumn{2}{c}{$h=0.05$} & \\
\cmidrule(lr){5-6}\cmidrule(lr){7-8}
$d$ & pair & $c$ & upper & $g(h)\le$ & $g(h)\ge$ & $g(h)\le$ & $g(h)\ge$ & width\\
\midrule
5 & 1 & $\phantom{-}0.539231$ & $\phantom{-}0.467683$ & $\phantom{-}0.465337$ & $\phantom{-}0.465326$ & $\phantom{-}0.466492$ & $\phantom{-}0.366491$ & $2.4\times10^{-3}$\\
5 & 2 & $\phantom{-}0.459280$ & $\phantom{-}0.457887$ & $\phantom{-}0.442429$ & $\phantom{-}0.442418$ & $\phantom{-}0.450720$ & $\phantom{-}0.440719$ & $1.6\times10^{-2}$\\
5 & 3 & $-0.333377$ & $-0.365833$ & $-0.457212$ & $-0.457223$ & $-0.409526$ & $-0.410527^\top$ & $4.5\times10^{-2}$\\
5 & 4 & $-0.401497$ & $-0.430082$ & $-0.460383$ & $-0.460394$ & $-0.443685$ & $-0.443696$ & $1.4\times10^{-2}$\\
5 & 5 & $\phantom{-}0.587104$ & $\phantom{-}0.548729$ & $\phantom{-}0.531162$ & $\phantom{-}0.531151$ & $\phantom{-}0.540134$ & $\phantom{-}0.539133$ & $9.6\times10^{-3}$\\
5 & 6 & $\phantom{-}0.532754$ & $\phantom{-}0.519922$ & $\phantom{-}0.500647$ & $\phantom{-}0.500636$ & $\phantom{-}0.511408$ & $\phantom{-}0.511397$ & $8.6\times10^{-3}$\\
5 & 7 & $\phantom{-}1.431824$ & $\phantom{-}1.428992$ & $\phantom{-}1.377996$ & $\phantom{-}1.377985$ & $\phantom{-}1.406985$ & $\phantom{-}1.405984^\top$ & $2.4\times10^{-2}$\\
5 & 8 & $-0.152645$ & $-0.171842$ & $-0.210279$ & $-0.211280$ & $-0.189358$ & $-0.189369$ & $1.8\times10^{-2}$\\
5 & 9 & $\phantom{-}0.481345$ & $\phantom{-}0.474726$ & $\phantom{-}0.404665$ & $\phantom{-}0.404654$ & $\phantom{-}0.440660$ & $\phantom{-}0.430659$ & $4.5\times10^{-2}$\\
5 & 10 & $\phantom{-}0.857180$ & $\phantom{-}0.843598$ & $\phantom{-}0.839946$ & $\phantom{-}0.838945$ & $\phantom{-}0.841758$ & $\phantom{-}0.840757$ & $2.9\times10^{-3}$\\
\midrule
10 & 1 & $\phantom{-}0.567910$ & $\phantom{-}0.559043$ & $\phantom{-}0.532613$ & $\phantom{-}0.522612^{+}$ & $\phantom{-}0.546054$ & $\phantom{-}0.446053$ & $3.7\times10^{-2}$\\
10 & 2 & $\phantom{-}0.716928$ & $\phantom{-}0.674041$ & $\phantom{-}0.629877$ & $\phantom{-}0.629866$ & $\phantom{-}0.653399$ & $\phantom{-}0.553398$ & $4.5\times10^{-2}$\\
10 & 3 & $\phantom{-}0.802565$ & $\phantom{-}0.788865$ & $\phantom{-}0.777934$ & $\phantom{-}0.727933^{+}$ & $\phantom{-}0.783449$ & $\phantom{-}0.683448$ & $6.1\times10^{-2}$\\
10 & 4 & $\phantom{-}0.253414$ & $\phantom{-}0.245147$ & $\phantom{-}0.216605$ & $\phantom{-}0.196604^{+}$ & $\phantom{-}0.231252$ & $\phantom{-}0.081251^\top$ & $4.9\times10^{-2}$\\
10 & 5 & $-0.168419$ & $-0.197966$ & $-0.228868$ & $-0.228879$ & $-0.212874$ & $-0.312875$ & $3.1\times10^{-2}$\\
10 & 6 & $\phantom{-}0.673707$ & $\phantom{-}0.638392$ & $\phantom{-}0.635617$ & $\phantom{-}0.625616^\top$ & $\phantom{-}0.636998$ & $\phantom{-}0.536997^\top$ & $1.3\times10^{-2}$\\
10 & 7 & $\phantom{-}1.200095$ & $\phantom{-}1.185699$ & $\phantom{-}1.136618$ & $\phantom{-}1.136607^\top$ & $\phantom{-}1.161028$ & $\phantom{-}1.061027^\top$ & $5.0\times10^{-2}$\\
10 & 8 & $\phantom{-}0.198883$ & $\phantom{-}0.142660$ & $\phantom{-}0.077000$ & $\phantom{-}0.056999^{+}$ & $\phantom{-}0.110582$ & $\phantom{-}0.010581^{+}$ & $8.6\times10^{-2}$\\
10 & 9 & $\phantom{-}0.953083$ & $\phantom{-}0.944065$ & $\phantom{-}0.932144$ & $\phantom{-}0.882143^{\top+}$ & $\phantom{-}0.938291$ & $\phantom{-}0.738290^\top$ & $6.2\times10^{-2}$\\
10 & 10 & $\phantom{-}0.670254$ & $\phantom{-}0.387498$ & $\phantom{-}0.369426$ & $\phantom{-}0.359425^{+}$ & $\phantom{-}0.378543$ & $\phantom{-}0.278542^{+}$ & $2.9\times10^{-2}$\\
\bottomrule
\end{tabular}
\end{table}

\subsection{Comparison with Algorithm~(L)}\label{subsec:comparison}

Algorithm~(L) of~\cite[Sec.~4.6]{GLP} certifies lower bounds by polytopes for the exponential family $\{e^{\tau A_i}\}$. We ran a comparison of our simplest implementation of it with the closure code and budgets of the Euler searches.

\begin{figure}[!htbp]
\centering
\includegraphics[width=\linewidth]{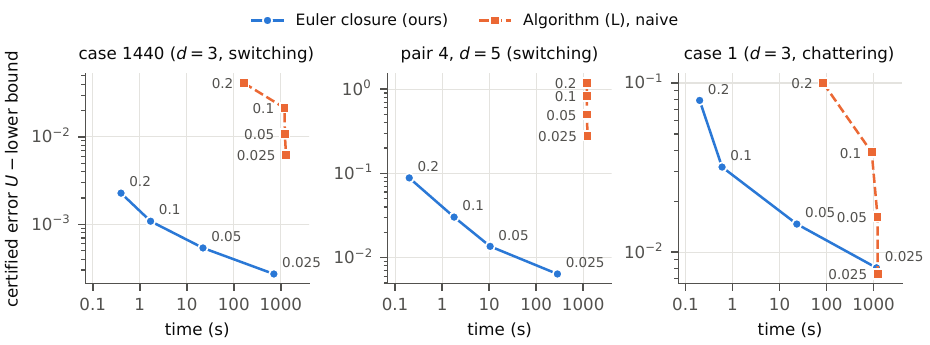}
\caption{Certified bracket width versus time (\cref{obs:glp}), log-log: Euler closure and our naive Algorithm~(L). Labels give $h$ or $\tau$; times below $0.1$\,s are drawn at $0.1$\,s.}
\label{fig:wp}
\end{figure}

\begin{observation}\label{obs:glp}
Neither the Euler approach nor Algorithm~(L) dominates on the test and scaling pairs. Under the reported implementations and budgets, Euler performs better on the switching benchmarks in \cref{fig:wp}; the best periods found can exceed the reach of short exponential words. At $\tau=0.1,0.05$, Algorithm~(L) gives the larger lower bound on $6/15$ pairs, versus Euler's $9/15$ (\cref{tab:glp}; reducible case~148 is excluded). Its tangent test improves neither the planar nor case~1 Euler polytope. Rigorous mixture or schedule bounds improve both methods' upper bounds.

In \cref{fig:wp}, Euler gaps decrease roughly linearly, with time rising from under one second to $280$--$1150$\,s. Algorithm~(L) is comparable for case~1 but worse for the switching pairs: pair~4's gap exceeds $0.27$, since length-$14$ words at $\tau=0.025$ cover periods only up to $0.35$, versus the best period near $3.3$. The figure uses fixed certified upper bounds and includes searches, failures and checks in the timings.

The lower bound of Algorithm~(L) is the largest $\alpha$ such that $(A_i-\alpha I)w$ points into $\mathcal P=\operatorname{conv}V+K$ at every vertex $w$, for every $i$. This is the bound $\check\alpha(\mathcal P)$ of~\cite[Prop.~8]{GLP}. It is \cref{prop:polytope} in the limit $h\to0$: if the condition holds, then $(I+s(A_i-\alpha I))\mathcal P\subseteq\mathcal P$ for small $s>0$, and \cref{prop:polytope}(i), applied to the family $\{A_i-\alpha I\}$ with $t=1$, gives $\lexp\ge\alpha$. 

We compute $\alpha$ by linear programming and then exactly, from rationalized multipliers, for both data readings. These lower bounds are certified but not independently replayed.

Closure code, budgets and hardware are shared, but our extra long-word search and naive implementation prevent attributing differences solely to discretization. \Cref{subsec:algorithmL,subsec:workprecision} specify the exponential-word search, the polytope closure and the timing protocol.
\end{observation}

\begin{table}[!htbp]
\centering\small
\setlength{\tabcolsep}{2pt}
\caption{Best bounds over steps $0.1,0.05$ (\cref{obs:glp}), with the larger lower bound in bold. Algorithm~(L) lower bounds are not independently replayed, and its upper values are estimates. $^g$Closure seeded with adaptive Gripenberg vertices at $h=0.1$; our search alone gives $0.254546,1.312394$. ``$c$ or schedule'' is a mixture or $^s$schedule bound. Width uses the larger lower bound and last upper column; upper values and widths are rounded up. Reducible case~148 was not run with Algorithm~(L).}
\label{tab:glp}
\begin{tabular}{@{}lcccccc@{}}
\toprule
& \multicolumn{2}{c}{lower bound (certified)} & \multicolumn{3}{c}{upper bound} &\\
\cmidrule(lr){2-3}\cmidrule(lr){4-6}
pair & Alg.~(L) & ours & Alg.~(L) & \shortstack{implicit\\Euler} & \shortstack{$c$ or\\schedule} & width\\
\midrule
planar & $\mathbf{-0.210758}$ & $-0.232706$ & $-0.194213$ & $-0.158706$ & $-0.194976$ & $1.6\times10^{-2}$\\
case 1 & $-1.028582$ & $\mathbf{-1.027040}$ & $-1.011626$ & $-0.999505$ & $-1.012395$ & $1.5\times10^{-2}$\\
case 1440 & $\phantom{-}0.334850$ & $\mathbf{\phantom{-}0.345146}$ & $\phantom{-}0.346401$ & $\phantom{-}0.346182$ & $\phantom{-}0.345675$$^{s}$ & $5.3\times10^{-4}$\\
case 148 & --- & $-1.247834$ & --- & $-1.238961$ & $-1.241336$$^{s}$ & $6.5\times10^{-3}$\\
\midrule
$d=5$, dense 1 & $-0.276163$ & $\mathbf{-0.222751}$ & $-0.213075$ & $-0.209005$ & $-0.213383$ & $9.4\times10^{-3}$\\
$d=5$, dense 2 & $-0.081527$ & $\mathbf{-0.059850}$ & $-0.028459$ & $-0.014798$ & $-0.029331$ & $3.1\times10^{-2}$\\
$d=5$, dense 3 & $\phantom{-}0.210187$ & $\mathbf{\phantom{-}0.222107}$ & $\phantom{-}0.233137$ & $\phantom{-}0.237887$ & $\phantom{-}0.232972$ & $1.1\times10^{-2}$\\
$d=5$, sparse 1 & $-0.204184$ & $\mathbf{-0.107049}$ & $-0.032087$ & $-0.001609$ & $-0.034259$ & $7.3\times10^{-2}$\\
$d=5$, sparse 2 & $\mathbf{-0.384238}$ & $-0.641285$ & $-0.102190$ & $-0.097694$ & $-0.102660$ & $2.9\times10^{-1}$\\
$d=5$, sparse 3 & $\mathbf{-0.049779}$ & $-0.071028$ & $\phantom{-}0.055587$ & $\phantom{-}0.112334$ & $\phantom{-}0.053208$ & $1.1\times10^{-1}$\\
\midrule
$d=10$, dense 1 & $\mathbf{\phantom{-}0.870186}$ & $\phantom{-}0.459247$ & $\phantom{-}0.990484$ & $\phantom{-}0.991569$ & $\phantom{-}0.990410$ & $1.3\times10^{-1}$\\
$d=10$, dense 2 & $\mathbf{\phantom{-}0.779913}$ & $\phantom{-}0.525978$ & $\phantom{-}0.857197$ & $\phantom{-}0.865203$ & $\phantom{-}0.856999$ & $7.8\times10^{-2}$\\
$d=10$, dense 3 & $\mathbf{\phantom{-}0.009496}$ & $-0.229448$ & $\phantom{-}0.066898$ & $\phantom{-}0.068764$ & $\phantom{-}0.066852$ & $5.8\times10^{-2}$\\
$d=10$, sparse 1 & $-0.064519$ & $\mathbf{\phantom{-}0.269192}^g$ & $\phantom{-}0.318964$ & $\phantom{-}0.334496$ & $\phantom{-}0.317685$ & $4.9\times10^{-2}$\\
$d=10$, sparse 2 & $\phantom{-}0.364290$ & $\mathbf{\phantom{-}1.356787}^g$ & $\phantom{-}1.570804$ & $\phantom{-}1.583464$ & $\phantom{-}1.570368$ & $2.2\times10^{-1}$\\
$d=10$, sparse 3 & $-0.451797$ & $\mathbf{-0.349697}$ & $-0.179980$ & $-0.143149$ & $-0.181089$ & $1.7\times10^{-1}$\\
\bottomrule
\end{tabular}
\end{table}

\FloatBarrier
\section{Open problems}\label{sec:open}

We conclude by summarizing some open questions that arise naturally from this paper. 

\begin{enumerate}[label=(\arabic*)]
\item \textbf{General explicit convergence.} Does $g(h)\to\lexp$ as $h\downarrow0$ for every finite family of Metzler matrices? The sufficient cone criterion of \cref{cor:eulerface} settles the five test pairs; its failure leaves convergence undecided.

\item \textbf{Sharp first-order asymptotics.} Under \cref{thm:rateconst}, with $\lexp=c$ and $\eta''(a^\ast)>0$, where $\eta(a)=\eta(S_a)$, do $g(h)=c-Qh/2+O(h^2)$ and $G(h)=c+Qh/2+O(h^2)$ hold? Matching lower estimates are missing. For three or more modes, identifying optimal coefficients must also address word ordering (\cref{app:rateproof}).

\item \textbf{Gripenberg extensions.} Can the fixed-antinorm result of \cref{cor:gripenberg} be extended to reducible families? What conditions ensure validity and convergence of the adaptive lower estimates (\cref{rem:adaptive})?

\item \textbf{Computational cost.} Which structural assumptions yield bounds on vertex counts and computational cost for a prescribed dimension and certified bracket width? The obstruction in \cref{ex:bcv} motivates nonpolyhedral antinorms. Combining the constructions of \cref{obs:glp} is a separate practical direction; it cannot remove a shared face cap.

\item \textbf{Other discretizations and cones.} Can cone-preserving higher-order discretizations yield rigorous two-sided bounds with smaller discretization gaps? On proper nonpolyhedral cones, flow invariance need not imply invariance under $I+hA$.
\end{enumerate}

\section*{Acknowledgments}
The author is a member of INdAM-GNCS (Gruppo Nazionale di Calcolo Scientifico).

\section*{Statements and declarations}
\subsection*{Data and code availability}

The datasets generated and analysed during this study, together with the computational code, are available in the \href{https://github.com/francescomaiale/certified-lower-exponents}{\texttt{certified-lower-exponents} repository}, at verified commit \href{https://github.com/francescomaiale/certified-lower-exponents/tree/a10a464e1f7a9272426fff289db4448cb9408b88}{\texttt{a10a464e1f7a}}. The repository contains the input matrices, complete numerical records, stored certificates, and independent Python and MATLAB verification scripts.

The repository also contains the original Gripenberg run records and final vertex sets. The repository is publicly accessible. The main code uses the MIT License and the optional Gripenberg implementation retains GPL-3.0.

Detailed Gripenberg runs appear in \cref{app:gm}; \cref{app:protocols,app:samples} give the experimental settings, verification conventions, and scaling and switching samples.

Commands and driver options are in \nolinkurl{reproduce/README.md}, with instructions for the Gripenberg runs in \nolinkurl{reproduce/gm/README.md}. Detailed results are in \nolinkurl{reproduce/results/}, including \nolinkurl{gm_original/} and \nolinkurl{larger_samples/}. The verification run of 30 October 2026 passed all 14 Python regression tests and replayed all 568 supplied certificates. Replays check saved witnesses; rerunning searches can change timings and final estimates.

Algorithm~(L)'s lower bounds are checked exactly when generated but lack independent replay witnesses; its upper bounds and the word and mixture estimates remain floating point.

\subsection*{Author contributions}
F.P.M. conceived the study, developed the mathematical results and proofs, developed and validated the software, carried out the numerical experiments, prepared the figures and tables, and wrote and revised the manuscript. F.P.M. approved the final manuscript and takes responsibility for all aspects of the work.

\subsection*{Funding}
No funding was received for this research.

\subsection*{Competing interests}
The author has no competing interests to declare.

\subsection*{AI assistance}
Codex (OpenAI) assisted with numerical experiments, code development, and presentation\slash grammar editing. The author takes full responsibility for the article; in particular, verified bounds rest on the exact checks of \cref{sec:gripenberg}.

\appendix
\setcounter{equation}{0}
\renewcommand{\theequation}{\Alph{section}\arabic{equation}}
\renewcommand{\theHequation}{appendix.\thesection.\arabic{equation}}

\section{Data of the test pairs}\label[appendix]{app:data}

\begin{sloppypar}
Cases~77, 148, 677, 1195 and~1440 retain the identifiers used in the repository. Case~1 is an additional test pair, specified to six decimal places. For the other cases, each entry is printed with the fewest digits needed to recover its stored double-precision (binary64) value when read back.

The bounds in \cref{sec:gripenberg} are valid whether the printed decimals are treated as exact numbers or converted to binary64. The computations in \cref{obs:cone} use the latter interpretation.

The planar pairs are given in \cref{ex:sharp,obs:sandwich}.
\end{sloppypar}

\noindent\textbf{Case 1}:
\begin{align*}
A_1&=\begin{psmallmatrix}-1.730493&0.996436&0.459716\\0.69104&-2.643051&0.03405\\0.84589&0.587882&-2.309322\end{psmallmatrix},\\*
A_2&=\begin{psmallmatrix}-2.892226&0.839125&0.466303\\0.127203&0.220544&0.195653\\0.06192&0.598392&-2.23932\end{psmallmatrix}.
\end{align*}

\noindent\textbf{Case 77}:
\begin{align*}
A_1&=\begin{psmallmatrix}-0.10421255427821108&0&0.5995296020142661\\0&-3.983676417656142&10.44547437798186\\0.6089541200819454&0.41616545550352757&-3.2474691918642016\end{psmallmatrix},\\*
A_2&=\begin{psmallmatrix}-0.6965629728040771&0&0.23901622681524412\\0&0.02489242262644542&0\\0.18904586288556002&0.8602987239031196&-2.0161017416771956\end{psmallmatrix}.
\end{align*}

\noindent\textbf{Case 148}:
\begin{align*}
A_1&=\begin{psmallmatrix}-2.854194178918802&0.24749928971794594&0\\0.56180752192356&-1.3059431283783423&0.20794658684868464\\0&0&-3.669288482126051\end{psmallmatrix},\\*
A_2&=\begin{psmallmatrix}-3.653224173747234&0&0.49425460335236504\\2.342863144895686&-1.3680096637401422&0\\0.7499705168800616&0&0.05001607999557478\end{psmallmatrix}.
\end{align*}

\noindent\textbf{Case 677}:
\begin{align*}
A_1&=\begin{psmallmatrix}-3.5751753658589895&0&0.7873089953050745\\3.872088293234137&-2.4704665923380382&3.292200232648817\\1.6772091050048494&0.4874274815499069&-1.7576243608193685\end{psmallmatrix},\\*
A_2&=\begin{psmallmatrix}-3.608486201678475&0&4.290473593750214\\0&-0.23709653436107514&4.089406788141199\\0.7796352280553934&0&-3.3591039935935543\end{psmallmatrix}.
\end{align*}

\noindent\textbf{Case 1195}:
\begin{align*}
A_1&=\begin{psmallmatrix}-3.2011480664004672&0&0\\0&-1.601165653605288&1.1392361293277327\\0.11566889906380853&2.7419424885042747&-2.670478064763273\end{psmallmatrix},\\*
A_2&=\begin{psmallmatrix}0.09930018711201494&0&0.271862921876863\\0.831655604528253&-2.3694492859575673&2.394339449701255\\5.327089103180349&1.129064780325959&-3.977234504357898\end{psmallmatrix}.
\end{align*}

\noindent\textbf{Case 1440}:
\begin{align*}
A_1&=\begin{psmallmatrix}-2.163042579047693&0.33559573959929034&0.17548453245292883\\0.7960011589145731&-1.0721730232799498&0.7689150093877756\\0.3079265062822859&0.7102921800837112&-0.11334509645810664\end{psmallmatrix},\\*
A_2&=\begin{psmallmatrix}-0.365512369671777&0.895993665091137&0.051482815948108995\\0.24535640370572198&-0.3430416586652343&0.3203492321022128\\0.05334247416806037&0.7483745160181475&-0.29764076048251553\end{psmallmatrix}.
\end{align*}

\section{Runs of the Gripenberg algorithms}\label[appendix]{app:gm}

\Cref{cor:gripenberg} guarantees fixed-antinorm convergence. \cref{tab:gmruns,tab:gmscaling,tab:gmswitching} report adaptive runs, whose lower estimates require a posteriori verification (\cref{rem:adaptive}).

\paragraph{Normalization and accuracy} Runs use transposed families and $\delta=h^2/10$, a discrete-gap tolerance rather than a continuous-time width; precise settings are in \cref{subsec:gmsettings}.

For a fixed antinorm, scaling $(\mathcal F,\delta)$ to $(\gamma\mathcal F,\gamma\delta)$ scales both bounds by $\gamma>0$. Adaptive acceptance requires $a(z)\le1$ and presumes $\lsr(\mathcal F)\approx1$~\cite[Sec.~4.2]{GM}; we normalize as $e^{-h\hat\lambda}\cB_h$ and $(1-h\hat\lambda)\cR_h$, then rescale. Without normalization, case~1440's lower estimate remained nonpositive at $h=0.1,0.05$, although $\lexp\approx0.3457$. Transposition preserves primitivity, $\lsr$ and, under \cref{cor:gripenberg}, the cone condition. On $\gamma\mathcal F$, tolerance $\delta$ becomes $\delta/\gamma$ in that corollary, contributing $\delta/(\gamma h)$ to the explicit-rate error; without normalization this corresponds to $\delta/h=h/10$.

The code's inherited joint-spectral-radius error correction inflates antinorms by order $10^{-15}$ here. This cannot validate lower bounds as exact polytope checks do.

\begin{observation}[Gripenberg runs]\label{obs:gmruns}
Adaptive runs often exhaust their budgets at finer steps. All four fixed-antinorm explicit runs do so with weaker verified lower bounds, whereas their implicit estimates are within $2\times10^{-4}$ of \cref{tab:gripenberg}. \cref{tab:gmruns,tab:gmscaling,tab:gmswitching} give the counts and seeded-closure improvements.
\end{observation}

\paragraph{The obstruction in case~148} Both modes are reducible. Since the third row of $A_1$ is $(0,0,\nu)$, with $\nu=(A_1)_{33}\approx-3.6692885<\lexp$ and $\eta(A_1)\approx-1.220815$, $A_1^\top e_3=\nu e_3$ is a secondary eigenvector in $K$, violating the cone condition~\cite{GM}. For $P=(I+hA_1)^\top$, $\mu=1+h\nu$ and every positive antinorm $a$, including adaptive ones,
\[
 a(P^k)\le\frac{a(P^ke_3)}{a(e_3)}=\mu^k,
\]
so $q(P^k)\le\mu$ in~\cite{GM}. Since $\mu+\delta$ lies below every upper bound reached, $P^k$ is never discarded in exact arithmetic: $t_k\le\mu$ and $(t_k-1)/h\le\nu$. The observed plateau is $-3.6692885$, as in the critical example of~\cite[Sec.~5.2]{GM}; upper bounds remain valid, but $s_k-t_k>5\times10^{-2}$ at every reported step.

More generally, with a fixed positive antinorm the Gripenberg lower bounds never exceed
\[
 \lsr_{\min}(\mathcal F)=\lim_{n\to\infty}\min_{|w|=n}\min_j\bigl(\one^\top M_we_j\bigr)^{1/n},
\]
the discrete counterpart of $\lexp_{\min}$ and of \cref{lem:bounded}(ii). Indeed, the proof of the lower-bound lemma of~\cite[Sec.~3.1]{GM} gives $t_k\le\liminf_N\alpha_N(\mathcal F)$, where $\alpha_N(\mathcal F)=\min_{|w|=N}a(M_w)^{1/N}$, without the cone condition; and $c_1\one^\top x\le a(x)\le c_2\one^\top x$ on $K$ gives
\[
 \alpha_N(\mathcal F)\le(c_2/c_1)^{1/N}\min_{|w|=N}\min_j\bigl(\one^\top M_we_j\bigr)^{1/N}.
\]
Without transposition, $A_2e_2=(A_2)_{22}e_2$ gives the better cap $(A_2)_{22}=\lexp_{\min}\approx-1.368010$ (\cref{obs:cone}). Transposition preserves the lower spectral radius, so one can retain the better orientation. Even then, positive antinorms leave width at least $\lexp-(A_2)_{22}\ge0.12$, since \cref{rem:certified} gives $\lexp\ge-1.247834$. Its polytopes avoid the axis of $e_2$; their antinorms vanish there and escape this cap.

\FloatBarrier
\subsection{Implementation, settings and detailed results}\label{subsec:gmsettings}
The following results complement \cref{obs:gmruns}; verification procedures are given in \cref{app:protocols}.

All Gripenberg runs use the updated \href{https://github.com/francescomaiale/subradius_computation}{\texttt{subradius\_computation} repository} of Guglielmi and Maiale, fixed at commit \texttt{181d6b5}. Runs use transposed families, coordinate seeds, accuracy $\delta=h^2/10$ and a budget of $2\times10^4$ products, which can be exceeded while completing a level. Time caps are $7200$\,s for test pairs and $1800$\,s for random pairs. The fixed-antinorm runs discussed in \cref{obs:gmruns} use $a(x)=\one^\top x$ and $3600$\,s. Time caps are checked between products; twelve runs ran concurrently. Accuracy $\delta$ means a reported discrete gap for the supplied, possibly normalized family, not a certified continuous-time bracket width. The explicit and implicit families are normalized as $e^{-h\hat\lambda}\cB_h$ and $(1-h\hat\lambda)\cR_h$ and then rescaled. We use $\hat\lambda=0.345674608$ for case 1440, no normalization for the other test pairs, the mixture value for the scaling sample and the certified schedule bound for the switching sample. The preceding analysis explains the normalization, rounding caveat and face obstruction in case 148.

\emph{Test pairs} (\cref{tab:gmruns}). The planar pair and cases 1 and 1440 reach accuracy $\delta$ at $h=0.1$ and $0.05$. At finer steps the explicit runs exhaust their product or time budgets and fall well below the certified bounds of \cref{tab:gripenberg}. Every reported lower estimate is certified a posteriori by \cref{prop:polytope}, with $3$ to $225$ vertices.

\emph{Random pairs} (\cref{tab:gmscaling}). At $h=0.1$, four of six dimension-$5$ explicit runs reach accuracy $\delta$; none does at $h=0.05$. Only one of the $24$ dimension-$10$ runs, explicit and implicit, reaches it. Nine of the other eleven explicit runs discard no products, exploring complete trees of depth $9$ or $10$. Our lower bounds exceed the estimates in $18$ of $24$ rows. Seeding with the final vertices of the adaptive Gripenberg runs certifies two further estimates, $0.269192$ and $1.356787$, used in \cref{tab:glp}; the remaining four higher estimates are uncertified.

\emph{Switching pairs} (\cref{tab:gmswitching}). Seven of ten dimension-$5$ explicit runs reach accuracy $\delta$ at $h=0.1$; all others time out below our certified bounds. All implicit runs time out above the schedule bounds.

\emph{Fixed antinorm.} All four explicit runs exhaust the product budget before accuracy $\delta$; their certified estimates lie below the adaptive and our bounds. Implicit estimates are within $2\times10^{-4}$ of \cref{tab:gripenberg}. Slow fixed-antinorm convergence is also discussed in \cite[Sec.~5.1]{GM}.

\paragraph{Certification of the outputs} A final adaptive vertex set can be checked directly by \cref{prop:polytope}; pruning and rescaling need not preserve invariance. Alternatively, it can seed a new closure. Seeded closures certify the estimate, rounded down to six decimals, in $45/47$ explicit test, switching and fixed-antinorm runs. Failures within $600$\,s are case~148 and dimension-$10$ pair~3 at $h=0.1$; none improves \cref{tab:gripenberg} and \cref{tab:switching}. Failure is inconclusive since separate closures certify every test-pair estimate, as reported above.

\begingroup
\small
\setlength{\tabcolsep}{4pt}
\begin{longtable}{@{}lccccc@{}}
\caption{Adaptive Gripenberg runs (\cref{obs:gmruns}). Lower: explicit estimate; upper: best implicit product, both converted to continuous-time rates. Depth gives the two maximum word lengths; time is their sum. $^\ast$Product budget and $^\dagger$time limit reached before accuracy $\delta$.}\label{tab:gmruns}\\
\toprule
pair & $h$ & lower & upper & depth & time (s)\\
\midrule
\endfirsthead
\caption[]{Continued from the previous page.}\\
\toprule
pair & $h$ & lower & upper & depth & time (s)\\
\midrule
\endhead
\bottomrule
\endfoot
\bottomrule
\endlastfoot
planar & 0.1 & $-0.281832$ & $-0.121884$ & 11/10 & 12\\*
 & 0.05 & $-0.237069$ & $-0.157436$ & 13/25 & 43\\*
 & 0.025 & $-0.238368^\ast$ & $-0.176426$ & 20/19 & 4445\\*
 & 0.0125 & $-1.019681^\ast$ & $-0.185482^\ast$ & 15/15 & 6189\\
\midrule
case 1 & 0.1 & $-1.054233$ & $-0.987914$ & 9/10 & 14\\*
 & 0.05 & $-1.031873$ & $-0.999444$ & 16/33 & 215\\*
 & 0.025 & $-1.100284^\dagger$ & $-1.005758^\dagger$ & 14/14 & 16289\\
\midrule
1440 & 0.1 & $\phantom{-}0.336905$ & $\phantom{-}0.347799$ & 11/12 & 36\\*
 & 0.05 & $\phantom{-}0.342769$ & $\phantom{-}0.347957$ & 18/19 & 3235\\*
 & 0.025 & $\phantom{-}0.116748^\dagger$ & $\phantom{-}0.348080^\dagger$ & 13/13 & 14400\\
\midrule
148 & 0.1 & $-3.669288^\ast$ & $-1.235611^\ast$ & 14/15 & 1426\\*
 & 0.05 & $-3.669288^\ast$ & $-1.235823^\ast$ & 14/14 & 1376\\*
 & 0.025 & $-3.669288^\ast$ & $-1.235803^\ast$ & 14/14 & 258\\
\end{longtable}
\endgroup

\begingroup
\small
\setlength{\tabcolsep}{4pt}
\begin{longtable}{@{}llcccc@{}}
\caption{Random-pair explicit estimates (\cref{obs:gmruns}), floating point; $c$ is the numerical mixture minimum; $^\ast$time limit before accuracy $\delta$. ``Ours'' certifies $g(h)$; $^\dagger$transferred from $h=0.1$. The last column certifies the adaptive Gripenberg estimate by a closure seeded with its final vertices; --- means no closure within $300$\,s.}\label{tab:gmscaling}\\
\toprule
$d$ & pair ($c$) & $h$ & estimate & ours & seeded bound\\
\midrule
\endfirsthead
\caption[]{Continued from the previous page.}\\
\toprule
$d$ & pair ($c$) & $h$ & estimate & ours & seeded bound\\
\midrule
\endhead
\bottomrule
\endfoot
\bottomrule
\endlastfoot
5 & dense 1 ($-0.213383$) & 0.1 & $-0.232425$ & $-0.222751$ & $-0.232425$\\*
 &  & 0.05 & $-0.292711^\ast$ & $-0.222751^\dagger$ & $-0.292712$\\
5 & dense 2 ($-0.029331$) & 0.1 & $-0.069430$ & $-0.059850$ & $-0.069431$\\*
 &  & 0.05 & $-0.138457^\ast$ & $-0.059850^\dagger$ & $-0.138458$\\
5 & dense 3 ($+0.232972$) & 0.1 & $\phantom{-}0.194335^\ast$ & $\phantom{-}0.222107$ & $\phantom{-}0.194334$\\*
 &  & 0.05 & $\phantom{-}0.057287^\ast$ & $\phantom{-}0.222107^\dagger$ & ---\\
5 & sparse 1 ($-0.034259$) & 0.1 & $-0.112654$ & $-0.107049$ & $-0.112655$\\*
 &  & 0.05 & $-0.814649^\ast$ & $-0.107049^\dagger$ & $-0.814650$\\
5 & sparse 2 ($-0.102660$) & 0.1 & $-0.798666^\ast$ & $-0.641285$ & $-0.798667$\\*
 &  & 0.05 & $-1.014139^\ast$ & $-0.641285^\dagger$ & $-1.014139$\\
5 & sparse 3 ($+0.053208$) & 0.1 & $-0.080404$ & $-0.071028$ & $-0.080404$\\*
 &  & 0.05 & $-0.076264^\ast$ & $-0.071028^\dagger$ & $-0.076264$\\
\midrule
10 & dense 1 ($+0.990410$) & 0.1 & $\phantom{-}0.533412^\ast$ & $\phantom{-}0.459247$ & ---\\*
 &  & 0.05 & $\phantom{-}0.046459^\ast$ & $\phantom{-}0.459247^\dagger$ & ---\\
10 & dense 2 ($+0.856999$) & 0.1 & $\phantom{-}0.783322^\ast$ & $\phantom{-}0.525978$ & ---\\*
 &  & 0.05 & $\phantom{-}0.541437^\ast$ & $\phantom{-}0.525978^\dagger$ & ---\\
10 & dense 3 ($+0.066852$) & 0.1 & $-0.164566^\ast$ & $-0.229448$ & ---\\*
 &  & 0.05 & $-0.612417^\ast$ & $-0.229448^\dagger$ & $-0.612417$\\
10 & sparse 1 ($+0.317685$) & 0.1 & $\phantom{-}0.269193$ & $\phantom{-}0.254546$ & $\phantom{-}0.269192$\\*
 &  & 0.05 & $-1.092236^\ast$ & $\phantom{-}0.254546^\dagger$ & $-1.092236$\\
10 & sparse 2 ($+1.570368$) & 0.1 & $\phantom{-}1.356787^\ast$ & $\phantom{-}1.312394$ & $\phantom{-}1.356787$\\*
 &  & 0.05 & $\phantom{-}0.762998^\ast$ & $\phantom{-}1.312394^\dagger$ & $\phantom{-}0.762998$\\
10 & sparse 3 ($-0.181089$) & 0.1 & $-0.917161^\ast$ & $-0.349697$ & $-0.917161$\\*
 &  & 0.05 & $-1.723107^\ast$ & $-0.349697^\dagger$ & $-1.723108$\\
\end{longtable}
\endgroup

\begingroup
\small
\setlength{\tabcolsep}{4pt}
\begin{longtable}{@{}rcrcrcr@{}}
\caption{Explicit estimates from the adaptive Gripenberg algorithm for \cref{tab:switching}, converted to $g(h)$, and run times (floating point). $^\ast$Time limit before accuracy $\delta$; dimension~$10$ was not run at $h=0.05$.}\label{tab:gmswitching}\\
\toprule
& \multicolumn{2}{c}{$d=5$, $h=0.1$} & \multicolumn{2}{c}{$d=5$, $h=0.05$} & \multicolumn{2}{c}{$d=10$, $h=0.1$}\\
\cmidrule(lr){2-3}\cmidrule(lr){4-5}\cmidrule(lr){6-7}
pair & estimate & time (s) & estimate & time (s) & estimate & time (s)\\
\midrule
\endfirsthead
\caption[]{Continued from the previous page.}\\
\toprule
& \multicolumn{2}{c}{$d=5$, $h=0.1$} & \multicolumn{2}{c}{$d=5$, $h=0.05$} & \multicolumn{2}{c}{$d=10$, $h=0.1$}\\
\cmidrule(lr){2-3}\cmidrule(lr){4-5}\cmidrule(lr){6-7}
pair & estimate & time (s) & estimate & time (s) & estimate & time (s)\\
\midrule
\endhead
\bottomrule
\endfoot
\bottomrule
\endlastfoot
1 & $\phantom{-}0.114886^\ast$ & 1800 & $\phantom{-}0.080790^\ast$ & 2703 & $-0.391747^\ast$ & 3544\\*
2 & $\phantom{-}0.431960$ & 103 & $\phantom{-}0.021179^\ast$ & 1801 & $-1.797610^\ast$ & 1800\\*
3 & $-0.466853$ & 778 & $-1.057132^\ast$ & 1801 & $\phantom{-}0.151578^\ast$ & 1800\\
4 & $-0.463750$ & 421 & $-0.971803^\ast$ & 4148 & $-0.876970^\ast$ & 5116\\
5 & $\phantom{-}0.520597$ & 288 & $\phantom{-}0.400567^\ast$ & 2252 & $-1.331466^\ast$ & 1801\\
6 & $\phantom{-}0.419929^\ast$ & 2010 & $-0.459780^\ast$ & 1800 & $-0.766025^\ast$ & 1800\\
7 & $\phantom{-}1.366459$ & 72 & $\phantom{-}0.628847^\ast$ & 3368 & $\phantom{-}1.079967^\ast$ & 2511\\
8 & $-0.218978$ & 760 & $-1.557098^\ast$ & 1801 & $-0.488982^\ast$ & 1800\\*
9 & $\phantom{-}0.394178$ & 1156 & $-0.450092^\ast$ & 1801 & $\phantom{-}0.504566^\ast$ & 3581\\*
10 & $\phantom{-}0.766664^\ast$ & 1800 & $-0.036211^\ast$ & 1826 & $-0.018742^\ast$ & 1801\\
\end{longtable}
\endgroup

\FloatBarrier
\section{First-order coefficients}\label[appendix]{app:tables}

For the planar pair, case~1 and the original twelve pairs of \cref{obs:scaling}, use $a^\ast,c,Q$ and Perron vectors normalized by $u^\top v=1$ as in \cref{thm:rateconst}. The finite-word estimates $g_{\rm est},G_{\rm est}$ define
\[
 s_-:=\frac{c-g_{\rm est}}h,\qquad
 s_+:=\frac{G_{\rm est}-c}h,\qquad
 s_-+s_+=\frac{G_{\rm est}-g_{\rm est}}h.
\]
Optimal chattering and the conjectured expansions of \cref{sec:open} would give limits $Q/2$ for the exact slopes $(c-g(h))/h$, $(G(h)-c)/h$, and $Q$ for their sum. Finite-word slopes require both search errors to be $o(h)$ to share these limits. Cancellation can conceal those errors, so the table also reports separate slopes at $h=0.025$. Over the fourteen pairs, the maximum relative deviations of $s_-+s_+$ from $Q$ are $2.52\%$, $0.64\%$ and $0.31\%$ at $h=0.1,0.05,0.025$, and at the finest step the midpoint $(g_{\rm est}+G_{\rm est})/2$ differs from $c$ by at most $0.83h^2$. \Cref{thm:rateconst} proves only one-sided bounds; these diagnostics support, but prove neither, matching expansions and optimal chattering, and better words remain possible.

\begingroup
\footnotesize
\setlength{\tabcolsep}{3pt}
\begin{longtable}{@{}lccccccc@{}}
\caption{Fourteen-pair rate diagnostics (\cref{obs:gripenberg}). Word estimates bound the discrete minima from above up to rounding; slopes use the unrounded repository records.}\label{tab:rates}\\
\toprule
& & & \multicolumn{3}{c}{$s_-+s_+$ (target $Q$)} & \multicolumn{2}{c}{\shortstack{target $Q/2$\\$h=0.025$}}\\
\cmidrule(lr){4-6}\cmidrule(l){7-8}
pair & $a^\ast$ & $Q$ & $h=0.1$ & $h=0.05$ & $h=0.025$ & $s_-$ & $s_+$\\
\midrule
\endfirsthead
\caption[]{Continued from the previous page.}\\
\toprule
& & & \multicolumn{3}{c}{$s_-+s_+$ (target $Q$)} & \multicolumn{2}{c}{\shortstack{target $Q/2$\\$h=0.025$}}\\
\cmidrule(lr){4-6}\cmidrule(l){7-8}
pair & $a^\ast$ & $Q$ & $h=0.1$ & $h=0.05$ & $h=0.025$ & $s_-$ & $s_+$\\
\midrule
\endhead
\bottomrule
\endfoot
\bottomrule
\endlastfoot
planar & $0.3905$ & $1.4780$ & $1.4861$ & $1.4800$ & $1.4785$ & $0.7466$ & $0.7320$\\
case 1 & $0.7890$ & $0.5480$ & $0.5619$ & $0.5507$ & $0.5485$ & $0.2830$ & $0.2655$\\
\midrule
$d=5$, dense 1 & $0.8548$ & $0.1785$ & $0.1785$ & $0.1785$ & $0.1785$ & $0.0901$ & $0.0884$\\
$d=5$, dense 2 & $0.1906$ & $0.5888$ & $0.5901$ & $0.5898$ & $0.5889$ & $0.2964$ & $0.2926$\\
$d=5$, dense 3 & $0.2520$ & $0.1975$ & $0.1973$ & $0.1974$ & $0.1975$ & $0.0989$ & $0.0985$\\
$d=5$, sparse 1 & $0.2068$ & $1.3307$ & $1.3273$ & $1.3282$ & $1.3304$ & $0.6703$ & $0.6601$\\
$d=5$, sparse 2 & $0.8950$ & $0.2024$ & $0.2026$ & $0.2025$ & $0.2025$ & $0.1022$ & $0.1003$\\
$d=5$, sparse 3 & $0.4736$ & $2.3887$ & $2.4183$ & $2.4010$ & $2.3903$ & $1.2025$ & $1.1878$\\
\midrule
$d=10$, dense 1 & $0.2085$ & $0.0447$ & $0.0445$ & $0.0446$ & $0.0447$ & $0.0219$ & $0.0228$\\
$d=10$, dense 2 & $0.4063$ & $0.3178$ & $0.3191$ & $0.3183$ & $0.3179$ & $0.1565$ & $0.1615$\\
$d=10$, dense 3 & $0.4197$ & $0.0770$ & $0.0768$ & $0.0770$ & $0.0770$ & $0.0386$ & $0.0384$\\
$d=10$, sparse 1 & $0.8253$ & $0.7134$ & $0.6981$ & $0.7091$ & $0.7123$ & $0.3654$ & $0.3469$\\
$d=10$, sparse 2 & $0.3296$ & $0.4656$ & $0.4732$ & $0.4686$ & $0.4671$ & $0.2199$ & $0.2472$\\
$d=10$, sparse 3 & $0.5997$ & $1.5941$ & $1.6197$ & $1.6006$ & $1.5959$ & $0.8187$ & $0.7772$\\
\end{longtable}
\endgroup

\FloatBarrier
\section{Proof of the rate constants}\label[appendix]{app:rateproof}

We keep the notation of \cref{thm:rateconst} and its proof: $m=2$, $S_a=aA_1+(1-a)A_2$, the Perron vectors $u_a,v_a$ of $S_a$ with $u_a^\top v_a=1$, $\kappa(a)=u_a^\top(aA_1^2+(1-a)A_2^2)v_a$, and the words $w(p,q)$ defined by \eqref{eq:chr_word}. Every $S_a$ with $0<a<1$ is irreducible.

\begin{lemma}\label{lem:wordexp}
Let $\mathcal I\subset(0,1)$ be a compact interval. There are $T_1,h_3,C>0$ such that, for every $h\in(0,h_3]$ and all integers $1\le p<q$ with $b=p/q\in \mathcal I$ and $T:=qh\le T_1$, the word $w=w(p,q)$, defined by \eqref{eq:chr_word}, satisfies
\[
 \left|\frac{\log\rho(B_w)}T-\eta(b)+\frac h2\kappa(b)\right|\le Ch^2,\qquad
 \left|\frac{\log\rho(R_w)}T-\eta(b)-\frac h2\kappa(b)\right|\le Ch^2 .
\]
\end{lemma}

\begin{proof}
Fix $T_1>0$ and take $h_3K_A\le1/2$. Write $S=S_b$, $\eta=\eta(b)$, $u=u_b$ and $v=v_b$, choosing the Perron vectors continuously with $\one^\top v_a=1$ and $u_a^\top v_a=1$. Their norms are bounded on $\mathcal I$. We use the induced $\ell^1$ norm; constants $C$ may change from line to line and depend only on $A_1,A_2,\mathcal I,T_1$, never on $p,q,h,T$.

Let $\sigma=w_j$ on $[t_{j-1},t_j)$, where $t_j=jh$, and let $\Phi(t,s)$ be its transition matrix. If $W(s)$ is the time spent in mode~$1$ up to $s$, put $\delta(s)=W(s)-bs$ and $F=A_1-A_2$. The prefix counts give $\delta(t_k)=h(\lfloor kb\rfloor-kb)$; since $|\delta'|\le1$, we have $|\delta|\le2h$ and $\delta(0)=\delta(T)=0$. Thus $A_\sigma-S=\Theta'$ almost everywhere, where $\Theta=\delta F=O(h)$.

For $0\le s\le t\le T$, set $\Theta_s(r)=\Theta(r)-\Theta(s)$. Integrating the identity $\frac{\mathrm d}{\mathrm dr}(e^{(t-r)S}\Phi(r,s))=e^{(t-r)S}\Theta_s'(r)\Phi(r,s)$ by parts gives
\[
 \Phi(t,s)-e^{(t-s)S}
 =\Theta_s(t)\Phi(t,s)
 +\int_s^t e^{(t-r)S}\bigl(S\Theta_s(r)-\Theta_s(r)A_{\sigma(r)}\bigr)\Phi(r,s)\,\mathrm dr.
\]
The transition matrices and exponentials have norm at most $e^{K_AT_1}$, while $\Theta_s=O(h)$, so $\|\Phi(t,s)-e^{(t-s)S}\|\le Ch$. Over the full period the boundary term vanishes, giving $\|E_1\|\le ChT$ for $E_1=\Phi(T,0)-e^{TS}$. Projecting this identity and using $\Phi(r,0)v=e^{r\eta}v+O(h)$ yields
\[
 u^\top E_1v
 =e^{T\eta}\int_0^T u^\top\bigl(S\Theta-\Theta A_\sigma\bigr)v\,\mathrm dr+O(h^2T).
\]
The replacement error is $O(h^2T)$ because an $O(h)$ vector error multiplies an $O(h)$ factor under the integral. Now $S\Theta-\Theta A_\sigma=\delta[S,F]-\delta\delta'F^2$, and
\[
 u^\top[S,F]v=\eta u^\top Fv-\eta u^\top Fv=0,\qquad
 \int_0^T\delta\delta'\,\mathrm dr
 =\tfrac12\bigl(\delta(T)^2-\delta(0)^2\bigr)=0.
\]
Hence $u^\top E_1v=O(h^2T)$. Both cancellations hold for every $b\in \mathcal I$, not only at a minimizing mixture.

Take $H_j=B_{w_j}$ for every $j$, or $H_j=R_{w_j}$ for every $j$, and set $G_j=e^{hA_{w_j}}$. Taylor and Neumann expansions (minus for explicit Euler and plus for implicit Euler) give
\[
D_j:=H_j-G_j=\mp h^2A_{w_j}^2/2+O(h^3).
\]
Since $\|G_j\|\le e^{hK_A}$ and $\|H_j\|\le e^{2hK_A}$, all products over times at most $T_1$ are uniformly bounded. Exploiting the telescopic sum yields
\[
 E_2:=H_q\cdots H_1-\Phi(T,0)
 =\sum_{j=1}^qH_q\cdots H_{j+1}D_jG_{j-1}\cdots G_1,
 \qquad \|E_2\|\le Cqh^2=ChT.
\]
The same estimate on each suffix, combined with the exact-flow comparison, yields
\[
 u^\top H_q\cdots H_{j+1}=e^{(T-t_j)\eta}u^\top+O(h),\qquad
 G_{j-1}\cdots G_1v=e^{t_{j-1}\eta}v+O(h).
\]
In the telescoping sum these replacements and the local Taylor remainder contribute $O(h^3)$ per term, hence $O(h^2T)$ in total. As $\sum_jA_{w_j}^2=q(bA_1^2+(1-b)A_2^2)$,
\[
 u^\top E_2v
 =\mp e^{(T-h)\eta}\frac{h^2}{2}\sum_{j=1}^q u^\top A_{w_j}^2v+O(h^2T)
 =\mp e^{T\eta}\frac{hT}{2}\kappa(b)+O(h^2T),
\]
where replacing $e^{-h\eta}$ by $1$ also costs $O(h^2T)$. Consequently
\[
 X:=H_q\cdots H_1=e^{TS}+E,\qquad
 \|E\|\le ChT,\qquad
 u^\top Ev=\mp e^{T\eta}\frac{hT}{2}\kappa(b)+O(h^2T).
\]
The eigenvalue gap of $e^{TS}$ shrinks as $T\downarrow0$. To remove this dependence, compare
\[
 M_{b,T}:=\frac{e^{TS_b}-I}{T},\qquad
 \widehat M_T:=\frac{X-I}{T}=M_{b,T}+\frac ET.
\]
Set $M_{b,0}=S_b$. This is a continuous family of irreducible Metzler matrices on $\mathcal I\times[0,T_1]$, with Perron vectors $u_b,v_b$ and eigenvalue $\theta_{b,T}=(e^{T\eta(b)}-1)/T$, extended by $\theta_{b,0}=\eta(b)$. Its Perron eigenvalue has a uniform real-part gap by irreducibility and compactness. Since $E/T=O(h)$, simple eigenvalue perturbation, after reducing $h_3$, gives uniformly
\[
 \eta(\widehat M_T)=\theta_{b,T}+\frac{u^\top Ev}{T}+O(h^2).
\]
The Euler factors are nonnegative for $h$, so $X\ge0$ and $\eta(\widehat M_T)=(\rho(X)-1)/T$. Thus
\[
 \rho(X)=e^{T\eta}+u^\top Ev+O(h^2T)
 =e^{T\eta}\bigl(1\mp\tfrac{hT}{2}\kappa(b)+O(h^2T)\bigr).
\]
The quadratic remainder in the scalar logarithm is $O(h^2T^2)$. Dividing by $T\le T_1$ gives
\[
 \frac{\log\rho(X)}T=\eta\mp\frac h2\kappa(b)+O(h^2),
\]
uniformly even as $T\downarrow0$, proving both estimates.
\end{proof}

\emph{The frequency average of \cref{rem:freqavg}.}
At the minimizing mixture, set $S=S_{a^\ast}$, $c=\eta(S)$ and
$\Pi=vu^\top$. On the invariant subspace $\ker u^\top$, the matrix
$N=(S-cI)|_{\ker u^\top}$ has spectrum in $\operatorname{Re}z<0$,
and for $\omega\ne0$,
\[
 (c+\imath \omega-S)^{-1}
 =\frac{\Pi}{\imath \omega}+(\imath \omega-N)^{-1}(I-\Pi).
\]
After multiplication by $u^\top F$ and $Fv$, the first term vanishes
because $u^\top Fv=0$. The scalar identity
$\mathrm{PV}\int_{\mathbb{R}}(\imath\omega-y)^{-1}\,\mathrm d\omega=\pi$
for $\operatorname{Re}y<0$ gives
$\mathrm{PV} \int_{\mathbb{R}}(\imath\omega-N)^{-1}\,\mathrm d\omega=\pi I$
on the stable subspace. Real matrices make $\psi$ even, and the
leading resolvent term is purely imaginary, so
$\psi(\omega)=O(\omega^{-2})$. Since the stable resolvent is bounded
near zero, $\psi$ is integrable and
\[
 \int_0^\infty\psi(\omega)\,\mathrm d\omega
 =\frac{\pi}{2}u^\top F(I-\Pi)Fv
 =\frac{\pi}{2}u^\top F^2v.
\]
Multiplying by $2a^\ast(1-a^\ast)/\pi$ gives the formula for $Q$
in \cref{rem:freqavg}.

\smallskip
\emph{Three or more modes.}
For a word of period $T>0$, let $W_i(s)$ be the occupation time of
mode~$i$ up to $s$, and set $b_i=W_i(T)/T$ and $S=\sum_i b_iA_i$.
Assume that $S$ is irreducible, and let $u,v$ be its left and right
Perron vectors, normalized by $u^\top v=1$. Set
$\tilde F_i=A_i-A_m$ and $\delta_i(s)=W_i(s)-b_is$ for $i<m$. Then
$A_\sigma-S=\Theta'$ with
$\Theta=\sum_{i<m}\delta_i\tilde F_i$ and
$\Theta(0)=\Theta(T)=0$. The Perron projection of $[S,\Theta]$
still vanishes, but the remaining term is
\[
 -\int_0^T u^\top\Theta\Theta'v\,\mathrm dr
 =-\frac{1}{2}\sum_{i<k}\mathcal A_{ik}\,
   u^\top[\tilde F_i,\tilde F_k]v,
 \qquad
 \mathcal A_{ik}
 =\int_0^T(\delta_i\delta_k'-\delta_k\delta_i')\,\mathrm dr,
\]
where $1\le i<k<m$. Integration by parts eliminates the diagonal
terms and pairs the others into commutators. Here $\mathcal A_{ik}$
is twice the signed area of the discrepancy path in coordinates
$i,k$. If $\delta_i=O(h)$, boundedness of $\delta_i'$ gives
$\mathcal A_{ik}=O(hT)$. This term depends on letter order and need
not vanish. Thus, for three or more modes, an additional
ordering-dependent contribution must be controlled when identifying
the first-order coefficient.

\FloatBarrier
\section{Proof of the restriction lemma} \label[appendix]{app:faceproof}

\begin{proof}[Proof of \cref{lem:restriction}]
Fix a nonempty $J$ and $\mu\in\mathcal M$. For simplicity, we divide the proof.

\smallskip
\emph{Proof of (a) and (b).} The set $X_J$ need not be closed, so invariance of $\mu$ does not by itself imply invariance of its restriction. The following cutoff estimates show both that no mass flows out of $F_J$ and that the boundary terms vanish when restricting to $X_J$. Set $s(p)=\sum_{k\notin J}p_k$. Then $\nabla s\cdot F_i=\ell_i+r_i$, with the outflow
\[
\ell_i(p)=\sum_{k\notin J,\,j\in J}(A_i)_{kj}p_j\ge0, \qquad |r_i(p)|\le2K_As(p).
\]
Similarly, for $j\in J$, $\nabla p_j\cdot F_i=\iota_{ij}+r_{ij}$ with 
\[
\iota_{ij}(p)=\sum_{l\ne j}(A_i)_{jl}p_l\ge0, \qquad |r_{ij}(p)|\le2K_Ap_j.
\]
For $\varepsilon>0$ consider the $C^1$-regular cut-off function $\psi_\varepsilon(r)=(r/\varepsilon)(2-r/\varepsilon)$ on $[0,\varepsilon]$, $\psi_\varepsilon=1$ on $[\varepsilon,\infty)$ and $\psi_\varepsilon(r)=2r/\varepsilon$ for $r<0$. It satisfies $0\le\psi'_\varepsilon\le2/\varepsilon$ on $[0,\infty)$, $\psi_\varepsilon'(0)=2/\varepsilon$, and $\psi'_\varepsilon(r)\,r\le2(r/\varepsilon)\mathbf 1_{r<\varepsilon}$ for $r\ge0$ as well as $0\le\psi_\varepsilon\le1$ on $[0,\infty)$. By \eqref{eq:invariance} with $f=\psi_\varepsilon\circ s$,
\[
 \int\psi'_\varepsilon(s)\,\ell\, \mathrm d\mu=-\int\psi'_\varepsilon(s)\,r\, \mathrm d\mu\le2K_A\int\psi'_\varepsilon(s)\,s\, \mathrm d\mu\le4K_A\int\frac s\varepsilon\mathbf 1_{s<\varepsilon}\,\mathrm d\mu .
\]
Hence $\int_{\{s=0\}}\ell\, \mathrm d\mu\le2K_A\varepsilon$, so $\ell=0$ $\mu$-almost everywhere on $\{s=0\}$, and $\int\psi'_\varepsilon(s)|\nabla s\cdot F|\,\mathrm d\mu\le8K_A\int\frac s\varepsilon\mathbf 1_{s<\varepsilon}\, \mathrm d\mu\to0$ as $\varepsilon\to0$, by dominated convergence. In the same way $\int\psi'_\varepsilon(p_j)|\nabla p_j\cdot F|\, \mathrm d\mu\to0$ for $j\in J$. This gives (a) and (b):
\begin{itemize}[leftmargin=*]
\item[(a)] For $\mu$-almost every $(p,i)$ with $p\in X_J$ we have $i\in S(J)$. Otherwise $(A_i)_{kj}>0$ for some $k\notin J$, $j\in J$, and $\ell_i(p)\ge(A_i)_{kj}p_j>0$. In particular $\mu(X_J)=0$ if $S(J)=\emptyset$.

\item[(b)] It holds $\int_{X_J}\nabla f\cdot F\,\mathrm d\mu=0$ for every $f\in C^1(\R^d)$. Apply \eqref{eq:invariance} to $f\chi_\varepsilon$ with $\chi_\varepsilon=(1-\psi_\varepsilon(s))\prod_{j\in J}\psi_\varepsilon(p_j)$, which tends to the indicator of $X_J$ pointwise on $\Delta$, with $|\chi_\varepsilon|\le1$. Then $|\nabla\chi_\varepsilon\cdot F_i|\le\psi'_\varepsilon(s)|\nabla s\cdot F_i|+\sum_{j\in J}\psi'_\varepsilon(p_j)|\nabla p_j\cdot F_i|$, whose integral tends to $0$.
\end{itemize}

\emph{Proof of (c).} Let $S(J)\ne\emptyset$. We show $\int_{X_J}\gamma\,\mathrm d\mu\ge\lexp_J\,\mu(X_J)$. If $J=\{j\}$, then $X_J=\{e_j\}$ and $\gamma_i(e_j)=(A_i)_{jj}\ge\lexp_J$ for $i\in S(J)$, so (a) gives the claim. Let $|J|\ge2$ and $\nu'<\lexp_J$, and define, for $x\in F_J$, the function
\[
 v(x)=\inf\bigl\{e^{-\nu't}\one^\top\Phi_\sigma(t)x \ : \ t\ge0,\ \sigma\text{ with values in }S(J)\bigr\}.
\]
Then $v$ is concave and $v\le\one^\top x$. Moreover $v(x)\ge\min_{j\in J}x_j>0$ on $X_J$, because $\one^\top\Phi_\sigma(t)x\ge\min_{j\in J}x_j\,\|\Phi_\sigma(t)|_{F_J}\|_1\ge\min_{j\in J}x_j\,e^{\lexp_Jt}$ by \cref{prop:basic}(i). Finally $v(\Phi_\sigma(t)x)\ge e^{\nu't}v(x)$ for signals with values in $S(J)$. Fix $\theta\in(0,1]$ and let $w=v^\theta$ on $X_J$, a relatively open convex subset of the affine hull of $F_J\cap\Delta$, whose tangent space $T_J=\{y:\ y_k=0\ (k\notin J),\ \one^\top y=0\}$ contains $F_i(p)$ for $p\in X_J$ and $i\in S(J)$. Since $r\mapsto r^\theta$ is concave and nondecreasing, $w$ is concave, hence locally Lipschitz; and $0<w\le1$. On $\{p\in X_J:\ p_j\ge\delta/2\ \forall j\in J\}$ we have $v\ge\delta/2$, so $w$ is Lipschitz there. By homogeneity and \eqref{eq:growth}, along the projective solution of a mode $i\in S(J)$, which stays in $X_J$, we have 
\[
w(p(t))\ge w(p)\exp\left(\theta\int_0^t(\nu'-\gamma_i(p(s)))\, \mathrm ds\right).
\]
Since $p(t)=p+tF_i(p)+O(t^2)$ and $w$ is locally Lipschitz, at every point $p$ where $w$ is differentiable this gives
\begin{equation}\label{eq:supergrad}
 \nabla w(p)\cdot F_i(p)\ge\theta(\nu'-\gamma_i(p))\,w(p)\qquad\text{for every }i\in S(J),
\end{equation}
with $\nabla w$ the gradient in $T_J$. To use restricted invariance, we cut off near the boundary and mollify $w$ to obtain a $C^1$ test function. The preceding estimates will remove the resulting boundary terms.
For $\delta>0$ let $\zeta_\delta\in C^1(\R)$ satisfy $\zeta_\delta=0$ on $(-\infty,\delta]$, $\zeta_\delta=1$ on $[2\delta,\infty)$, $0\le\zeta_\delta\le1$ and $0\le\zeta_\delta'\le2/\delta$, and let $\kappa_\delta(p)=\prod_{j\in J}\zeta_\delta(p_j)$, so $0\le\kappa_\delta\le1$. Mollify $w$ in the affine hull of $X_J$, with $w_\omega=w*\rho_\omega$. For $\omega<\delta/2$ it is smooth near $\{\kappa_\delta>0\}\cap X_J$, $0<w_\omega\le1$ there, and $\nabla w_\omega=(\nabla w)*\rho_\omega$. By \eqref{eq:supergrad} and the Lipschitz continuity of $F_i$, $\gamma_i$ and $w$ on $\{p\in X_J \ :\ p_j\ge\delta/2\ \forall j\in J\}$, there is $c_\delta$, depending also on $\theta$ and $\nu'$, with
\[
 \nabla w_\omega\cdot F_i\ge\theta(\nu'-\gamma_i)\,w_\omega-c_\delta\omega\qquad\text{on }\{\kappa_\delta>0\}\cap X_J,\ i\in S(J).
\]
Let $G=\kappa_\delta w_\omega$ on the affine hull of $F_J\cap\Delta$, extended by $0$ and note that it is $C^1$-regular. Set
\[
 H(p):=\phi(\one^\top p_J)\,G\bigl(p_J/\one^\top p_J\bigr),
\]
when $\one^\top p_J\ne0$, and set $H(p)=0$ when $\one^\top p_J=0$. Here $\phi\in C^1(\R)$ vanishes near $0$ and satisfies $\phi(1)=1$, so $H\in C^1(\R^d)$, and on $X_J$, for $i\in S(J)$, $\nabla H\cdot F_i=\nabla G\cdot F_i$, because $\one^\top(F_i)_J=0$ there. By (b) and (a),
\[
 0=\int_{X_J}\nabla H\cdot F\, \mathrm d\mu\ge\theta\int_{X_J}\kappa_\delta(\nu'-\gamma)\,w_\omega\, \mathrm d\mu-c_\delta\omega\,\mu(X_J)-\sum_{j\in J}\int_{X_J}|\zeta'_\delta(p_j)|\,|\nabla p_j\cdot F|\, \mathrm d\mu .
\]
Since $|\zeta'_\delta|\le8\psi'_{4\delta}$ on $[0,\infty)$, the last sum tends to $0$ as $\delta\to0$, by the estimate above. Letting $\omega\to0$ and then $\delta\to0$ gives $\int_{X_J}(\gamma-\nu')\,v^\theta\,\mathrm d\mu\ge0$. Letting $\theta\to0$, with $v^\theta\to1$ on $X_J$ and dominated convergence, gives $\int_{X_J}\gamma\,\mathrm d\mu\ge\nu'\mu(X_J)$. Now let $\nu'\uparrow\lexp_J$.
\end{proof}

\FloatBarrier
\section{Construction and experimental protocols}\label[appendix]{app:protocols}
This appendix records the procedures behind \cref{sec:gripenberg}. The sample generators and detailed scaling and switching results are in \cref{app:samples}.

\subsection{Notation and verification}
For a Metzler family $\{A_i\}$, $\lexp$ is its lower Lyapunov exponent, $g(h)$ the explicit Euler lower bound and $G(h)$ the implicit Euler upper bound, with $B_i=I+hA_i$, $R_i=(I-hA_i)^{-1}$, $\cB_h=\{B_i\}$ and $\cR_h=\{R_i\}$. Write $S_\alpha=\sum_i\alpha_iA_i$ for a convex mixture. Word estimates are numerical upper estimates of the discrete lower spectral radii, not certified approximations of them. A target gap below such an estimate is not a certified two-sided error for $\lexp$.

Certified bracket endpoints use exact rational checks and independent Python and MATLAB replay. They hold for both exact-decimal and binary64 readings of the input data. The Gripenberg algorithm of Guglielmi and Maiale is discussed in \cref{app:gm}. Its adaptive variant supplies the estimates in \cref{tab:gmruns,tab:gmscaling,tab:gmswitching}; lower estimates require a posteriori certification. A failed or unfinished closure is inconclusive. Timings are wall-clock times from single-threaded runs, several executed concurrently on a 16-core Apple-silicon machine under normal load.

\subsection{Certificate construction and verification}

\paragraph{Lower certificate search} To find a certificate at a target level $t$, we close a seed set under the maps $v\mapsto B_iv/t$. Start with coordinate vectors and add $w=B_iv/t$ when the LP $\max\{\one^\top\mu:\ V\mu\le w,\ \mu\ge0\}<1$ detects $w\notin\operatorname{conv}V+K$. Discard redundant vertices. Termination gives the inequalities in \cref{prop:polytope} in floating point.

Bisect between a certified level and the smaller of $1+hc$ and the adaptive Gripenberg word bound, with $c=\min_\alpha\eta(S_\alpha)$. The starting level is the smallest column sum of $B_i$, or an estimate whose closure finishes. Probes stopped by time or vertex caps are inconclusive.

\paragraph{Exact verification} Rationalized LP weights give an exact level $t_*$; a dual vector excludes invariance above $\tau_{\rm up}$. All $342$ lower certificates have $\tau_{\rm up}/t_*-1<5\times10^{-8}$, measuring sharpness for that polytope, not accuracy for $\lexp$. Checks are independently replayed.

\paragraph{Upper word certificates} We search finite words using the settings in \cref{subsec:searchsettings}. We certify the best word by \cref{prop:upper}: rationalize its Perron vector $x_0$ and $r$, and round $(1+2^{-40})R_{w_j}x_{j-1}$ upward to dyadic $x_j$. The inverse-free check raises the word rate by less than $10^{-6}$.

\paragraph{Mixture and schedule certificates} Mixture bounds use Collatz--Wielandt at a rational $\alpha$ near the numerical minimizer. For optimized two-segment schedules with rational durations, choose a diagonal shift $\gamma$ and enclose $e^{tA_i}$ by the Taylor series of $2^{-j}t(A_i+\gamma I)\ge0$ with a tail bound, square $j$ times and multiply by an upper bound for $e^{-\gamma t}$, rounding upward. Collatz--Wielandt bounds $\rho(\Phi)$, and an artanh series with remainder bounds $\log\rho(\Phi)/T$. To verify a proposed upper growth rate $U$, the independent replay uses uniformization and compares an upper bound for $\rho(\Phi)$ with a lower bound for $e^{UT}$.

\subsection{Test-pair certificate records}
Lower certificates have $169$ to $3131$ vertices and hold for exact-decimal and binary64 data. Case~1440 at $h=0.025$ reuses the $h=0.05$ polytope (\cref{prop:polytope}(ii)). Bisections took $8$ to $4391$s per row; upper certificates use words of length $5$ to $141$. 

For case~1 at $h=0.1$, the certified $-1.044234$ is close to the separate floating-point word estimate $-1.044232$.

For case~1440 at every reported step and case~148 at $h\le0.05$, certified upper bounds beat the numerical mixture minima by up to $2.2\times10^{-3}$ and $3.1\times10^{-3}$, respectively. This comparison does not prove $\lexp<c$ without a global lower bound for $c$.

\subsection{Experimental setting} 

Times refer to Python, except for transposed, extended and larger-sample searches in \cref{obs:scaling}, which used MATLAB. Runs were single-threaded on a $16$-core Apple-silicon machine, several concurrently. The repository supplies the MATLAB package, complete run records and reproduction instructions.

\Cref{app:gm} gives the Gripenberg implementation, settings and fixed/adaptive comparison.

Every adaptive Gripenberg lower estimate on the test pairs is certified a posteriori; seven of the $54$ random-pair estimates remain uncertified within the budget (\cref{app:gm}). Adaptive outputs guide the certificate search, but claimed lower bounds require \cref{prop:polytope}.

\subsection{Algorithm (L) protocol}\label{subsec:algorithmL}
\paragraph{Upper bound} Among all products of the matrices $e^{\tau A_i}$ of length at most $16$, $14$ or $12$ (for $d\le3$, $d=5$ and $d=10$), take the product $\Pi$ of length $n$ that minimizes $\rho(\Pi)^{1/n}$. Its rate $\ln\rho(\Pi)/(n\tau)$ is an upper bound for $\lexp$, computed in floating point.

\paragraph{Polytope} Close the Perron vector of $\Pi$ under the maps $e^{\tau A_i}e^{-\tau(\beta-\nu)}$, where $\beta$ is that rate, and find by bisection the smallest $\nu\ge0$ for which the closure finishes. We use the closure code of \cref{prop:polytope} and the same budgets: $300$\,s per closure and $1800$\,s per step.

\subsection{Algorithm (L) comparison conventions}
Over steps $0.1,0.05$, Euler gives the larger lower bound on $9/15$ pairs and Algorithm~(L) on $6/15$ (\cref{tab:glp}, excluding reducible case~148). The tangent test improves neither the planar nor case~1 Euler polytope. Algorithm~(L)'s numerical upper bound improves on implicit Euler on $14$ pairs; case~1440 uses an Euler word of length $50$, beyond its exhaustive search. Certified mixtures or schedules improve both.

Closure code, budgets and hardware are shared, with twelve concurrent runs, but our additional long-word search prevents attributing differences solely to discretization. Our naive implementation may differ from that of \cite{GLP}; only the Euler lower bounds are independently replayed.

\subsection{Work-precision protocol}\label{subsec:workprecision}
\cref{fig:wp} uses fixed certified upper bounds: schedules for case~1440 and dimension-$5$ pair~4 of \cref{tab:switching}, a mixture for case~1. At $h=0.2,0.1,0.05,0.025$, the explicit-word search supplies a Perron orbit for the levels of the switching sample in \cref{obs:scaling}; Algorithm~(L) uses $\tau=h$ and $1200$\,s of closures. Times include searches, failures and checks, with eight concurrent Python runs. MATLAB reproduces the twelve Euler bounds; budget-limited Algorithm~(L) results vary slightly by machine.

\subsection{Search settings}\label{subsec:searchsettings}
\paragraph{Euler word search} Search all words of length at most $16$, $14$ or $12$ (for $d\le3$, $d=5$, $d=10$), two-block words $1^a2^b$ with $a,b\le\min(4/h,400)$, and balanced words of frequencies $p/q$ with $q\le\min(4/h,320)$, then apply single-letter flips, insertions and deletions.

\paragraph{Scaling closures} Budgets were $300$\,s per probe and $1800$\,s per row, plus $20$--$300$\,s for exact checks, with twelve concurrent runs.

\paragraph{Switching closures} Here $g_{\rm up}$ is the numerical explicit-word estimate and $\varepsilon$ the target tolerance below it. The search minimizes $\rho(B_w)^{1/|w|}$ over two-block words with durations within four steps of the optimized schedule. The best Perron orbit seeds closures at $1+h(g_{\rm up}-\varepsilon)$ for $\varepsilon=10^{-5},10^{-3},10^{-2},10^{-1}$, with respective budgets $600,600,300,300$\,s and at most $20000$ vertices. Transposed searches on the $23$ rows reaching only $0.1$ or no level, with extra targets $0.15,0.2$ and $300$\,s per closure, improve $11$ rows. For dimension-$10$ rows still at $\varepsilon\ge10^{-2}$, extended MATLAB searches use both orientations, five two-block Perron orbits, final vertices from the adaptive Gripenberg algorithm on the transpose at $h=0.1$, and certified $h=0.1$ vertices at $h=0.05$. Targets $0.15,0.1,0.05,0.02,0.01,10^{-3},10^{-5}$ below the current level receive $1800$\,s per closure; improvements are marked $^+$.

\FloatBarrier
\section{Scaling and switching samples}\label[appendix]{app:samples}
This appendix gives the sample construction and detailed results supporting \cref{obs:scaling}; search and verification protocols are in \cref{app:protocols}.

\subsection{Sample construction and success counts}
Transferring certificates improves every direct $h=0.05$ bound for the original six dimension-$5$ pairs.

The initial sample has three dense and three sparse pairs in each of dimensions $5,10$. Dense off-diagonal entries are uniform in $[0,2/(d-1)]$, diagonals in $[-3,1]$. Sparse entries are nonzero with probability $0.5$ ($d=5$) or $0.3$ ($d=10$), then log-normal with $\sigma=1.5$ and expected row sum $2.46$; diagonals are uniform in $[-4,1]$. We require irreducible modes and $\alpha^\ast\in[0.05,0.95]$ ($13$ rejections). No improving two-segment schedule was found; optimal chattering is unproved.

The sample was extended in MATLAB to ten dense and ten sparse pairs per dimension, plus three dense and two sparse pairs in dimension~$20$ (sparse probability $0.15$). All bracket endpoints pass replay. Upper mixture bounds use Collatz--Wielandt, except dimension-$10$ sparse~8, which uses an improving schedule.

At $h=0.1$, lower bounds lie within $4.4\times10^{-3}$ of the numerical word estimates (\cref{tab:scaling,tab:scalingnew}) for $18/20$ dimension-$5$ pairs, but for only $1/20$ dimension-$10$ pairs (sparse~7). Dimension-$20$ certificates use $1845$--$2810$ vertices, with widths $0.33$--$1.6$. Mixture bounds improve on implicit Euler for the original twelve pairs (\cref{tab:glp}).

In these experiments, finer steps increase search cost; bracket widths also include the schedule's suboptimality.

Using the sparse generator above with irreducible modes, we kept pairs whose optimized two-segment schedule beats the numerical value of $c$ by $10^{-3}$, i.e. $16/300$ draws in dimension~$5$ and $28/300$ in dimension~$10$. The first ten of each form \cref{tab:switching}, with best periods $0.6$--$10.9$.

Target tolerance $\varepsilon$ means level $g_{\rm up}-\varepsilon$ below a floating-point word estimate, not a certified two-sided gap. Rounding the lower endpoint down and word value up each adds at most $10^{-6}$. Schedule bounds and polytopes pass exact checks and replay for both data readings; $\lexp<c$ remains a numerical comparison.
The direct, transposed and extended searches, including seed choices, target levels and budgets, are specified in \cref{subsec:searchsettings}. The lower endpoints in \cref{tab:switching} report separate searches at each step; retaining transferred certificates prevents the best available lower bound from deteriorating as the step decreases (\cref{prop:polytope}(ii)).

The sample was extended to $30$ qualifying pairs in each dimension (among $643$ and $343$ draws), plus five in dimension~$20$ (among $76$, sparse probability $0.15$), using the same MATLAB protocol without vertices from Gripenberg runs on new pairs. \cref{tab:switchingnew} completes the switching sample.

At $h=0.1$, target $10^{-3}$ is reached for all $30$ dimension-$5$ pairs and $14/30$ dimension-$10$ pairs, with median bracket widths $1.4\times10^{-2}$ and $3.0\times10^{-2}$. At $h=0.05$, only $5/30$ dimension-$10$ pairs reach $10^{-2}$. Dimension~$20$ closes only at target $0.1$: for all five pairs at $h=0.1$ and for one at $h=0.05$, with widths $0.11$--$0.13$.

\subsection{Detailed results}
Here $^\top$ marks transposition, $^+$ extended search, and $^\dagger$ transfer from a larger step, justified by \cref{prop:polytope}(ii). A dash denotes failure to close within the budget. Upper bounds and widths are rounded upward; lower endpoints downward.

\begingroup
\small
\setlength{\tabcolsep}{4pt}
\begin{longtable}{@{}lcccccc@{}}
\caption{Dimension-$5$ pairs of \cref{obs:scaling}. Here $c$ and $g(h)\le$ are floating-point mixture and adaptive Gripenberg word estimates, rounded up; $^\ast$time limit for the adaptive algorithm. Lower and upper are certified Euler endpoints (\cref{prop:polytope,prop:upper}); width is their difference and $|V|$ the vertex count. $^\dagger$The $h=0.1$ certificate, valid at $h=0.05$ by \cref{prop:polytope}(ii), improves the direct finer-step certificate in every row.}\label{tab:scaling}\\
\toprule
pair ($c$) & $h$ & lower & $g(h)\le$ & upper & width & $|V|$\\
\midrule
\endfirsthead
\caption[]{Continued from the previous page.}\\
\toprule
pair ($c$) & $h$ & lower & $g(h)\le$ & upper & width & $|V|$\\
\midrule
\endhead
\bottomrule
\endfoot
\bottomrule
\endlastfoot
dense 1 ($-0.213383$) & 0.1 & $-0.222751$ & $-0.222636$ & $-0.204815$ & $1.8\times10^{-2}$ & 1755\\*
 & 0.05 & $-0.222751^\dagger$ & $-0.217883^\ast$ & $-0.209005$ & $1.4\times10^{-2}$ & 1755\\
\midrule
dense 2 ($-0.029331$) & 0.1 & $-0.059850$ & $-0.059459$ & $-0.000609$ & $6.0\times10^{-2}$ & 1391\\*
 & 0.05 & $-0.059850^\dagger$ & $-0.044288^\ast$ & $-0.014798$ & $4.6\times10^{-2}$ & 1391\\
\midrule
dense 3 ($+0.232972$) & 0.1 & $\phantom{-}0.222107$ & $\phantom{-}0.223019^\ast$ & $\phantom{-}0.242747$ & $2.1\times10^{-2}$ & 1146\\*
 & 0.05 & $\phantom{-}0.222107^\dagger$ & $\phantom{-}0.228093^\ast$ & $\phantom{-}0.237887$ & $1.6\times10^{-2}$ & 1146\\
\midrule
sparse 1 ($-0.034259$) & 0.1 & $-0.107049$ & $-0.102688$ & $\phantom{-}0.029982$ & $1.4\times10^{-1}$ & 1118\\*
 & 0.05 & $-0.107049^\dagger$ & $-0.067500^\ast$ & $-0.001609$ & $1.1\times10^{-1}$ & 1118\\
\midrule
sparse 2 ($-0.102660$) & 0.1 & $-0.641285$ & $-0.113168^\ast$ & $-0.092940$ & $5.5\times10^{-1}$ & 1802\\*
 & 0.05 & $-0.641285^\dagger$ & $-0.107746^\ast$ & $-0.097694$ & $5.5\times10^{-1}$ & 1802\\
\midrule
sparse 3 ($+0.053208$) & 0.1 & $-0.071028$ & $-0.070350$ & $\phantom{-}0.170806$ & $2.5\times10^{-1}$ & 1454\\*
 & 0.05 & $-0.071028^\dagger$ & $-0.007714^\ast$ & $\phantom{-}0.112334$ & $1.9\times10^{-1}$ & 1454\\
\end{longtable}
\endgroup

\begingroup
\footnotesize
\setlength{\tabcolsep}{3pt}
\begin{longtable}{@{}rlccccccc@{}}
\caption{Additional pairs of \cref{obs:scaling}. Here $c$ is a certified mixture upper bound, or $^s$a certified schedule bound below the numerical mixture minimum. Lower and upper are certified Euler endpoints; $^\dagger$transferred from $h=0.1$, improving the direct certificate. $g(h)\le$ is the floating-point explicit-word estimate, rounded up. Width is $c$ minus the larger lower bound. Upper bounds and widths are rounded up; all endpoints pass both replays.}\label{tab:scalingnew}\\
\toprule
& & & \multicolumn{3}{c}{$h=0.1$} & \multicolumn{2}{c}{$h=0.05$} & \\
\cmidrule(lr){4-6}\cmidrule(lr){7-8}
$d$ & pair & $c$ & lower & $g(h)\le$ & upper & lower & upper & width\\
\midrule
\endfirsthead
\caption[]{Continued from the previous page.}\\
\toprule
& & & \multicolumn{3}{c}{$h=0.1$} & \multicolumn{2}{c}{$h=0.05$} & \\
\cmidrule(lr){4-6}\cmidrule(lr){7-8}
$d$ & pair & $c$ & lower & $g(h)\le$ & upper & lower & upper & width\\
\midrule
\endhead
\bottomrule
\endfoot
\bottomrule
\endlastfoot
5 & dense 4 & $\phantom{-}0.224753$ & $\phantom{-}0.208178$ & $\phantom{-}0.208195$ & $\phantom{-}0.240933$ & $\phantom{-}0.208178^{\dagger}$ & $\phantom{-}0.232908$ & $1.7\times10^{-2}$\\
5 & dense 5 & $-0.102452$ & $-0.109222$ & $-0.109221$ & $-0.096018$ & $-0.109222^{\dagger}$ & $-0.099172$ & $6.8\times10^{-3}$\\
5 & dense 6 & $\phantom{-}0.590670$ & $\phantom{-}0.561428$ & $\phantom{-}0.561429$ & $\phantom{-}0.620237$ & $\phantom{-}0.561428^{\dagger}$ & $\phantom{-}0.605415$ & $3.0\times10^{-2}$\\
5 & dense 7 & $\phantom{-}0.553355$ & $\phantom{-}0.536324$ & $\phantom{-}0.536325$ & $\phantom{-}0.570945$ & $\phantom{-}0.543527$ & $\phantom{-}0.562075$ & $9.9\times10^{-3}$\\
5 & dense 8 & $\phantom{-}0.956855$ & $\phantom{-}0.940609$ & $\phantom{-}0.940616$ & $\phantom{-}0.974893$ & $\phantom{-}0.940609^{\dagger}$ & $\phantom{-}0.965630$ & $1.7\times10^{-2}$\\
5 & dense 9 & $\phantom{-}0.338856$ & $\phantom{-}0.290975$ & $\phantom{-}0.290976$ & $\phantom{-}0.386052$ & $\phantom{-}0.307913$ & $\phantom{-}0.362520$ & $3.1\times10^{-2}$\\
5 & dense 10 & $\phantom{-}0.491873$ & $\phantom{-}0.476472$ & $\phantom{-}0.476547$ & $\phantom{-}0.508426$ & $\phantom{-}0.476472^{\dagger}$ & $\phantom{-}0.500015$ & $1.6\times10^{-2}$\\
5 & sparse 4 & $-0.267567$ & $-0.405495$ & $-0.405494$ & $-0.151325$ & $-0.405495^{\dagger}$ & $-0.207181$ & $1.4\times10^{-1}$\\
5 & sparse 5 & $-0.111087$ & $-0.279163$ & $-0.279162$ & $\phantom{-}0.029473$ & $-0.193416$ & $-0.037738$ & $8.3\times10^{-2}$\\
5 & sparse 6 & $\phantom{-}0.236626$ & $\phantom{-}0.225171$ & $\phantom{-}0.225360$ & $\phantom{-}0.246479$ & $\phantom{-}0.225171^{\dagger}$ & $\phantom{-}0.241436$ & $1.2\times10^{-2}$\\
5 & sparse 7 & $\phantom{-}2.858514$ & $\phantom{-}2.544951$ & $\phantom{-}2.544952$ & $\phantom{-}3.136416$ & $\phantom{-}2.703311$ & $\phantom{-}3.006918$ & $1.6\times10^{-1}$\\
5 & sparse 8 & $-0.389270$ & $-0.474814$ & $-0.415082$ & $-0.368156$ & $-0.474814^{\dagger}$ & $-0.378216$ & $8.6\times10^{-2}$\\
5 & sparse 9 & $\phantom{-}1.346248$ & $\phantom{-}1.189855$ & $\phantom{-}1.189856$ & $\phantom{-}1.452482$ & $\phantom{-}1.275491$ & $\phantom{-}1.403754$ & $7.1\times10^{-2}$\\
5 & sparse 10 & $-0.475098$ & $-0.571418$ & $-0.569407$ & $-0.388128$ & $-0.571418^{\dagger}$ & $-0.430868$ & $9.7\times10^{-2}$\\
\midrule
10 & dense 4 & $\phantom{-}0.396993$ & $\phantom{-}0.067175$ & $\phantom{-}0.393396$ & $\phantom{-}0.401010$ & $\phantom{-}0.067175^{\dagger}$ & $\phantom{-}0.398963$ & $3.3\times10^{-1}$\\
10 & dense 5 & $\phantom{-}0.592625$ & $\phantom{-}0.535052$ & $\phantom{-}0.579588$ & $\phantom{-}0.606755$ & $\phantom{-}0.535052^{\dagger}$ & $\phantom{-}0.599577$ & $5.8\times10^{-2}$\\
10 & dense 6 & $-0.014617$ & $-0.199153$ & $-0.059777$ & $\phantom{-}0.028727$ & $-0.199153^{\dagger}$ & $\phantom{-}0.007246$ & $1.9\times10^{-1}$\\
10 & dense 7 & $\phantom{-}0.769243$ & $\phantom{-}0.594096$ & $\phantom{-}0.759850$ & $\phantom{-}0.779437$ & $\phantom{-}0.594096^{\dagger}$ & $\phantom{-}0.774244$ & $1.8\times10^{-1}$\\
10 & dense 8 & $\phantom{-}1.052014$ & $\phantom{-}0.732335$ & $\phantom{-}1.045460$ & $\phantom{-}1.059949$ & $\phantom{-}0.732335^{\dagger}$ & $\phantom{-}1.055807$ & $3.2\times10^{-1}$\\
10 & dense 9 & $\phantom{-}0.004987$ & $-0.154668$ & $-0.019725$ & $\phantom{-}0.029610$ & $-0.154668^{\dagger}$ & $\phantom{-}0.017302$ & $1.6\times10^{-1}$\\
10 & dense 10 & $-0.184297$ & $-0.341127$ & $-0.193572$ & $-0.175481$ & $-0.341127^{\dagger}$ & $-0.179831$ & $1.6\times10^{-1}$\\
10 & sparse 4 & $\phantom{-}0.143602$ & $\phantom{-}0.029023$ & $\phantom{-}0.116893$ & $\phantom{-}0.176427$ & $\phantom{-}0.029023^{\dagger}$ & $\phantom{-}0.159333$ & $1.2\times10^{-1}$\\
10 & sparse 5 & $-0.583787$ & $-0.613114$ & $-0.591249$ & $-0.575413$ & $-0.613114^{\dagger}$ & $-0.579692$ & $3.0\times10^{-2}$\\
10 & sparse 6 & $-0.179129$ & $-0.749589$ & $-0.205900$ & $-0.152209$ & $-0.749589^{\dagger}$ & $-0.165618$ & $5.8\times10^{-1}$\\
10 & sparse 7 & $\phantom{-}0.084506$ & $\phantom{-}0.071806$ & $\phantom{-}0.073599$ & $\phantom{-}0.097158$ & $\phantom{-}0.073149$ & $\phantom{-}0.090545$ & $1.2\times10^{-2}$\\
10 & sparse 8 & $\phantom{-}0.973159^{s}$ & $\phantom{-}0.755608$ & $\phantom{-}0.957096$ & $\phantom{-}0.989161$ & $\phantom{-}0.755608^{\dagger}$ & $\phantom{-}0.981195$ & $2.2\times10^{-1}$\\
10 & sparse 9 & $\phantom{-}0.628628$ & $\phantom{-}0.344828$ & $\phantom{-}0.569675$ & $\phantom{-}0.687165$ & $\phantom{-}0.344828^{\dagger}$ & $\phantom{-}0.658158$ & $2.9\times10^{-1}$\\
10 & sparse 10 & $\phantom{-}1.002926$ & $\phantom{-}0.771833$ & $\phantom{-}0.944048$ & $\phantom{-}1.059675$ & $\phantom{-}0.771833^{\dagger}$ & $\phantom{-}1.031721$ & $2.4\times10^{-1}$\\
\midrule
20 & dense 1 & $\phantom{-}0.538291$ & $\phantom{-}0.210786$ & $\phantom{-}0.514021$ & $\phantom{-}0.564454$ & $\phantom{-}0.210786^{\dagger}$ & $\phantom{-}0.551118$ & $3.3\times10^{-1}$\\
20 & dense 2 & $\phantom{-}0.519703$ & $-0.142266$ & $\phantom{-}0.514174$ & $\phantom{-}0.525795$ & $-0.142266^{\dagger}$ & $\phantom{-}0.522678$ & $6.7\times10^{-1}$\\
20 & dense 3 & $\phantom{-}0.690889$ & $\phantom{-}0.117393$ & $\phantom{-}0.685092$ & $\phantom{-}0.697866$ & $\phantom{-}0.117393^{\dagger}$ & $\phantom{-}0.694230$ & $5.8\times10^{-1}$\\
20 & sparse 1 & $\phantom{-}0.557884$ & $-0.306752$ & $\phantom{-}0.524683$ & $\phantom{-}0.591857$ & $-0.306752^{\dagger}$ & $\phantom{-}0.574671$ & $8.7\times10^{-1}$\\
20 & sparse 2 & $\phantom{-}0.928956$ & $-0.633672$ & $\phantom{-}0.884728$ & $\phantom{-}0.979816$ & $-0.633672^{\dagger}$ & $\phantom{-}0.953409$ & $1.6$\\
\end{longtable}
\endgroup

\begingroup
\footnotesize
\setlength{\tabcolsep}{3pt}
\begin{longtable}{@{}rrccccccc@{}}
\caption{Further pairs of the switching sample of \cref{obs:scaling}: $11$--$30$ in dimensions $5$ and $10$, and five in dimension~$20$. Upper is a certified schedule bound, and $g(h)\ge$ the better certified lower bound from both orientations. Width is upper minus the larger lower endpoint. $^\top$Transposed family; $^+$extended search; --- no closure. MATLAB bracket endpoints hold for both data readings and pass both replays; $c$ and explicit-word upper estimates are floating point.}\label{tab:switchingnew}\\
\toprule
& & & & \multicolumn{2}{c}{$h=0.1$} & \multicolumn{2}{c}{$h=0.05$} & \\
\cmidrule(lr){5-6}\cmidrule(lr){7-8}
$d$ & pair & $c$ & upper & $g(h)\le$ & $g(h)\ge$ & $g(h)\le$ & $g(h)\ge$ & width\\
\midrule
\endfirsthead
\caption[]{Continued from the previous page.}\\
\toprule
& & & & \multicolumn{2}{c}{$h=0.1$} & \multicolumn{2}{c}{$h=0.05$} & \\
\cmidrule(lr){5-6}\cmidrule(lr){7-8}
$d$ & pair & $c$ & upper & $g(h)\le$ & $g(h)\ge$ & $g(h)\le$ & $g(h)\ge$ & width\\
\midrule
\endhead
\bottomrule
\endfoot
\bottomrule
\endlastfoot
5 & 11 & $\phantom{-}0.564144$ & $\phantom{-}0.562844$ & $\phantom{-}0.544576$ & $\phantom{-}0.544565$ & $\phantom{-}0.554479$ & $\phantom{-}0.544478$ & $1.9\times10^{-2}$\\
5 & 12 & $\phantom{-}0.077634$ & $\phantom{-}0.073309$ & $\phantom{-}0.032674$ & $\phantom{-}0.032663$ & $\phantom{-}0.054067$ & $\phantom{-}0.044066^{\top}$ & $3.0\times10^{-2}$\\
5 & 13 & $\phantom{-}0.497816$ & $\phantom{-}0.483853$ & $\phantom{-}0.356192$ & $\phantom{-}0.356181$ & $\phantom{-}0.421788$ & $\phantom{-}0.420787$ & $6.4\times10^{-2}$\\
5 & 14 & $-0.815606$ & $-0.825431$ & $-0.864856$ & $-0.864867$ & $-0.843581$ & $-0.843592$ & $1.9\times10^{-2}$\\
5 & 15 & $\phantom{-}0.115795$ & $\phantom{-}0.049136$ & $\phantom{-}0.018823$ & $\phantom{-}0.018812$ & $\phantom{-}0.034435$ & $\phantom{-}0.034424$ & $1.5\times10^{-2}$\\
5 & 16 & $\phantom{-}0.436569$ & $\phantom{-}0.399343$ & $\phantom{-}0.362758$ & $\phantom{-}0.362747$ & $\phantom{-}0.381893$ & $\phantom{-}0.381882$ & $1.8\times10^{-2}$\\
5 & 17 & $\phantom{-}0.774866$ & $\phantom{-}0.748884$ & $\phantom{-}0.701671$ & $\phantom{-}0.701660$ & $\phantom{-}0.726561$ & $\phantom{-}0.725560^{\top}$ & $2.4\times10^{-2}$\\
5 & 18 & $\phantom{-}0.448169$ & $\phantom{-}0.366762$ & $\phantom{-}0.350265$ & $\phantom{-}0.350254$ & $\phantom{-}0.358966$ & $\phantom{-}0.358955$ & $7.9\times10^{-3}$\\
5 & 19 & $\phantom{-}1.565797$ & $\phantom{-}1.533944$ & $\phantom{-}1.523421$ & $\phantom{-}1.523410$ & $\phantom{-}1.528635$ & $\phantom{-}1.528624$ & $5.4\times10^{-3}$\\
5 & 20 & $\phantom{-}0.715940$ & $\phantom{-}0.565468$ & $\phantom{-}0.555534$ & $\phantom{-}0.555523$ & $\phantom{-}0.560635$ & $\phantom{-}0.560624$ & $4.9\times10^{-3}$\\
5 & 21 & $\phantom{-}0.234460$ & $\phantom{-}0.176968$ & $\phantom{-}0.167223$ & $\phantom{-}0.167212$ & $\phantom{-}0.172150$ & $\phantom{-}0.172139^{\top}$ & $4.9\times10^{-3}$\\
5 & 22 & $\phantom{-}0.711039$ & $\phantom{-}0.709367$ & $\phantom{-}0.703899$ & $\phantom{-}0.703888$ & $\phantom{-}0.706719$ & $\phantom{-}0.706708$ & $2.7\times10^{-3}$\\
5 & 23 & $-0.243249$ & $-0.250637$ & $-0.280925$ & $-0.280936$ & $-0.264833$ & $-0.265834$ & $1.6\times10^{-2}$\\
5 & 24 & $\phantom{-}0.713176$ & $\phantom{-}0.685533$ & $\phantom{-}0.659217$ & $\phantom{-}0.658216$ & $\phantom{-}0.672270$ & $\phantom{-}0.672259$ & $1.4\times10^{-2}$\\
5 & 25 & $\phantom{-}1.271640$ & $\phantom{-}1.226099$ & $\phantom{-}1.197536$ & $\phantom{-}1.197525$ & $\phantom{-}1.212215$ & $\phantom{-}1.212204$ & $1.4\times10^{-2}$\\
5 & 26 & $\phantom{-}0.092427$ & $\phantom{-}0.088170$ & $\phantom{-}0.084819$ & $\phantom{-}0.084808$ & $\phantom{-}0.086573$ & $\phantom{-}0.086562$ & $1.7\times10^{-3}$\\
5 & 27 & $\phantom{-}1.043939$ & $\phantom{-}1.034476$ & $\phantom{-}1.000730$ & $\phantom{-}1.000719$ & $\phantom{-}1.018034$ & $\phantom{-}1.018023$ & $1.7\times10^{-2}$\\
5 & 28 & $\phantom{-}0.926204$ & $\phantom{-}0.917748$ & $\phantom{-}0.906064$ & $\phantom{-}0.905063$ & $\phantom{-}0.912019$ & $\phantom{-}0.911018$ & $6.8\times10^{-3}$\\
5 & 29 & $\phantom{-}0.182585$ & $\phantom{-}0.133088$ & $\phantom{-}0.117201$ & $\phantom{-}0.117190$ & $\phantom{-}0.125361$ & $\phantom{-}0.125350$ & $7.8\times10^{-3}$\\
5 & 30 & $-0.224996$ & $-0.332651$ & $-0.424698$ & $-0.424709$ & $-0.375756$ & $-0.376757$ & $4.5\times10^{-2}$\\
\midrule
10 & 11 & $\phantom{-}0.795549$ & $\phantom{-}0.784526$ & $\phantom{-}0.770406$ & $\phantom{-}0.750405^{+}$ & $\phantom{-}0.777396$ & --- & $3.5\times10^{-2}$\\
10 & 12 & $\phantom{-}0.130308$ & $\phantom{-}0.097810$ & $\phantom{-}0.027943$ & $\phantom{-}0.027932$ & $\phantom{-}0.063852$ & $\phantom{-}0.043851^{+}$ & $5.4\times10^{-2}$\\
10 & 13 & $\phantom{-}0.532599$ & $\phantom{-}0.508683$ & $\phantom{-}0.499747$ & $\phantom{-}0.489746^{\top+}$ & $\phantom{-}0.504189$ & $\phantom{-}0.404188^{\top}$ & $1.9\times10^{-2}$\\
10 & 14 & $\phantom{-}0.230392$ & $\phantom{-}0.213210$ & $\phantom{-}0.200303$ & $\phantom{-}0.190302^{\top}$ & $\phantom{-}0.207386$ & $\phantom{-}0.157385^{+}$ & $2.3\times10^{-2}$\\
10 & 15 & $\phantom{-}0.572938$ & $\phantom{-}0.571814$ & $\phantom{-}0.563308$ & $\phantom{-}0.562307^{+}$ & $\phantom{-}0.567561$ & $\phantom{-}0.467560$ & $9.6\times10^{-3}$\\
10 & 16 & $\phantom{-}0.940565$ & $\phantom{-}0.938854$ & $\phantom{-}0.934184$ & $\phantom{-}0.933183^{\top}$ & $\phantom{-}0.936637$ & $\phantom{-}0.926636^{\top}$ & $5.7\times10^{-3}$\\
10 & 17 & $\phantom{-}0.840727$ & $\phantom{-}0.834497$ & $\phantom{-}0.808640$ & $\phantom{-}0.808629$ & $\phantom{-}0.822385$ & $\phantom{-}0.812384$ & $2.3\times10^{-2}$\\
10 & 18 & $\phantom{-}0.063995$ & $\phantom{-}0.062223$ & $\phantom{-}0.045147$ & $-0.054854^{+}$ & $\phantom{-}0.053804$ & --- & $1.2\times10^{-1}$\\
10 & 19 & $\phantom{-}0.559509$ & $\phantom{-}0.556613$ & $\phantom{-}0.551435$ & $\phantom{-}0.550434$ & $\phantom{-}0.553980$ & $\phantom{-}0.533979^{\top+}$ & $6.2\times10^{-3}$\\
10 & 20 & $\phantom{-}0.237440$ & $\phantom{-}0.214496$ & $\phantom{-}0.176835$ & $\phantom{-}0.176824^{\top}$ & $\phantom{-}0.195603$ & $\phantom{-}0.185602^{\top+}$ & $2.9\times10^{-2}$\\
10 & 21 & $\phantom{-}0.836483$ & $\phantom{-}0.829528$ & $\phantom{-}0.815532$ & $\phantom{-}0.805531$ & $\phantom{-}0.822184$ & $\phantom{-}0.722183$ & $2.4\times10^{-2}$\\
10 & 22 & $\phantom{-}0.820069$ & $\phantom{-}0.811741$ & $\phantom{-}0.798140$ & $\phantom{-}0.788139^{\top+}$ & $\phantom{-}0.805005$ & $\phantom{-}0.705004^{\top}$ & $2.4\times10^{-2}$\\
10 & 23 & $\phantom{-}0.802194$ & $\phantom{-}0.753531$ & $\phantom{-}0.748609$ & $\phantom{-}0.748598$ & $\phantom{-}0.751044$ & $\phantom{-}0.651043$ & $5.0\times10^{-3}$\\
10 & 24 & $\phantom{-}0.979140$ & $\phantom{-}0.966475$ & $\phantom{-}0.964792$ & $\phantom{-}0.914791^{\top+}$ & $\phantom{-}0.965613$ & $\phantom{-}0.765612^{\top}$ & $5.2\times10^{-2}$\\
10 & 25 & $-0.128046$ & $-0.133086$ & $-0.159266$ & $-0.160267^{\top}$ & $-0.145737$ & $-0.155738^{\top}$ & $2.3\times10^{-2}$\\
10 & 26 & $\phantom{-}1.427628$ & $\phantom{-}1.330134$ & $\phantom{-}1.304476$ & $\phantom{-}1.294475^{+}$ & $\phantom{-}1.317307$ & $\phantom{-}1.297306^{\top+}$ & $3.3\times10^{-2}$\\
10 & 27 & $\phantom{-}0.545386$ & $\phantom{-}0.538185$ & $\phantom{-}0.408598$ & $\phantom{-}0.407597$ & $\phantom{-}0.481903$ & $\phantom{-}0.471902^{\top+}$ & $6.7\times10^{-2}$\\
10 & 28 & $\phantom{-}0.334115$ & $\phantom{-}0.271371$ & $\phantom{-}0.259274$ & $\phantom{-}0.258273^{\top}$ & $\phantom{-}0.265374$ & $\phantom{-}0.165373^{\top}$ & $1.4\times10^{-2}$\\
10 & 29 & $\phantom{-}0.863320$ & $\phantom{-}0.694745$ & $\phantom{-}0.675212$ & $\phantom{-}0.675201$ & $\phantom{-}0.685172$ & $\phantom{-}0.585171$ & $2.0\times10^{-2}$\\
10 & 30 & $\phantom{-}1.137295$ & $\phantom{-}1.135605$ & $\phantom{-}1.112501$ & $\phantom{-}1.062500^{+}$ & $\phantom{-}1.124068$ & $\phantom{-}0.974067^{\top}$ & $7.4\times10^{-2}$\\
\midrule
20 & 1 & $\phantom{-}1.173320$ & $\phantom{-}1.105091$ & $\phantom{-}1.096882$ & $\phantom{-}0.996881$ & $\phantom{-}1.101188$ & --- & $1.1\times10^{-1}$\\
20 & 2 & $\phantom{-}0.907818$ & $\phantom{-}0.906359$ & $\phantom{-}0.891846$ & $\phantom{-}0.791845$ & $\phantom{-}0.899688$ & --- & $1.2\times10^{-1}$\\
20 & 3 & $\phantom{-}0.718666$ & $\phantom{-}0.635962$ & $\phantom{-}0.604406$ & $\phantom{-}0.504405$ & $\phantom{-}0.620409$ & $\phantom{-}0.520408$ & $1.2\times10^{-1}$\\
20 & 4 & $\phantom{-}0.689138$ & $\phantom{-}0.670308$ & $\phantom{-}0.647472$ & $\phantom{-}0.547471$ & $\phantom{-}0.658777$ & --- & $1.3\times10^{-1}$\\
20 & 5 & $\phantom{-}0.708303$ & $\phantom{-}0.704794$ & $\phantom{-}0.691886$ & $\phantom{-}0.591885$ & $\phantom{-}0.698495$ & --- & $1.2\times10^{-1}$\\
\end{longtable}
\endgroup

\FloatBarrier
\bibliographystyle{sn-basic}
\bibliography{references}

\end{document}